\documentclass[trsc,sglanonrev,letterpaper,onecolumn]{informs4}
\OneAndAHalfSpacedXI
\renewcommand{\Dskips}[4]{\abovedisplayskip=4pt plus 2pt minus 1pt\abovedisplayshortskip=0pt plus 1pt\belowdisplayshortskip=2pt plus 1pt minus 1pt\belowdisplayskip=4pt plus 2pt minus 1pt}
\normalsize
\usepackage[letterpaper,margin=1in]{geometry}
\usepackage{fix-cm,anyfontsize}
\usepackage{amsmath,amssymb}
\usepackage{empheq}
\usepackage{float}
\usepackage{algorithm,algpseudocode}
\usepackage{tikz}
\usepackage{booktabs,multirow}
\renewcommand{\FigureCaptionFontStyle}{\sffamily\bfseries\fontsize{9}{11}\selectfont}
\renewcommand{\FigureNameFontStyle}{\sffamily\bfseries\fontsize{9}{11}\selectfont}
\renewcommand{\TableCaptionFontStyle}{\sffamily\bfseries\fontsize{9}{11}\selectfont}
\renewcommand{\TableNameFontStyle}{\sffamily\bfseries\fontsize{9}{11}\selectfont}

\usepackage{natbib}
\bibpunct[, ]{(}{)}{,}{a}{}{,}
\def\bibfont{\normalsize}

\usepackage{xcolor}

\usepackage[colorlinks=true,allcolors=blue,hypertexnames=false]{hyperref}

\TheoremsNumberedThrough

\begin{document}

\RUNAUTHOR{Xu, You, and Mo}
\RUNTITLE{DG-CGM for Integrated Railroad Blocking}
\ECRUNAUTHOR{}
\ECAUpunct{}

\TITLE{A Dual-Guided Column-Generation-Based Matheuristic for Integrated Railroad Blocking}

\ARTICLEAUTHORS{%
  \AUTHOR{Yida Xu\textsuperscript{\ensuremath{\dagger}}, Zhengzhong Ricky You\textsuperscript{\ensuremath{\dagger}}, Baichuan Mo\textsuperscript{*}}%
  \AFF{Department of Civil Engineering, Tsinghua University, Beijing 100084, China.\\
    \textsuperscript{\ensuremath{\dagger}}These authors contributed equally to this work.\\
    \textsuperscript{*}Corresponding author.}%
}

\ABSTRACT{
Integrated railroad blocking is a central planning problem in freight railway operations, which requires coordinated decisions on block design, blocking sequences, and block routing. The interactions among these decisions and the scale of real railway networks make direct mixed-integer optimization computationally challenging. 
Existing construction and local-search heuristics may overlook promising transportation paths, while expanding the path set increases the size of the resulting optimization problems. 
To solve this problem, we propose a \textit{dual-guided column-generation-based matheuristic} (DG-CGM) that uses linear programming (LP) dual information to identify promising transportation paths and incorporates them into restricted mixed-integer programs to improve feasible solutions.
Starting from an initial plan built by \textit{construction}, DG-CGM alternates \textit{destroy} and \textit{local repair} with \textit{path generation}, which uses dual prices from restricted LPs of a Dantzig--Wolfe reformulation to identify promising transportation paths.
When local improvements continue to stagnate, we further apply \textit{joint repair} to combine alternatives across demand groups.
To address the coupling between logical and physical path pricing, we introduce a \textit{joint reduced cost} that eliminates path-specific dual variables from pricing, enabling more efficient column identification through exact labeling algorithms. 
We evaluate the performance of DG-CGM on both generated instances and real-world cases.
On 150 generated instances with 50--150 nodes, DG-CGM achieves the optimal cost on 89 of the 91 instances solved to proven optimality by Gurobi, with a geometric mean runtime of only 5.2\% of Gurobi's and a mean relative cost reduction of 3.4\%.
On six full-scale scenarios derived from the 2026 INFORMS Railway Applications Section competition data, each with 47,193 nodes, DG-CGM obtains feasible solutions within two hours that serve 95.9\%--99.6\% of the requested freight volume, providing practical blocking plans for large-scale railway operations. The code, instances, and computational results are publicly available at \url{https://github.com/Zhengzhong-You/IRBP}.
}

\maketitle

\section{Introduction}
\label{sec:Introduction}

Railroad blocking consolidates compatible railcars into blocks that move together between classification yards. At Union Pacific's Bailey Yard alone, approximately 10,000 railcars are handled daily, with more than 3,000 sorted across 114 classification tracks \citep{union2025senior,visit2020north}. Improvements in blocking decisions can yield substantial economic benefits. For example, CSX reported that its Cumberland network optimization generated approximately \$15 million in direct savings and reduced railcar handlings by 59,000 \citep{csx2024csx}. At the network level, freight transportation requires three interdependent decisions on which yard-to-yard blocks to open (\textit{block design}), which sequence of compatible blocks each demand should use (\textit{blocking sequence}), and which physical route each open block should follow (\textit{block routing}). 
These three decisions define the integrated railroad blocking problem (IRBP) \citep{ahuja2007solving,informs20262026}. The IRBP seeks to minimize transportation and operating costs while satisfying minimum block loads and shared resource capacities, which couple the decisions across different demands.

Existing studies have developed arc-based mixed-integer programming (MIP) formulations and exact methods for railroad blocking \citep{newton1998constructing,arsalani2023two}. However, real-world instances can involve millions of decision variables. For example, the 2026 INFORMS Railway Applications Section (RAS) Problem Solving Competition provides a network with 47,193 physical nodes, 106,570 directed arcs, and more than 2,000 freight demands \citep{informs20262026}. At this scale, obtaining proven optimal solutions within practical time limits is challenging. Heuristic approaches often restrict candidate blocks or paths to reduce the model size, potentially excluding better solutions \citep{newton1998constructing,fang2023two}. Neighborhood search provides greater flexibility by repeatedly optimizing selected subsets of decisions \citep{ahuja2007solving,yaghini2021mathematical}. Nevertheless, its effectiveness depends on the available paths. Even an optimally solved neighborhood may miss improvements if suitable paths have not been generated, while adding too many paths increases the size of the resulting MIPs.

Linear programming (LP) dual prices and reduced costs identify promising decisions from network-wide information \citep{buson2014reduced,hewitt2013branch}. For the IRBP, such information can be obtained from the LP relaxation of a compact arc-based formulation (ARC). However, this relaxation is weak, which can limit the value of its dual information for guiding integer search. To obtain stronger LP-based guidance, we derive an exact Dantzig--Wolfe (DW) formulation by convexifying the demand and block substructures \citep{lubbecke2005selected,lubbecke2018strength}. We establish its integer equivalence to ARC and prove that its LP relaxation is at least as strong as that of ARC, with strict improvement possible. As shown in Section~\ref{subsec:genratedInstanceExperiments}, the mean lower-bound gap decreases from 67.5\% for ARC to 3.1\% for DW. Although the DW relaxation provides substantially tighter lower bounds, solving it involves challenging pricing subproblems because logical paths and physical block routes are coupled through the DW master constraints. To address this difficulty, we introduce a \textit{joint reduced cost} (JRC) measure that eliminates the dependence on path-specific dual variables in pricing while preserving the exact pricing optimality condition. This allows the corresponding labeling algorithms to identify improving columns more efficiently. In particular, a nonnegative minimum JRC remains a valid stopping criterion for exact pricing, allowing the converged DW relaxation to provide a valid lower bound.

Building on the DW formulation, we develop a \textit{dual-guided column-generation-based matheuristic} (DG-CGM). Its central idea is to construct a feasible blocking plan and progressively improve it by using DW dual information to generate promising transportation paths for restricted MIP optimization. The method has two main stages. In the first stage, \textit{construction} builds an initial feasible solution by grouping compatible demands around candidate blocks and consolidation hubs. The second stage combines \textit{destroy}, \textit{local repair}, \textit{path generation}, and \textit{joint repair}. Specifically, \textit{destroy} releases related demand assignments, and \textit{local repair} jointly reoptimizes their logical paths, served volumes, and associated physical routes. When this local search stagnates, \textit{path generation} solves a restricted DW LP over the current logical and physical paths and uses its dual prices to identify additional logical paths for repair. If subsequent rounds of local search do not improve the solution, \textit{joint repair} further optimizes selected demand groups jointly over their accumulated paths before the search returns to \textit{destroy} and \textit{local repair}.

Our contributions are threefold.
\begin{itemize}

\item We develop the ARC and DW formulations for the IRBP, prove their integer equivalence, and establish the strength of the DW relaxation. We further develop exact JRC-based pricing algorithms that simplify the search for improving columns while preserving pricing optimality. These exact methods also provide benchmarks for evaluating DG-CGM.

\item We develop DG-CGM, which constructs an initial feasible blocking plan and iteratively improves it through dual-guided \textit{path generation} and restricted MIP-based neighborhood search. \textit{Path generation} solves a restricted DW LP and returns the logical paths with the lowest JRC as additional candidates for subsequent \textit{local repair} and \textit{joint repair}. This iterative combination of column generation and integer repair allows the search to explore new blocking decisions while keeping the optimization problems manageable.

\item We conduct extensive computational experiments to evaluate the effectiveness of DG-CGM on generated instances calibrated to real railway data and full-scale RAS instances. On 150 generated instances, DG-CGM matches 89 of the 91 optima proven by ARC MIP within two hours, while reducing feasible costs by 3.4\% on average and up to 40.4\%. On larger instances, its mean gap to the DW lower bound after 300 seconds is at most 1.6\%, compared with 24.2\%--78.8\% for ARC MIP after two hours. We further conduct ablation studies to assess the contributions of the individual search components and validate their design choices. On six full-scale RAS configurations, DG-CGM obtains feasible solutions within two hours that serve 95.9\%--99.6\% of the requested freight volume, demonstrating its ability to construct high-service blocking plans on large railway networks.
 
\end{itemize}

The remainder of this paper is organized as follows. Section~\ref{sec:LiteratureReview} reviews related work. Section~\ref{sec:ProblemDefinition} defines the IRBP and presents ARC. Section~\ref{sec:DWReformulation} develops the DW formulation and pricing algorithms. Section~\ref{sec:PrimalHeuristic} presents DG-CGM, followed by computational results in Section~\ref{sec:HeuristicExperiments}. Section~\ref{sec:Conclusion} concludes the paper.

\section{Related Work}
\label{sec:LiteratureReview}

We review railroad blocking and integrated service network design in Section~\ref{subsec:LiteratureBlocking}, DW decomposition and column generation in Section~\ref{subsec:LiteraturePricing}, and heuristic and matheuristic methods in Section~\ref{subsec:LiteratureMatheuristics}.

\subsection{Railroad Blocking and Integrated Service Network Design}
\label{subsec:LiteratureBlocking}

Early railroad blocking studies focused on coordinating railcar classification and block formation across railway networks. For example, \citet{bodin1980model} modeled network-wide classification strategies as a nonlinear mixed-integer multicommodity flow problem. Extending this perspective to broader railway operations, \citet{crainic1984tactical} jointly considered freight routing, train services, and yard classification. Subsequently, \citet{keaton1989designing} incorporated train frequencies, car routing, and train-size and yard-capacity constraints, using Lagrangian relaxation and heuristics to construct feasible operating plans.

Later studies incorporated more detailed blocking constraints and developed corresponding solution methods. \citet{newton1998constructing} considered classification limits, shipment priorities, and mileage, handling, and delay costs, developing a column-generation-based branch-and-bound method. Building on this network design perspective, \citet{barnhart2000railroad} distinguished limits on the number of blocks formed at a yard from limits on the freight volume handled there, using Lagrangian relaxation. 

Subsequent research extended railroad blocking to integrated service network design, incorporating additional train operations and service requirements. \citet{zhu2014scheduled} jointly considered train scheduling, blocking, train makeup, and time-dependent shipment routing. More closely related to our setting, \citet{arsalani2023two} integrated train formation and shipment routing under link, classification, and path-length restrictions. More recently, \citet{xin2025rail} considered heterogeneous demands under mixed blocking policies, combining time-constrained express shipments with an intree policy for ordinary freight.

These studies address different combinations of blocking, routing, scheduling, and operational constraints. Building on this literature, the IRBP considered here jointly optimizes block design, logical demand paths, and physical block routes in a static setting. It allows continuous partial service while requiring each served demand to follow a single logical path. Minimum and maximum block loads, shared physical link capacities, yard handling limits, and classification track capacities further couple these decisions. 

\subsection{DW Decomposition and Column Generation}
\label{subsec:LiteraturePricing}

Arc-based formulations represent block or service decisions and shipment movements through variables associated with individual network arcs. Flow conservation, activation constraints, and shared capacities link these decisions. Such formulations provide a direct representation of the network, as illustrated by the integrated models of \citet{yaghini2021mathematical} and \citet{arsalani2023two}. Their size grows with the network and the number of commodities, while the quality of the LP bound depends on how the interactions among the discrete decisions are represented.

DW decomposition replaces selected substructures with columns representing complete feasible patterns, coordinated through a master problem \citep{dantzig1960decomposition,lubbecke2005selected}. When the columns represent integer-feasible subproblem solutions, the reformulation convexifies those substructures and can strengthen the original LP relaxation. The improvement depends on the chosen decomposition, not simply on replacing arc variables with path variables. \citet{lubbecke2018strength} characterize this dependence for stable set and related problems, showing that different decompositions can yield substantially different relaxation strengths. For railroad blocking, a path representation also allows routing feasibility to be handled within subproblems while the master coordinates shared resources.

However, the resulting master problem may contain too many columns to enumerate explicitly. To address this difficulty, column generation maintains a restricted master problem (RMP) containing only a subset of columns. At each iteration, the RMP is solved to obtain dual prices, which guide the pricing subproblems in identifying improving columns. These columns are then added to the RMP, and the process repeats until pricing confirms that no improving columns remain \citep{lubbecke2005selected}. In railroad blocking, \citet{newton1998constructing} used shortest-path pricing to generate shipment paths within a branch-and-bound framework, while \citet{voll2013branch} developed a branch-and-price method for a European variant of the problem.

Computational performance depends heavily on pricing efficiency. In our DW formulation, logical and physical paths are coupled through master constraints involving path selection and volume variables, which complicates the design of conventional pricing algorithms based on reduced costs. In contrast, the proposed JRC measure accounts for this coupling while eliminating the dependence on path-specific dual variables, enabling more efficient pricing.

\subsection{Heuristics and Matheuristics}
\label{subsec:LiteratureMatheuristics}

Problem-specific heuristics have been central to solving large railroad blocking and service network design problems. \citet{keaton1989designing} used Lagrangian relaxation together with heuristic adjustments to obtain operating plans satisfying train capacities. \citet{barnhart2000railroad} similarly developed a Lagrangian approach that separates the blocking problem into simpler subproblems and constructs feasible plans. Neighborhood search provides another way to improve a plan without optimizing all decisions simultaneously. \citet{ahuja2007solving} developed a very large-scale neighborhood search for railroad blocking and applied it to operating plans for major US railroads.

Metaheuristics broaden the exploration of blocking plans through population search or mechanisms for moving beyond local improvements. \citet{yaghini2012population} developed a genetic algorithm for selecting blocks and assigning commodities, and \citet{yaghini2013tabu} proposed a tabu search algorithm for the railroad blocking problem. Both studies evaluated solution quality and computation time against a MIP solver. Such methods make large problems accessible to heuristic search, but the moves and candidate decisions available within each neighborhood still determine which coordinated changes can be explored.

Matheuristics use mathematical optimization to guide this search or optimize selected parts of a solution. For scheduled rail service network design, \citet{zhu2014scheduled} combined slope scaling, dynamic block generation, memory-based perturbation, and optimization over large neighborhoods. Particularly relevant to our setting, \citet{yaghini2021mathematical} developed an LP-based neighborhood search for integrated train blocking and shipment routing. Their method repeatedly solves LP relaxations, modifies objective coefficients through a tabu-inspired penalization scheme, and constructs a restricted MIP from the resulting solutions. They report high-quality solutions within reasonable computation times. These approaches illustrate how network-wide optimization information can support the selection and joint adjustment of blocking decisions.

Candidate generation provides a complementary way to control the size of the optimization problems. \citet{fang2023two} generated promising candidate blocks using shortest-path information and optimized blocking and train design through a two-stage approach. In the related multicommodity fixed-charge network flow problem, \citet{hewitt2013branch} used column generation both to construct problem restrictions and to obtain bounds, together with local search around the incumbent. Although their algorithm is exact, its use of generated restrictions to obtain feasible solutions is closely related to the design of matheuristics. 

DG-CGM follows this optimization-guided search approach by combining DW-based \textit{path generation} with restricted MIP repair. Generated logical paths expand the choices available to selected demand groups, while repair jointly adjusts their service volumes and associated block paths. The connection between \textit{path generation} and neighborhood optimization is central to the method, since improving the available paths and coordinating their use are complementary parts of the integer search.

\section{Problem Definition and Compact Formulation}
\label{sec:ProblemDefinition}

We first define the IRBP in Section~\ref{subsec:ProblemDefinition}, describing its physical and logical networks, operating decisions, feasibility constraints, and objective function. Section~\ref{subsec:CompactModel} then presents the compact arc-based mixed-integer linear formulation and explains how its constraints represent the problem requirements.

\subsection{Problem Definition}
\label{subsec:ProblemDefinition}

In the IRBP, freight demands are assigned to sequences of blocks in what we call the logical network, while each open block follows a route through the physical railway network. The logical network $G^{\mathrm Y}=(\mathcal Y,\mathcal B)$ is a directed multigraph, where $\mathcal Y$ is the set of classification yards and $\mathcal B$ is the set of candidate blocks. Each block $b\in\mathcal B$ carries railcars from yard $i_b$ to a distinct yard $j_b$ without intermediate reclassification. Multiple blocks may connect the same pair of yards. Let $\mathcal K$ denote the set of freight demands, where each demand $k\in\mathcal K$ has origin $o_k$, destination $d_k\ne o_k$, requested volume $q_k>0$, and a set $\mathcal B_k\subseteq\mathcal B$ of compatible blocks. Some demands prohibit intermediate reclassification and can therefore use only direct blocks. A blocking plan selects a set $\mathcal B^{\mathrm O}\subseteq\mathcal B$ of open blocks (\textit{block design}) and assigns each demand a served volume $v_k\in[0,q_k]$ and, if $v_k>0$, a single logical path $p_k$ from $o_k$ to $d_k$ using blocks in $\mathcal B_k\cap\mathcal B^{\mathrm O}$ (\textit{blocking sequence}). The resulting load of block $b$ is $f_b=\sum_{k\in\mathcal K:\,b\in p_k}v_k$, where unserved demands have empty paths. Every open block must satisfy its minimum load requirement $f_b\ge\underline f_b>0$. Each yard $i\in\mathcal Y$ has $T_i$ classification track slots and handling capacity $H_i$. An open block $b$ occupies $\kappa_b\in\{0,1\}$ track slots at its origin yard, depending on its commodity type, while the total reclassified volume $v_i^{\mathrm H}$, comprising freight transferred between consecutive blocks at yard $i$, cannot exceed $H_i$. The total track slots occupied by open blocks originating at $i$ cannot exceed $T_i$.

The physical network $G=(\mathcal N,\mathcal E)$ consists of physical nodes $\mathcal N\supseteq\mathcal Y$ and links $\mathcal E$, where the nodes include classification yards and transit points used only for routing. Each link $e\in\mathcal E$ permits travel in one or both directions and has capacity $C_e$, shared across both directions when applicable. Let $\mathcal A_e$ denote its permitted directed traversal arcs, and define the directed physical graph $G^{\mathrm P}=(\mathcal N,\mathcal A)$, where $\mathcal A=\bigcup_{e\in\mathcal E}\mathcal A_e$. Each arc $a\in\mathcal A$ has tail $i_a$, head $j_a$, and length $D_a>0$. For each open block $b$, \textit{block routing} selects a physical path $r_b$ from $i_b$ to $j_b$ that carries its entire load $f_b$. The path must satisfy the length limit $\sum_{a\in r_b}D_a\le\overline D_b$, and the total loads of blocks traversing each physical link $e$ cannot exceed $C_e$. Unless stated otherwise, all logical and physical paths are elementary, i.e., each node is visited at most once along a path. These physical capacity constraints couple the routing decisions with the block loads determined in the logical network.

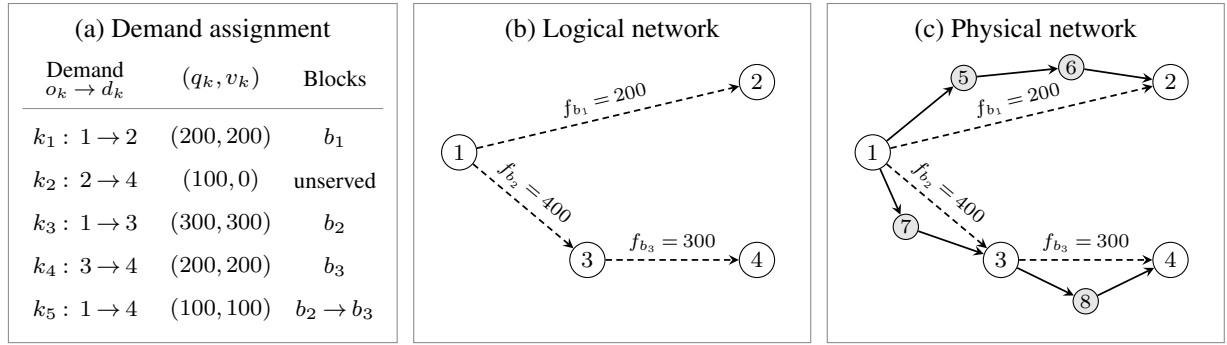
\begin{figure}[!ht]
    \caption{Demand assignment, logical network, and physical network.}
    \label{fig:BlockingDecisions}
    \centering
    \resizebox{\linewidth}{!}{%
\begin{tikzpicture}[
  x=1.08cm,y=-1.08cm,
  font=\fontsize{8.5}{9.5}\selectfont,
  yard/.style={circle,draw=black,fill=white,minimum size=4.6mm,inner sep=0pt},
  transit/.style={yard,fill=black!10,minimum size=3.3mm,font=\fontsize{7.5}{8.5}\selectfont},
  blockarc/.style={->,>=stealth,line width=.65pt,dash pattern=on 2.4pt off 1.6pt},
  blocklabel/.style={font=\fontsize{7.2}{8}\selectfont,inner sep=.5pt},
  paneltitle/.style={font=\fontsize{10}{11}\selectfont}
]
\useasboundingbox (0,0) rectangle (15,4.50);
\draw[black!45,line width=.4pt] (.15,.35) rectangle (4.85,4.45);
\draw[black!45,line width=.4pt] (5.05,.35) rectangle (9.85,4.45);
\draw[black!45,line width=.4pt] (10.05,.35) rectangle (14.85,4.45);
\node[paneltitle] at (2.5,.68) {(a) Demand assignment};
\node[paneltitle] at (7.45,.68) {(b) Logical network};
\node[paneltitle] at (12.45,.68) {(c) Physical network};
\node at (1.08,1.15) {Demand};
\node[font=\fontsize{7.5}{8.5}\selectfont] at (1.08,1.39) {\(o_k\to d_k\)};
\node at (2.72,1.25) {\((q_k,v_k)\)};
\node at (4.10,1.25) {Blocks};
\draw[black!45,line width=.35pt] (.32,1.62)--(4.68,1.62);
\foreach \id/\orig/\dest/\qty/\served/\blocks/\yy in {
  1/1/2/200/200/{b_1}/1.97,
  2/2/4/100/0/{\text{unserved}}/2.49,
  3/1/3/300/300/{b_2}/3.01,
  4/3/4/200/200/{b_3}/3.53,
  5/1/4/100/100/{b_2\to b_3}/4.05}{
  \node at (1.08,\yy) {\(k_{\id}:\ \orig\to\dest\)};
  \node at (2.72,\yy) {\((\qty,\served)\)};
  \node at (4.10,\yy) {\(\blocks\)};
}
\def\sharedyardcoordinates{1/.55/2.15,2/4.15/1.30,3/2.10/3.45,4/4.15/3.45}
\foreach \prefix/\offset in {L/5.05,P/10.05}{
  \begin{scope}[shift={(\offset,0)}]
    \foreach \id/\xx/\yy in \sharedyardcoordinates{
      \node[yard] (\prefix\id) at (\xx,\yy) {\(\id\)};
    }
  \end{scope}
}
\draw[blockarc] (L1)--node[midway,sloped,above=2pt,blocklabel] {\(f_{b_1}=200\)}(L2);
\draw[blockarc] (L1)--node[midway,sloped,above=2pt,blocklabel] {\(f_{b_2}=400\)}(L3);
\draw[blockarc] (L3)--node[midway,above=2pt,blocklabel] {\(f_{b_3}=300\)}(L4);
\draw[blockarc] (P1)--node[midway,sloped,above=2pt,blocklabel] {\(f_{b_1}=200\)}(P2);
\draw[blockarc] (P1)--node[midway,sloped,above=2pt,blocklabel] {\(f_{b_2}=400\)}(P3);
\draw[blockarc] (P3)--node[midway,above=2pt,blocklabel] {\(f_{b_3}=300\)}(P4);
\begin{scope}[shift={(10.05,0)}]
  \foreach \id/\xx/\yy in {5/1.65/1.25,6/2.95/1.12,7/.95/3.05,8/3.125/3.95}{
    \node[transit] (P\id) at (\xx,\yy) {\(\id\)};
  }
\end{scope}
\draw[->,>=stealth,line width=.65pt] (P1)--(P5);
\draw[->,>=stealth,line width=.65pt] (P5)--(P6);
\draw[->,>=stealth,line width=.65pt] (P6)--(P2);
\draw[->,>=stealth,line width=.65pt] (P1)--(P7);
\draw[->,>=stealth,line width=.65pt] (P7)--(P3);
\draw[->,>=stealth,line width=.65pt] (P3)--(P8);
\draw[->,>=stealth,line width=.65pt] (P8)--(P4);
\end{tikzpicture}
    }
    \par\smallskip
    {\footnotesize\raggedright \textit{Notes.} Dashed arrows represent open blocks, and solid arrows represent physical arcs. Classification track use is limited by the capacity at each yard (assume \(T_1=3\) and \(T_3=2\)). Yard 3 can reclassify at most 100 units of freight (assume \(H_3=100\)), and each shown arc represents a separate one-way link with capacity 500 (assume \(C_e=500\)). Assume unit arc lengths \(D_a=1\), minimum block loads \(\underline f_b=100\), and path length budgets \(\overline D_b=3\).\par}
\end{figure}

Figure~\ref{fig:BlockingDecisions} illustrates a feasible blocking plan with three open blocks. Demands $k_1$, $k_3$, $k_4$, and $k_5$ are fully served, while demand $k_2$ remains unserved. Demand $k_5$ follows the blocking sequence $b_2\to b_3$, giving block loads $f_{b_2}=v_{k_3}+v_{k_5}=400$ and $f_{b_3}=v_{k_4}+v_{k_5}=300$. All open blocks satisfy the minimum load requirement $f_b\ge\underline f_b=100$ and the physical path length limit $\sum_{a\in r_b}D_a\le\overline D_b=3$. Since blocks $b_1$ and $b_2$ originate at yard 1 and block $b_3$ at yard 3, they occupy two classification track slots at yard 1 and one at yard 3, within the respective capacities $T_1=3$ and $T_3=2$. Only demand $k_5$ requires intermediate reclassification, giving $v_3^{\mathrm H}=v_{k_5}=100\le H_3$, with no reclassification at the other yards. Finally, the load on every physical link is at most 400 units, below its capacity $C_e=500$.

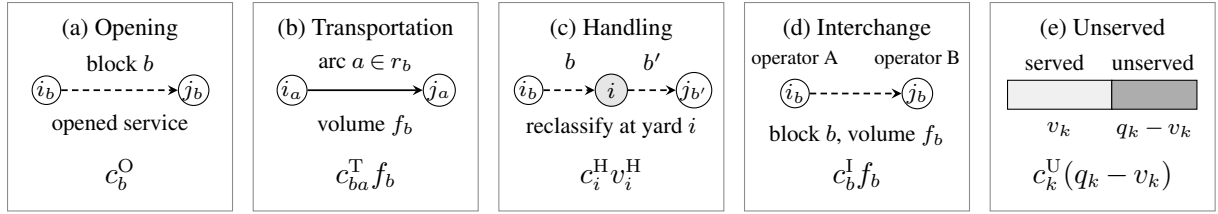
\begin{figure}[H]
    \caption{Decisions that generate the five cost components.}
    \label{fig:CostComponents}
    \centering
    \resizebox{\linewidth}{!}{%
\begin{tikzpicture}[
  x=1.08cm,y=-1.08cm,
  font=\fontsize{8.5}{9.5}\selectfont,
  point/.style={circle,draw=black,fill=white,minimum size=3.3mm,inner sep=.2pt},
  flow/.style={->,>=stealth,line width=.65pt},
  blockflow/.style={flow,dash pattern=on 2.4pt off 1.6pt},
  paneltitle/.style={font=\fontsize{9}{10}\selectfont},
  cost/.style={font=\fontsize{11}{12}\selectfont}
]
\useasboundingbox (0,5.60) rectangle (15,8.25);
\foreach \xx/\title in {.2/{(a) Opening},3.17/{(b) Transportation},6.14/{(c) Handling},9.11/{(d) Interchange},12.08/{(e) Unserved}}{
  \draw[black!45,line width=.4pt] (\xx,5.64) rectangle ({\xx+2.72},8.16);
  \node[paneltitle] at ({\xx+1.36},5.96) {\title};
}
\node[point] (Ostart) at (.67,6.70) {\(i_b\)};
\node[point] (Oend) at (2.45,6.70) {\(j_b\)};
\draw[blockflow] (Ostart)--(Oend);
\node at (1.56,6.37) {block \(b\)};
\node at (1.56,7.12) {opened service};
\node[cost] at (1.56,7.72) {\(c_b^{\mathrm O}\)};
\node[point] (Tstart) at (3.64,6.70) {\(i_a\)};
\node[point] (Tend) at (5.42,6.70) {\(j_a\)};
\draw[flow] (Tstart)--(Tend);
\node at (4.53,6.37) {arc \(a\in r_b\)};
\node at (4.53,7.12) {volume \(f_b\)};
\node[cost] at (4.53,7.72) {\(c_{ba}^{\mathrm T}f_b\)};
\node[point] (Hstart) at (6.50,6.70) {\(i_b\)};
\node[point,fill=black!10,minimum size=4.2mm] (Hyard) at (7.50,6.70) {\(i\)};
\node[point] (Hend) at (8.50,6.70) {\(j_{b'}\)};
\draw[blockflow] (Hstart)--(Hyard);
\draw[blockflow] (Hyard)--(Hend);
\node at (7.00,6.35) {\(b\)};
\node at (8.00,6.35) {\(b'\)};
\node at (7.50,7.15) {reclassify at yard \(i\)};
\node[cost] at (7.50,7.72) {\(c_i^{\mathrm H}v_i^{\mathrm H}\)};
\node[font=\fontsize{7.5}{8.5}\selectfont] at (9.72,6.33) {operator A};
\node[font=\fontsize{7.5}{8.5}\selectfont] at (11.20,6.33) {operator B};
\node[point] (Istart) at (9.72,6.74) {\(i_b\)};
\node[point] (Iend) at (11.20,6.74) {\(j_b\)};
\draw[blockflow] (Istart)--(Iend);
\node at (10.47,7.25) {block \(b\), volume \(f_b\)};
\node[cost] at (10.47,7.72) {\(c_b^{\mathrm I}f_b\)};
\node at (12.86,6.36) {served};
\node at (14.05,6.36) {unserved};
\draw[fill=black!6,line width=.4pt] (12.29,6.58) rectangle (13.55,6.93);
\draw[fill=black!32,line width=.4pt] (13.55,6.58) rectangle (14.58,6.93);
\node at (12.92,7.17) {\(v_k\)};
\node at (14.07,7.17) {\(q_k-v_k\)};
\node[cost] at (13.44,7.72) {\(c_k^{\mathrm U}(q_k-v_k)\)};
\end{tikzpicture}
    }
\end{figure}

The objective minimizes the total cost of the blocking plan, comprising the five components illustrated in Figure~\ref{fig:CostComponents}. Opening block $b$ incurs a \textit{fixed opening cost} $c_b^{\mathrm O}$, as shown in panel~(a). The \textit{transportation cost} $c_{ba}^{\mathrm T}\ge0$ applies per unit of block load carried over arc $a\in r_b$ in panel~(b), while the \textit{handling cost} $c_i^{\mathrm H}$ applies per unit of freight reclassified at intermediate yard $i$ in panel~(c). An \textit{interchange} occurs when freight is transferred between railroad operators, such as from Union Pacific to BNSF Railway, and incurs a per-unit cost $c_b^{\mathrm I}$ for block $b$, as shown in panel~(d). Finally, the \textit{unserved demand cost} $c_k^{\mathrm U}$ applies to each of the $q_k-v_k$ units of demand $k$ left unserved, as shown in panel~(e).

\subsection{Arc-Based Formulation}
\label{subsec:CompactModel}

We formulate the IRBP as an arc-based MIP, referred to as ARC, that jointly determines block design, blocking sequences, and block routing. Superscripts $\mathrm D$ and $\mathrm B$ distinguish variables associated with demand blocking sequences and physical block routes, respectively. In the logical network, for each demand $k$ and yard $i\in\mathcal Y$, let $\mathcal B_k^+(i)$ and $\mathcal B_k^-(i)$ denote the sets of compatible blocks leaving and entering $i$. For block design, binary variable $y_b$ indicates whether block $b$ is opened. For each demand $k$, binary variables $x_k^{\mathrm D}$ and $x_{kb}^{\mathrm D}$ indicate whether the demand is assigned a logical path and whether that path uses block $b\in\mathcal B_k$, respectively, while continuous variable $g_{kb}$ represents the volume of demand $k$ transported on that block. The served volumes $v_k$, block loads $f_b$, and reclassified volumes $v_i^{\mathrm H}$ defined in Section~\ref{subsec:ProblemDefinition} are also continuous variables. The set $\mathcal K_b=\{k\in\mathcal K:b\in\mathcal B_k\}$ contains the demands compatible with block $b$, and $\overline f_b=\sum_{k\in\mathcal K_b}q_k$ provides an upper bound on its load. Auxiliary continuous variables $\ell_{ki}^{\mathrm D}$ impose an ordering on the yards visited by each demand to exclude cycles in logical paths.

In the physical network, for each node $n\in\mathcal N$, let $\mathcal A^+(n)$ and $\mathcal A^-(n)$ denote the sets of arcs leaving and entering $n$. For each block $b$, binary variable $x_{ba}^{\mathrm B}$ indicates whether its physical path uses arc $a\in\mathcal A$, while continuous variable $w_{ba}$ represents the volume of block $b$ transported on that arc. Auxiliary continuous variables $\ell_{bn}^{\mathrm B}$ impose an ordering on the nodes visited by each block to exclude cycles. The indicator $\mathbf1\{\cdot\}$ equals one when its condition holds and zero otherwise. 

The ARC formulation is given in Eq.~\eqref{eq:CompactModel}. Objective~\eqref{eq:CompactObjective} minimizes the total cost of block opening, transportation, handling, interchange, and unserved demand. 
For blocking sequences, Eqs.~\eqref{eq:LogicalBalance}--\eqref{eq:LogicalFlowBalance} enforce flow conservation for the selected logical path and its associated freight volume, while Eqs.~\eqref{eq:LogicalDegrees}--\eqref{eq:LogicalOrder} prevent branching and eliminate cycles in the selected logical paths.
Eqs.~\eqref{eq:LogicalFlowBounds} and~\eqref{eq:LogicalActivation} connect logical path selections $x_{kb}^{\mathrm D}$ with freight volumes $g_{kb}$ and block openings $y_b$, while Eq.~\eqref{eq:ServedVolume} imposes lower and upper bounds on served volumes $v_k$.
For block routing, Eqs.~\eqref{eq:PhysicalBalance}--\eqref{eq:PhysicalDegrees} and \eqref{eq:PhysicalOrder} ensure that each open block follows exactly one elementary physical path from $i_b$ to $j_b$, subject to the length budget in Eq.~\eqref{eq:PhysicalRouteLength}.
Eq.~\eqref{eq:BlockVolume} aggregates the freight volumes assigned by compatible demands, 
while Eqs.~\eqref{eq:PhysicalFlowUpper}--\eqref{eq:PhysicalFlowLower} ensure that all freight assigned to each block follows its selected physical path.
Eq.~\eqref{eq:BlockVolumeBound} gives the lower and upper bounds on the freight volume of each block, while Eq.~\eqref{eq:ArcHandlingVolume} calculates the volume reclassified when freight transfers between blocks at intermediate yards.
Eqs.~\eqref{eq:ArcTrackCapacity}, \eqref{eq:ArcLinkCapacity}, and~\eqref{eq:ArcHandlingCapacity} limit classification track usage, aggregate physical link volumes, and yard handling volumes, respectively.
Finally, Eqs.~\eqref{eq:CompactOrderDomains}--\eqref{eq:CompactBinaryDomains} specify the variable domains.
We denote the optimal objective values of ARC and its LP relaxation by $z_{\mathrm{arc}}^*$ and $z_{\mathrm{arc}}^{\mathrm{lp}}$, respectively.

\noindent\begin{minipage}{\linewidth}
\begin{subequations}
\label{eq:CompactModel}
\begingroup
\normalsize
\setlength{\abovedisplayskip}{3pt}
\setlength{\belowdisplayskip}{3pt}
\setlength{\jot}{1pt}
\begin{empheq}[right={\mathmakebox[58pt][l]{}}]{alignat=1}
\mathmakebox[24pt][l]{\min}&\mathmakebox[340.5pt][l]{\displaystyle \sum_{b\in\mathcal B}c_b^{\mathrm O}y_b
+\sum_{b\in\mathcal B}\sum_{a\in\mathcal A}c_{ba}^{\mathrm T}w_{ba}
+\sum_{i\in\mathcal Y}c_i^{\mathrm H}v_i^{\mathrm H}
+\sum_{b\in\mathcal B}c_b^{\mathrm I}f_b
+\sum_{k\in\mathcal K}c_k^{\mathrm U}(q_k-v_k)}
\label{eq:CompactObjective}
\end{empheq}
\begin{empheq}[right={\mathmakebox[58pt][l]{\empheqrbrace\;\text{\scriptsize Assignment}}}]{alignat=2}
\mathmakebox[24pt][l]{\text{s.t.}}&\mathmakebox[248pt][l]{\displaystyle \sum_{b\in\mathcal B_k^+(i)}x_{kb}^{\mathrm D}
-\sum_{b\in\mathcal B_k^-(i)}x_{kb}^{\mathrm D}
=x_k^{\mathrm D}\bigl(\mathbf 1\{i=o_k\}-\mathbf 1\{i=d_k\}\bigr),}&\quad&\mathmakebox[82pt][l]{\forall k\in\mathcal K,\ i\in\mathcal Y,}
\label{eq:LogicalBalance}\\
\mathmakebox[24pt][l]{}&\mathmakebox[248pt][l]{\displaystyle \sum_{b\in\mathcal B_k^+(i)}g_{kb}
-\sum_{b\in\mathcal B_k^-(i)}g_{kb}
=v_k\bigl(\mathbf 1\{i=o_k\}-\mathbf 1\{i=d_k\}\bigr),}&\quad&\mathmakebox[82pt][l]{\forall k\in\mathcal K,\ i\in\mathcal Y,}
\label{eq:LogicalFlowBalance}\\
\mathmakebox[24pt][l]{}&\mathmakebox[248pt][l]{\displaystyle \sum_{b\in\mathcal B_k^+(i)}x_{kb}^{\mathrm D}\le x_k^{\mathrm D},\qquad
\sum_{b\in\mathcal B_k^-(i)}x_{kb}^{\mathrm D}\le x_k^{\mathrm D},}&\quad&\mathmakebox[82pt][l]{\forall k\in\mathcal K,\ i\in\mathcal Y,}
\label{eq:LogicalDegrees}\\
\mathmakebox[24pt][l]{}&\mathmakebox[248pt][l]{\displaystyle \ell_{k j_b}^{\mathrm D}\ge\ell_{k i_b}^{\mathrm D}+1-|\mathcal Y|(1-x_{kb}^{\mathrm D}),}&\quad&\mathmakebox[82pt][l]{\forall k\in\mathcal K,\ b\in\mathcal B_k,}
\label{eq:LogicalOrder}\\
\mathmakebox[24pt][l]{}&\mathmakebox[248pt][l]{\displaystyle 0\le g_{kb}\le q_kx_{kb}^{\mathrm D},}&\quad&\mathmakebox[82pt][l]{\forall k\in\mathcal K,\ b\in\mathcal B_k,}
\label{eq:LogicalFlowBounds}\\
\mathmakebox[24pt][l]{}&\mathmakebox[248pt][l]{\displaystyle 0\le v_k\le q_kx_k^{\mathrm D},}&\quad&\mathmakebox[82pt][l]{\forall k\in\mathcal K,}
\label{eq:ServedVolume}\\
\mathmakebox[24pt][l]{}&\mathmakebox[248pt][l]{\displaystyle x_{kb}^{\mathrm D}\le y_b,}&\quad&\mathmakebox[82pt][l]{\forall k\in\mathcal K,\ b\in\mathcal B_k.}
\label{eq:LogicalActivation}
\end{empheq}
\begin{empheq}[right={\mathmakebox[58pt][l]{\empheqrbrace\;\text{\scriptsize Paths}}}]{alignat=2}
\mathmakebox[24pt][l]{}&\mathmakebox[248pt][l]{\displaystyle \sum_{a\in\mathcal A^+(n)}x_{ba}^{\mathrm B}
-\sum_{a\in\mathcal A^-(n)}x_{ba}^{\mathrm B}
=y_b\bigl(\mathbf 1\{n=i_b\}-\mathbf 1\{n=j_b\}\bigr),}&\quad&\mathmakebox[82pt][l]{\forall b\in\mathcal B,\ n\in\mathcal N,}
\label{eq:PhysicalBalance}\\
\mathmakebox[24pt][l]{}&\mathmakebox[248pt][l]{\displaystyle \sum_{a\in\mathcal A^+(n)}x_{ba}^{\mathrm B}\le y_b,\qquad
\sum_{a\in\mathcal A^-(n)}x_{ba}^{\mathrm B}\le y_b,}&\quad&\mathmakebox[82pt][l]{\forall b\in\mathcal B,\ n\in\mathcal N,}
\label{eq:PhysicalDegrees}\\
\mathmakebox[24pt][l]{}&\mathmakebox[248pt][l]{\displaystyle \sum_{a\in\mathcal A}D_ax_{ba}^{\mathrm B}\le\overline D_by_b,}&\quad&\mathmakebox[82pt][l]{\forall b\in\mathcal B,}
\label{eq:PhysicalRouteLength}\\
\mathmakebox[24pt][l]{}&\mathmakebox[248pt][l]{\displaystyle f_b=\sum_{k\in\mathcal K_b}g_{kb},}&\quad&\mathmakebox[82pt][l]{\forall b\in\mathcal B,}
\label{eq:BlockVolume}\\
\mathmakebox[24pt][l]{}&\mathmakebox[248pt][l]{\displaystyle \ell_{b j_a}^{\mathrm B}\ge\ell_{b i_a}^{\mathrm B}+1-|\mathcal N|(1-x_{ba}^{\mathrm B}),}&\quad&\mathmakebox[82pt][l]{\forall b\in\mathcal B,\ a\in\mathcal A,}
\label{eq:PhysicalOrder}\\
\mathmakebox[24pt][l]{}&\mathmakebox[248pt][l]{\displaystyle 0\le w_{ba}\le \overline f_bx_{ba}^{\mathrm B},\qquad w_{ba}\le f_b,}&\quad&\mathmakebox[82pt][l]{\forall b\in\mathcal B,\ a\in\mathcal A,}
\label{eq:PhysicalFlowUpper}\\
\mathmakebox[24pt][l]{}&\mathmakebox[248pt][l]{\displaystyle w_{ba}\ge f_b-\overline f_b(y_b-x_{ba}^{\mathrm B}),}&\quad&\mathmakebox[82pt][l]{\forall b\in\mathcal B,\ a\in\mathcal A,}
\label{eq:PhysicalFlowLower}\\
\mathmakebox[24pt][l]{}&\mathmakebox[248pt][l]{\displaystyle \underline f_by_b\le f_b\le \overline f_by_b,}&\quad&\mathmakebox[82pt][l]{\forall b\in\mathcal B.}
\label{eq:BlockVolumeBound}
\end{empheq}
\begin{empheq}[right={\mathmakebox[58pt][l]{\empheqrbrace\;\text{\scriptsize Capacity}}}]{alignat=2}
\mathmakebox[24pt][l]{}&\mathmakebox[248pt][l]{\displaystyle v_i^{\mathrm H}=\sum_{k\in\mathcal K}
\left(\sum_{b\in\mathcal B_k^+(i)}g_{kb}-\mathbf 1\{i=o_k\}v_k\right),}&\quad&\mathmakebox[82pt][l]{\forall i\in\mathcal Y,}
\label{eq:ArcHandlingVolume}\\
\mathmakebox[24pt][l]{}&\mathmakebox[248pt][l]{\displaystyle \sum_{b\in\mathcal B:i_b=i}\kappa_by_b\le T_i,}&\quad&\mathmakebox[82pt][l]{\forall i\in\mathcal Y,}
\label{eq:ArcTrackCapacity}\\
\mathmakebox[24pt][l]{}&\mathmakebox[248pt][l]{\displaystyle \sum_{b\in\mathcal B}\sum_{a\in\mathcal A_e}w_{ba}\le C_e,}&\quad&\mathmakebox[82pt][l]{\forall e\in\mathcal E,}
\label{eq:ArcLinkCapacity}\\
\mathmakebox[24pt][l]{}&\mathmakebox[248pt][l]{\displaystyle v_i^{\mathrm H}\le H_i,}&\quad&\mathmakebox[82pt][l]{\forall i\in\mathcal Y.}
\label{eq:ArcHandlingCapacity}
\end{empheq}
\begin{empheq}[right={\mathmakebox[58pt][l]{\empheqrbrace\;\text{\scriptsize Domains}}}]{alignat=2}
\mathmakebox[24pt][l]{}&\mathmakebox[248pt][l]{\displaystyle 0\le\ell_{ki}^{\mathrm D}\le|\mathcal Y|-1,\qquad
0\le\ell_{bn}^{\mathrm B}\le|\mathcal N|-1,}&\quad&\mathmakebox[82pt][l]{}
\label{eq:CompactOrderDomains}\\
\mathmakebox[24pt][l]{}&\mathmakebox[248pt][l]{\displaystyle x_{kb}^{\mathrm D},\ x_k^{\mathrm D},\ y_b,\ x_{ba}^{\mathrm B}\in\{0,1\}.}&\quad&\mathmakebox[82pt][l]{}
\label{eq:CompactBinaryDomains}
\end{empheq}
\endgroup
\end{subequations}
\end{minipage}
\par\smallskip

\section{Dantzig--Wolfe Formulation and Column Generation}
\label{sec:DWReformulation}

The weak LP relaxation of ARC, as reflected in the large lower-bound gaps reported in Section~\ref{subsec:genratedInstanceExperiments}, limits its usefulness for evaluating heuristic solutions and may provide less informative dual prices for guiding integer search. To address this issue, we further develop a DW formulation in Section~\ref{subsec:MasterReformulation}. Its LP relaxation is solved to optimality through column generation using the logical and physical path pricing algorithms developed in Section~\ref{subsec:ExactPricing}.

\subsection{Dantzig--Wolfe Formulation}
\label{subsec:MasterReformulation}

The DW formulation represents demand blocking sequences and block routing decisions using complete logical and physical paths, respectively. For each demand $k\in\mathcal K$, let $\mathcal P_k$ denote the set of all elementary logical paths from $o_k$ to $d_k$ using compatible blocks in $\mathcal B_k$.
Similarly, for each block $b\in\mathcal B$, let $\mathcal R_b$ contain all elementary physical paths from $i_b$ to $j_b$ satisfying the length budget $\overline D_b$, with $\mathcal R_b=\emptyset$ if $\overline f_b<\underline f_b$.

To represent demand assignments and block routing using complete paths, we introduce binary variables $x_{kp}^{\mathrm D}$ and $x_{br}^{\mathrm B}$ to indicate whether demand $k$ selects logical path $p$ and block $b$ follows physical path $r$, respectively, while continuous variables $v_{kp}^{\mathrm D}$ and $v_{br}^{\mathrm B}$ denote the corresponding freight volumes. Let $\iota_{bp}$ equal one if logical path $p$ uses block $b$, and $\delta_{ar}$ equal one if physical path $r$ uses arc $a$, and zero otherwise. The coefficient $\eta_{ip}=\mathbf1\{i\ne o_k\}\sum_{b\in\mathcal B_k^+(i)}\iota_{bp}$ indicates whether freight is reclassified at yard $i$ along logical path $p$. We preserve the cost structure of ARC by expressing its volume-dependent terms at the path level. Specifically, the logical path unit cost $c_{kp}^{\mathrm D}=\sum_{i\in\mathcal Y}c_i^{\mathrm H}\eta_{ip}-c_k^{\mathrm U}$ accounts for handling costs at intermediate yards and avoided unserved demand costs, while the physical path unit cost $c_{br}^{\mathrm B}=c_b^{\mathrm I}+\sum_{a\in\mathcal A}c_{ba}^{\mathrm T}\delta_{ar}$ includes interchange and transportation costs.

The DW formulation is given in Eq.~\eqref{eq:DWMaster}. Objective~\eqref{eq:DWObjective} minimizes the same total cost of block opening, transportation, handling, interchange, and unserved demand in terms of path variables.
Since logical and physical paths are predefined in $\mathcal P_k$ and $\mathcal R_b$, path connectivity, elementarity, and physical path length limits are incorporated directly into these path sets.
Eq.~\eqref{eq:DWBlockBalance} ensures that the freight volume assigned to each block matches the volume transported along its physical route.
For resource capacities, Eqs.~\eqref{eq:DWLinkCapacity}--\eqref{eq:DWTrackCapacity} limit aggregate physical link volumes, yard handling volumes, and classification track usage, respectively.
Eqs.~\eqref{eq:DWDemandVolume}--\eqref{eq:DWBlockPathVolume} connect path selections with their corresponding freight volumes, imposing demand quantity bounds on logical paths and minimum and maximum load bounds on physical paths.
Eqs.~\eqref{eq:DWDemandConvexity}--\eqref{eq:DWBlockConvexity} restrict each demand to at most one logical path and each block to at most one physical path, while Eq.~\eqref{eq:DWServiceActivation} ensures that selected logical paths use only open blocks.
Finally, Eqs.~\eqref{eq:DWBinaryDomains}--\eqref{eq:DWVolumeDomains} specify the variable domains.
We obtain the LP relaxation by replacing the binary domains with $x_{kp}^{\mathrm D},x_{br}^{\mathrm B}\ge0$. We denote the optimal objective values of DW and its LP relaxation by $z_{\mathrm{dw}}^*$ and $z_{\mathrm{dw}}^{\mathrm{lp}}$, respectively.

\begin{theorem}
\label{thm:DantzigWolfeReformulation}
The \textup{DW} formulation~\eqref{eq:DWMaster} is an exact mixed-integer Dantzig--Wolfe reformulation of the \textup{ARC} formulation~\eqref{eq:CompactModel}, so \(z_{\mathrm{arc}}^*=z_{\mathrm{dw}}^*\). Moreover, every feasible solution of the DW LP relaxation can be mapped to a feasible solution of the ARC LP relaxation with the same objective value, implying \(z_{\mathrm{arc}}^{\mathrm{lp}}\le z_{\mathrm{dw}}^{\mathrm{lp}}\). There exist instances for which this inequality is strict.
\end{theorem}
\proof{Proof.} See Section~\ref{ec:DantzigWolfeReformulation} of the electronic companion. \Halmos
\endproof

\noindent\begin{minipage}{\linewidth}
\begin{subequations}
\label{eq:DWMaster}
\begingroup
\normalsize
\setlength{\abovedisplayskip}{3pt}
\setlength{\belowdisplayskip}{3pt}
\setlength{\jot}{1pt}
\begin{empheq}[right={\mathmakebox[58pt][l]{}}]{alignat=2}
\mathmakebox[60pt][l]{\min}&\mathmakebox[210pt][l]{\displaystyle \sum_{k\in\mathcal K}q_kc_k^{\mathrm U}+\sum_{k\in\mathcal K}\sum_p c_{kp}^{\mathrm D}v_{kp}^{\mathrm D}+\sum_{b\in\mathcal B}\sum_r(c_b^{\mathrm O}x_{br}^{\mathrm B}+c_{br}^{\mathrm B}v_{br}^{\mathrm B})} &\quad&\mathmakebox[94pt][l]{}
\label{eq:DWObjective}
\end{empheq}
\begin{empheq}[right={\mathmakebox[58pt][l]{\empheqrbrace\;\text{\scriptsize Balance}}}]{alignat=2}
\mathmakebox[60pt][l]{\text{s.t.}}&\mathmakebox[210pt][l]{\displaystyle \sum_r v_{br}^{\mathrm B}-\sum_{k\in\mathcal K_b}\sum_p\iota_{bp}v_{kp}^{\mathrm D}=0,}&\quad&\mathmakebox[94pt][l]{\forall b\in\mathcal B,}
\label{eq:DWBlockBalance}
\end{empheq}
\begin{empheq}[right={\mathmakebox[58pt][l]{\empheqrbrace\;\text{\scriptsize Capacity}}}]{alignat=2}
\mathmakebox[60pt][l]{}&\mathmakebox[210pt][l]{\displaystyle \sum_{b\in\mathcal B}\sum_r\left(\sum_{a\in\mathcal A_e}\delta_{ar}\right)v_{br}^{\mathrm B}\le C_e,}&\quad&\mathmakebox[94pt][l]{\forall e\in\mathcal E,}
\label{eq:DWLinkCapacity}\\
\mathmakebox[60pt][l]{}&\mathmakebox[210pt][l]{\displaystyle \sum_{k\in\mathcal K}\sum_p\eta_{ip}v_{kp}^{\mathrm D}\le H_i,}&\quad&\mathmakebox[94pt][l]{\forall i\in\mathcal Y,}
\label{eq:DWHandlingCapacity}\\
\mathmakebox[60pt][l]{}&\mathmakebox[210pt][l]{\displaystyle \sum_{b\in\mathcal B:i_b=i}\kappa_b\sum_r x_{br}^{\mathrm B}\le T_i,}&\quad&\mathmakebox[94pt][l]{\forall i\in\mathcal Y,}
\label{eq:DWTrackCapacity}
\end{empheq}
\begin{empheq}[right={\mathmakebox[58pt][l]{\empheqrbrace\;\text{\scriptsize Volume}}}]{alignat=2}
\mathmakebox[60pt][l]{}&\mathmakebox[210pt][l]{\displaystyle 0\le v_{kp}^{\mathrm D}\le q_kx_{kp}^{\mathrm D},}&\quad&\mathmakebox[94pt][l]{\forall k\in\mathcal K,\ p\in\mathcal P_k,}
\label{eq:DWDemandVolume}\\
\mathmakebox[60pt][l]{}&\mathmakebox[210pt][l]{\displaystyle \underline f_bx_{br}^{\mathrm B}\le v_{br}^{\mathrm B}\le \overline f_bx_{br}^{\mathrm B},}&\quad&\mathmakebox[94pt][l]{\forall b\in\mathcal B,\ r\in\mathcal R_b,}
\label{eq:DWBlockPathVolume}
\end{empheq}
\begin{empheq}[right={\mathmakebox[58pt][l]{\empheqrbrace\;\text{\scriptsize Convexity}}}]{alignat=2}
\mathmakebox[60pt][l]{}&\mathmakebox[210pt][l]{\displaystyle \sum_p x_{kp}^{\mathrm D}\le1,}&\quad&\mathmakebox[94pt][l]{\forall k\in\mathcal K,}
\label{eq:DWDemandConvexity}\\
\mathmakebox[60pt][l]{}&\mathmakebox[210pt][l]{\displaystyle \sum_r x_{br}^{\mathrm B}\le1,}&\quad&\mathmakebox[94pt][l]{\forall b\in\mathcal B,}
\label{eq:DWBlockConvexity}
\end{empheq}
\begin{empheq}[right={\mathmakebox[58pt][l]{\empheqrbrace\;\text{\scriptsize Activation}}}]{alignat=2}
\mathmakebox[60pt][l]{}&\mathmakebox[210pt][l]{\displaystyle \sum_p\iota_{bp}x_{kp}^{\mathrm D}\le\sum_r x_{br}^{\mathrm B},}&\quad&\mathmakebox[94pt][l]{\forall k\in\mathcal K,\ b\in\mathcal B_k,}
\label{eq:DWServiceActivation}
\end{empheq}
\begin{empheq}[right={\mathmakebox[58pt][l]{\empheqrbrace\;\text{\scriptsize Domains}}}]{alignat=2}
\mathmakebox[60pt][l]{}&\mathmakebox[210pt][l]{\displaystyle x_{kp}^{\mathrm D},\ x_{br}^{\mathrm B}\in\{0,1\},}&\quad&\mathmakebox[94pt][l]{}
\label{eq:DWBinaryDomains}\\
\mathmakebox[60pt][l]{}&\mathmakebox[210pt][l]{\displaystyle v_{kp}^{\mathrm D},\ v_{br}^{\mathrm B}\ge0.}&\quad&\mathmakebox[94pt][l]{}
\label{eq:DWVolumeDomains}
\end{empheq}
\endgroup
\end{subequations}
\end{minipage}
\par\smallskip

\subsection{Column Generation}
\label{subsec:ExactPricing}

The DW formulation contains a large number of logical and physical path variables, making direct solution of its LP relaxation impractical. We therefore use column generation, which maintains an RMP and iteratively adds improving columns identified by pricing subproblems. Specifically, Section~\ref{subsubsec:ExactDemandPricing} develops logical pricing to generate blocking sequences for each demand, while Section~\ref{subsubsec:ExactBlockPricing} presents physical pricing to identify feasible routes for each block. Both procedures use the dual prices of the RMP to identify improving columns, which are then added for reoptimization.

When both pricing procedures are solved exactly and no improving columns remain, the resulting DW LP optimum provides a valid lower bound for the IRBP. 
When heuristic pricing is used, the generated columns can still provide promising paths for improving integer solutions, although the resulting LP solution does not certify a valid lower bound. 
This distinction allows column generation to support both exact bound computation and the heuristic search developed in Section~\ref{sec:PrimalHeuristic}.

Let \(\sigma_b\), \(\mu_e\), \(\mu_i^{\mathrm H}\), \(\mu_i^{\mathrm T}\), \(\pi_k\), \(\pi_b\), and \(\gamma_{kb}\) denote the dual variables associated with Eqs.~\eqref{eq:DWBlockBalance}, \eqref{eq:DWLinkCapacity}, \eqref{eq:DWHandlingCapacity}, \eqref{eq:DWTrackCapacity}, \eqref{eq:DWDemandConvexity}, \eqref{eq:DWBlockConvexity}, and \eqref{eq:DWServiceActivation}, respectively. The block-balance duals \(\sigma_b\) are unrestricted in sign, while \(\mu_e,\mu_i^{\mathrm H},\mu_i^{\mathrm T},\pi_k,\pi_b,\gamma_{kb}\le0\).

\subsubsection{Logical Path Generation.}
\label{subsubsec:ExactDemandPricing}

For each demand $k\in\mathcal K$ and logical path $p\in\mathcal P_k$, define
\begin{equation}
\bar c_{kp}^{\mathrm D,x}=-\pi_k-\sum_{b\in\mathcal B_k}\gamma_{kb}\iota_{bp},
\qquad
\bar c_{kp}^{\mathrm D,v}=c_{kp}^{\mathrm D}+\sum_{b\in\mathcal B_k}\sigma_b\iota_{bp}-\sum_{i\in\mathcal Y}\mu_i^{\mathrm H}\eta_{ip}.
\label{eq:LogicalAdjustedCost}
\end{equation}
For a path included in the RMP, let $\lambda_{kp}\le0$ denote the dual variable associated with $v_{kp}^{\mathrm D}-q_kx_{kp}^{\mathrm D}\le0$ in Eq.~\eqref{eq:DWDemandVolume}. Using the coefficients in Eq.~\eqref{eq:LogicalAdjustedCost}, the reduced costs of the corresponding selection and volume variables are
\begin{equation}
\begin{aligned}
\operatorname{rc}(x_{kp}^{\mathrm D})
&=\bar c_{kp}^{\mathrm D,x}+q_k\lambda_{kp},\\
\operatorname{rc}(v_{kp}^{\mathrm D})
&=\bar c_{kp}^{\mathrm D,v}-\lambda_{kp}.
\end{aligned}
\label{eq:LogicalReducedCost}
\end{equation}

Conventional pricing uses the reduced costs in Eq.~\eqref{eq:LogicalReducedCost} to identify improving logical paths. Labeling algorithms typically extend partial paths while accumulating reduced costs, allowing unpromising labels to be discarded before reaching the destination. However, applying this approach directly to our formulation is challenging because each logical path introduces both a selection variable and a volume variable, linked by $v_{kp}^{\mathrm D}\le q_kx_{kp}^{\mathrm D}$, whose reduced costs depend on the path-specific dual variable $\lambda_{kp}$. Since this constraint has a zero right-hand side and involves no existing RMP variables, $\lambda_{kp}\le0$ can be chosen without affecting the current RMP dual objective or existing reduced costs. Thus, conventional pricing would require finding a suitable $\lambda_{kp}$ for each candidate path such that the reduced costs of both variables are nonnegative whenever possible. This path-dependent choice complicates the incremental evaluation of reduced costs and the design of dominance rules in labeling algorithms.

To address this difficulty, we define the JRC to combine the reduced costs of the selection and volume variables without depending on $\lambda_{kp}$. Specifically, the JRC of logical path $p$ is defined as
\begin{equation}
\begin{aligned}
\widetilde c_{kp}^{\mathrm D}
&=\bar c_{kp}^{\mathrm D,x}
+q_k\bar c_{kp}^{\mathrm D,v}.
\end{aligned}
\label{eq:LogicalFullLoadScore}
\end{equation}
By Proposition~\ref{prop:ExactPricing}, a nonnegative minimum JRC in Eq.~\eqref{eq:LogicalFullLoadScore} over all feasible logical paths for every demand certifies that no improving logical columns remain. 
Moreover, the JRC can be decomposed into an initial cost \(\widetilde c_k^{\mathrm D,\mathrm{init}}=-q_kc_k^{\mathrm U}-\pi_k\) and a block extension cost
\begin{equation}
\widetilde c_{kb}^{\mathrm D}
=q_k\sigma_b-\gamma_{kb}+q_k\mathbf1\{i_b\ne o_k\}(c_{i_b}^{\mathrm H}-\mu_{i_b}^{\mathrm H}).
\label{eq:LogicalExtensionCost}
\end{equation}
Using Eq.~\eqref{eq:LogicalExtensionCost}, the JRC of logical path \(p\) is \(\widetilde c_{kp}^{\mathrm D}=\widetilde c_k^{\mathrm D,\mathrm{init}}+\sum_{b\in\mathcal B_k}\iota_{bp}\widetilde c_{kb}^{\mathrm D}\).

We use a labeling algorithm to find minimum-JRC elementary logical paths. For each demand \(k\), a label \(L=(i,S,c)\) records the current yard \(i\), visited yards \(S\), and accumulated JRC \(c\). We initialize \(L_0=(o_k,\{o_k\},\widetilde c_k^{\mathrm D,\mathrm{init}})\). For each compatible block $b\in\mathcal B_k^+(i)$ with $j_b\notin S$, the label is extended to $L'=(j_b,S\cup\{j_b\},c+\widetilde c_{kb}^{\mathrm D})$. A label reaching $d_k$ represents a complete logical path with JRC $c$. The dominance rule used to reduce the labeling search space is given in Proposition~\ref{prop:LogicalDominance}.

\begin{proposition}
\label{prop:LogicalDominance}
Fix demand \(k\). For two feasible labels \(L_1=(i,S_1,c_1)\) and \(L_2=(i,S_2,c_2)\) at the same yard, write \(L_1\preceq L_2\) if \(S_1\subseteq S_2\) and \(c_1\le c_2\). Then \(L_1\) dominates \(L_2\), and discarding \(L_2\) while retaining \(L_1\) preserves the minimum JRC over all feasible logical paths for demand \(k\).
\end{proposition}
\proof{Proof.} See Section~\ref{ec:LogicalDominance} of the electronic companion. \Halmos
\endproof

\subsubsection{Physical Path Generation.}
\label{subsubsec:ExactBlockPricing}

For each block \(b\in\mathcal B\) and physical path \(r\in\mathcal R_b\), define
\begin{equation}
\begin{aligned}
\bar c_{br}^{\mathrm B,x}
&=c_b^{\mathrm O}-\pi_b-\kappa_b\mu_{i_b}^{\mathrm T}
+\sum_{k\in\mathcal K_b}\gamma_{kb},\\
\bar c_{br}^{\mathrm B,v}
&=c_{br}^{\mathrm B}-\sigma_b
-\sum_{a\in\mathcal A}\mu_{e_a}\delta_{ar},
\end{aligned}
\label{eq:PhysicalReducedCost}
\end{equation}
where \(e_a\) denotes the unique physical link containing arc \(a\). For each block \(b\) and physical path \(r\in\mathcal R_b\), we define its JRC at load \(t\) using the coefficients in Eq.~\eqref{eq:PhysicalReducedCost} as
\begin{equation}
\widetilde c_{br}^{\mathrm B}(t)
=\bar c_{br}^{\mathrm B,x}+t\bar c_{br}^{\mathrm B,v}.
\label{eq:PhysicalFullLoadScore}
\end{equation}
A selected physical path carries a load \(t\in[\underline f_b,\overline f_b]\). Since the JRC in Eq.~\eqref{eq:PhysicalFullLoadScore} is affine in \(t\), Proposition~\ref{prop:ExactPricing} establishes that it suffices to consider the two endpoints \(t\in\{\underline f_b,\overline f_b\}\). If both endpoint JRCs are nonnegative, path \(r\) cannot yield an improving column and need not be added to the RMP. Thus, physical pricing is complete when this condition holds for every physical path.

Moreover, the JRC can be decomposed into an initial cost \(\widetilde c_{bt}^{\mathrm B,\mathrm{init}}=c_b^{\mathrm O}+t(c_b^{\mathrm I}-\sigma_b)-\pi_b-\kappa_b\mu_{i_b}^{\mathrm T}+\sum_{k\in\mathcal K_b}\gamma_{kb}\) and an arc extension cost
\begin{equation}
\widetilde c_{bta}^{\mathrm B}
=t\left(c_{ba}^{\mathrm T}-\mu_{e_a}\right),
\label{eq:BlockExtensionCost}
\end{equation}
Using Eq.~\eqref{eq:BlockExtensionCost}, the JRC of physical path \(r\) is \(\widetilde c_{br}^{\mathrm B}(t)=\widetilde c_{bt}^{\mathrm B,\mathrm{init}}+\sum_{a\in\mathcal A}\delta_{ar}\widetilde c_{bta}^{\mathrm B}\). Since the initial cost is independent of \(r\) and the arc extension costs are proportional to \(t\), both endpoint JRCs rank physical paths identically. Thus, one search at \(\underline f_b\) identifies a minimum-JRC path for both loads, and its JRC at \(\overline f_b\) can be evaluated directly.

We use a labeling algorithm to find minimum-JRC elementary physical paths. For each block \(b\) with \(\overline f_b\ge\underline f_b\), since \(\mathcal R_b=\emptyset\) otherwise, a label \(L=(n,S,d,c)\) records the current node \(n\), visited nodes \(S\), accumulated physical length \(d\), and JRC \(c\). We initialize \(L_0=(i_b,\{i_b\},0,\widetilde c_{bt}^{\mathrm B,\mathrm{init}})\). For each arc \(a\in\mathcal A^+(n)\) with \(j_a\notin S\), the label extends to \(L'=(j_a,S\cup\{j_a\},d+D_a,c+\widetilde c_{bta}^{\mathrm B})\), subject to the length budget \(\overline D_b\). A label reaching \(j_b\) represents a complete physical path with JRC \(c\). The dominance rule used to reduce the labeling search space is given in Proposition~\ref{prop:PhysicalDominance}.

\begin{proposition}
\label{prop:PhysicalDominance}
Fix block \(b\), load \(t\), and the current dual prices. For feasible labels \(L_1=(n,S_1,d_1,c_1)\) and \(L_2=(n,S_2,d_2,c_2)\) at the same node, write \(L_1\preceq L_2\) if \(S_1\subseteq S_2\), \(d_1\le d_2\), and \(c_1\le c_2\). Then \(L_1\) dominates \(L_2\), and discarding \(L_2\) while retaining \(L_1\) preserves the minimum JRC over all feasible physical paths for block \(b\) at load \(t\).
\end{proposition}
\proof{Proof.} See Section~\ref{ec:PhysicalDominance} of the electronic companion. \Halmos
\endproof

To further reduce the search space, we use shortest-path distances to identify labels that cannot satisfy the physical path length budget. Let \(\underline D_b(n)\) denote the shortest directed length from node \(n\) to \(j_b\), ignoring visited nodes, with value \(+\infty\) if unreachable. These distances are computed using Dijkstra's algorithm on the reversed physical graph. We discard a label when \(d+\underline D_b(n)>\overline D_b\) and reject an extension along arc \(a\) when \(d+D_a+\underline D_b(j_a)>\overline D_b\). Proposition~\ref{prop:CompletionBounds} establishes the validity of these tests. Since the distances depend only on the physical graph and destination, they can be reused across pricing calls.

\begin{proposition}
\label{prop:CompletionBounds}
The distances \(\underline D_b(n)\) are lower bounds on every feasible suffix length. The resulting prefix and extension tests discard only prefixes without a feasible completion and preserve the minimum physical path JRC.
\end{proposition}
\proof{Proof.} See Section~\ref{ec:CompletionBounds} of the electronic companion. \Halmos
\endproof

\subsubsection{Overall Procedure.}
\label{subsubsec:PricingProcedure}
Algorithm~\ref{alg:PricingFramework} solves the RMP to optimality and uses its dual prices for exact logical and physical pricing. The returned sets \(\mathcal P^+\) and \(\mathcal R^+\) contain logical pairs \((k,p)\) and physical pairs \((b,r)\) with negative JRC. Each path is added once with its selection and volume variables, local bounds, and coupling coefficients. We reoptimize and repeat until both sets are empty. 

\begin{algorithm}[htbp]
\small
\caption{Exact Column Generation}
\label{alg:PricingFramework}
\begin{algorithmic}[1]
\State Initialize a feasible RMP.
\Repeat
    \State Solve the RMP and obtain dual prices.
    \State \(\mathcal P^+\gets\Call{ExactLogicalPricing}{}\).
    \State \(\mathcal R^+\gets\Call{ExactPhysicalPricing}{}\).
    \State Add the generated paths and associated constraints to the RMP.
\Until{\(\mathcal P^+=\emptyset\) and \(\mathcal R^+=\emptyset\)}
\State \Return optimal DW LP solution and dual prices.
\end{algorithmic}
\end{algorithm}

\begin{proposition}
\label{prop:ExactPricing}
Suppose the RMP is solved to optimality. If \(\widetilde c_{kp}^{\mathrm D}\ge0\) for every demand \(k\) and logical path \(p\in\mathcal P_k\), and \(\widetilde c_{br}^{\mathrm B}(t)\ge0\) for every block \(b\), physical path \(r\in\mathcal R_b\), and \(t\in\{\underline f_b,\overline f_b\}\), then the RMP solution is optimal for the DW LP relaxation, and its objective value provides a valid lower bound on $z_{\mathrm{dw}}^*$.
\end{proposition}
\proof{Proof.} See Section~\ref{ec:ExactPricing} of the electronic companion. \Halmos
\endproof
\par

\section{Dual-Guided Column-Generation-Based Matheuristic}
\label{sec:PrimalHeuristic}

\begin{figure}[!htbp]
\caption{The framework of DG-CGM.}
\label{fig:PrimalFramework}
\centering
\vspace{1mm}
\begin{tikzpicture}[
 x=1cm,y=1cm,
 font=\normalfont\rmfamily\fontsize{10}{12}\selectfont,
 text=black,
 flow/.style={->,>=stealth,line width=.5pt,rounded corners=1.5pt},
 rule/.style={line width=.5pt,rounded corners=1.5pt},
 stage/.style={draw=black!75,line width=.5pt,rounded corners=1.5pt,
  fill=black!2,minimum width=2.50cm,minimum height=.70cm,
  inner sep=3pt,align=center},
 terminal/.style={draw=black!75,line width=.5pt,rounded corners=8pt,
  minimum width=2.50cm,minimum height=.65cm,inner sep=3pt,align=center},
 question/.style={font=\normalfont\rmfamily\fontsize{10}{12}\selectfont,align=center},
 branch/.style={font=\normalfont\rmfamily\fontsize{9}{11}\selectfont,fill=white,inner sep=1.5pt}
]
\node[terminal] (construct) at (-2.65,4.40) {\textit{Construction}};
\node[stage,minimum width=1.60cm] (destroy) at (-2.65,2.60) {\textit{Destroy}};
\node[stage] (local) at (.40,2.60) {\textit{Local repair}};
\node[stage] (pg) at (4.50,2.60) {\textit{Path generation}};
\node[stage] (joint) at (8.60,2.60) {\textit{Joint repair}};

\foreach \xx/\name/\words in {.40/localtest/{Local\\stagnation?},4.50/roundtest/{Round\\stagnation?},8.60/endtest/{Terminate?}}{
 \draw[line width=.5pt,draw=black!75,fill=white]
  (\xx,1.60)--({\xx+1.35},.90)--(\xx,.20)--({\xx-1.35},.90)--cycle;
 \node[question] (\name) at (\xx,.90) {\words};
}
\node[terminal,minimum width=3.90cm] (finish) at (8.60,-.65) {Return best feasible plan};

\draw[flow] (construct.south)--(destroy.north);
\draw[flow] (destroy.east)--(local.west);
\draw[flow] (local.south)--(.40,1.60);
\draw[flow] (-.95,.90)--(-3.95,.90)--(-3.95,2.60)--(destroy.west);
\node[branch] at (-2.45,1.12) {No};
\draw[flow] (1.75,.90)--(3.15,.90);
\node[branch] at (2.45,1.12) {Yes};
\draw[flow] (4.50,1.60)--(pg.south);
\node[branch,anchor=west] at (4.60,1.92) {No};
\draw[flow] (5.85,.90)--(6.85,.90)--(6.85,2.60)--(joint.west);
\node[branch] at (6.35,1.12) {Yes};
\draw[flow] (joint.south)--(8.60,1.60);
\draw[flow] (8.60,.20)--(finish.north);
\node[branch,anchor=west] at (8.70,-.04) {Yes};

\draw[rule] (pg.north)--(4.50,3.55);
\draw[rule] (9.95,.90)--(10.65,.90)--(10.65,3.55);
\node[branch] at (10.30,1.12) {No};
\draw[rule] (10.65,3.55)--(-2.65,3.55);
\fill (4.50,3.55) circle (.7pt);
\fill (-2.65,3.55) circle (.7pt);
\end{tikzpicture}
\end{figure}
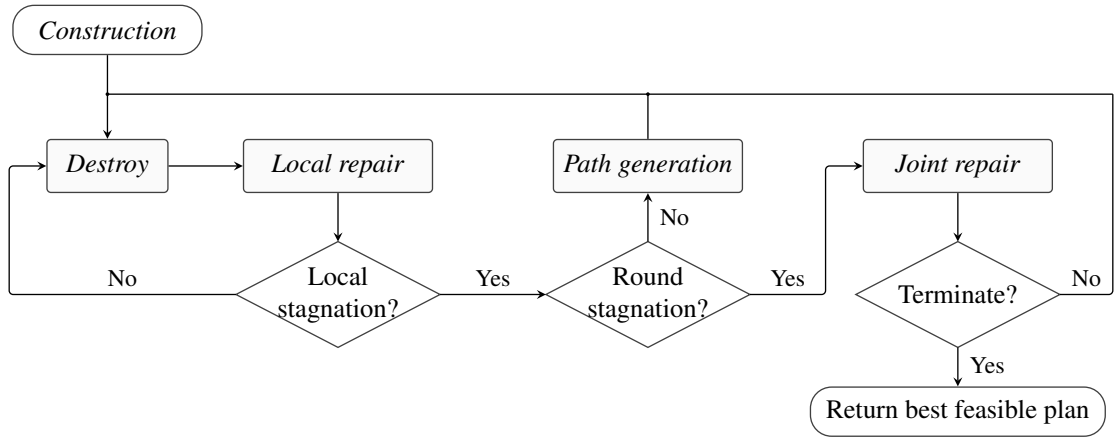

DG-CGM improves feasible blocking plans by combining restricted MIP-based neighborhood search with dual-guided \textit{path generation}. The framework consists of five components, as illustrated in Figure~\ref{fig:PrimalFramework}. First, \textit{construction} builds an initial feasible plan and candidate physical paths (Section~\ref{subsec:HeuristicCandidates}). The algorithm then repeatedly applies \textit{destroy} and \textit{local repair} to reconsider demand assignments and associated physical paths (Sections~\ref{subsec:HeuristicDestroyRepair} and~\ref{subsec:HeuristicSmallMIP}). When local improvement stagnates, \textit{path generation} uses restricted DW LPs and logical pricing to introduce additional logical paths for repair (Section~\ref{subsec:HeuristicLPSupplementation}). If several rounds of \textit{path generation} and local search yield no improvement, \textit{joint repair} combines the accumulated paths across demand groups and reallocates shared resources (Section~\ref{subsec:HeuristicFullMIP}). Both procedures return the search to \textit{destroy} and \textit{local repair}. The process terminates when the time limit is reached or a complete cycle yields no new LP or MIP solve and no new candidate path, returning the best feasible solution found.

\subsection{\textit{Construction}}
\label{subsec:HeuristicCandidates}

Let $\xi$ denote a feasible blocking plan and $z(\xi)$ its total cost. Algorithm~\ref{alg:HeuristicConstruction} summarizes the \textit{construction} procedure, starting from the all-unserved feasible plan. Each iteration processes an ordered sequence $\mathcal G$ of demand sets, beginning with individual demands and then considering groups of demands for freight consolidation (Section~\ref{subsubsec:ConstructionDemandSets}). For each demand set $\mathcal K_{\mathrm{cand}}\subseteq\mathcal K$, \textsc{Candidates} generates candidate assignments specifying the demands to adjust, their logical paths, and target served volumes (Section~\ref{subsubsec:ConstructionCandidateGeneration}). \textsc{Serve} then evaluates these candidates by determining feasible physical paths and actual served volumes, accepting the lowest-cost feasible plan if it improves the incumbent (Section~\ref{subsubsec:ConstructionJointEvaluation}). The procedure first operates over a restricted block domain and then expands the domain while retaining the constructed plan. Within each domain, iterations continue until the relative cost improvement $(z(\bar\xi)-z(\xi))/\max\{1,|z(\bar\xi)|\}\le\varepsilon_{\mathrm C}$, where $\bar\xi$ is the incumbent at the start of the iteration and $\varepsilon_{\mathrm C}\in[0,1)$ is the stopping tolerance.

\begin{algorithm}[!htbp]
\normalsize
\caption{\textit{Construction}}
\label{alg:HeuristicConstruction}
\begin{algorithmic}[1]
\setlength{\baselineskip}{12pt}
\State \(\xi\gets\) all-unserved feasible plan.
\ForAll{block domains, from restricted to expanded}
\Repeat
    \State \(\bar\xi\gets\xi\).
    \State Form the ordered demand sets \(\mathcal G\). \Comment{Section~\ref{subsubsec:ConstructionDemandSets}}
    \ForAll{\(\mathcal K_{\mathrm{cand}}\) in \(\mathcal G\), in order}
    \State \textsc{Candidates}: generate candidate assignments for \(\mathcal K_{\mathrm{cand}}\). \Comment{Section~\ref{subsubsec:ConstructionCandidateGeneration}}
    \State \textsc{Serve}: evaluate candidate assignments and update \(\xi\). \Comment{Section~\ref{subsubsec:ConstructionJointEvaluation}}
\EndFor
\Until{\((z(\bar\xi)-z(\xi))/\max\{1,|z(\bar\xi)|\}\le\varepsilon_{\mathrm C}\)}
\EndFor
\State \Return \(\xi\).
\end{algorithmic}
\end{algorithm}

\subsubsection{\textnormal{\textit{Demand Sets.}}}
\label{subsubsec:ConstructionDemandSets}

Demand grouping allows multiple demands to be considered jointly, creating opportunities for freight consolidation and satisfying minimum block loads. Accordingly, the sequence $\mathcal G$ comprises four types of demand sets, ordered by increasing complexity of their logical paths. \textit{Singleton sets} $\{k\}$ are processed first in input order. \textit{Direct groups} collect, for each candidate block $b$, demands satisfying $o_k=i_b$, $d_k=j_b$, and $b\in\mathcal B_k$. \textit{One-hub groups} collect demands that admit compatible logical paths $o_k\to h\to d_k$ through a common intermediate yard $h$, while \textit{two-hub groups} consider a common ordered pair of distinct intermediate yards $(h_1,h_2)$, following $o_k\to h_1\to h_2\to d_k$.

To limit the number of two-hub groups, we rank yards by $s_h=\sum_{k:\,h\in\{o_k,d_k\}}q_kc_k^{\mathrm U}$, which measures the total potential unserved demand cost associated with yard $h$. We form all ordered pairs of distinct yards among the $N_{\mathrm{hub}}$ highest-ranked yards, then add pairs $(i_b,j_b)$ from available blocks whose endpoints match the origin and destination of at least one demand, excluding pairs already considered.

\subsubsection{\textnormal{\textit{Candidate Assignments.}}}
\label{subsubsec:ConstructionCandidateGeneration}

During \textit{construction}, we use the \textit{estimated unit cost} (EUC) to rank candidate logical paths. We first define the block-level EUC $c_{kb}^{\mathrm C}$, which estimates the cost per unit of demand $k$ transported on block $b$. Let $\bar y_b$ denote the value of $y_b$ in the current solution, and let $r_b^{\mathrm{ref}}$ be the current physical path of an open block or a retained physical path with minimum unit cost for an unopened block. The block-level EUC is
\begin{equation}
c_{kb}^{\mathrm C}
=c_{b,r_b^{\mathrm{ref}}}^{\mathrm B}
+(1-\bar y_b)\frac{c_b^{\mathrm O}}{q_k}
+\mathbf1\{i_b\ne o_k\}c_{i_b}^{\mathrm H}.
\label{eq:ConstructionBlockUnitCost}
\end{equation}
The three terms account for transportation and interchange costs, the opening cost of an unopened block allocated over $q_k$, and intermediate-yard handling costs, respectively. Since a logical path consists of a sequence of blocks, we define its EUC as $\widehat c_k(p)=\sum_{b\in p}c_{kb}^{\mathrm C}$, which estimates the unit cost of serving demand $k$ along path $p$.

For each demand set $\mathcal K_{\mathrm{cand}}$, \textsc{Candidates} generates assignments of the form $\{(k,p_k,q_k):k\in U\}$, where $\emptyset\ne U\subseteq \mathcal K_{\mathrm{cand}}$ contains the demands to be reassigned, $p_k\in\mathcal P_k$ is the proposed logical path, and $q_k$ is the target served volume. For a singleton set $\mathcal K_{\mathrm{cand}}=\{k\}$, we enumerate compatible direct and one-hub paths, followed by minimum-EUC path searches over all available compatible blocks and, separately, over currently open blocks. For direct groups, each demand uses the designated block, while one-hub and two-hub groups follow their prescribed yard sequences using compatible blocks. Physical path candidates are generated as needed through shortest-path searches (see Section~\ref{ec:IncrementalEvaluation}).

When $\mathcal K_{\mathrm{cand}}$ is a multi-demand group, we construct candidate assignments by ranking demands according to the estimated unit savings of their proposed logical paths. Let $\bar v_k$ denote the incumbent served volume and $\bar c_k$ the unit variable cost of its incumbent logical and physical paths, excluding block opening costs, with $\bar c_k=0$ for unserved demands. The estimated unit saving is
\begin{equation}
\Delta_k=
\frac{\bar v_k \bar c_k
+(q_k-\bar v_k)c_k^{\mathrm U}}{q_k}
-\widehat c_k(p_k).
\label{eq:ConstructionUnitSaving}
\end{equation}
We order the demands $k_1,\ldots,k_{|\mathcal K_{\mathrm{cand}}|}$ by decreasing $\Delta_k$ and consider progressively larger subsets $\{k_1,k_2\},\ldots,\{k_1,\ldots,k_{|\mathcal K_{\mathrm{cand}}|}\}$, each with its proposed logical paths and target volumes. For each subset $U$, replacing the incumbent assignments with the proposed changes gives the block loads $f_b^{\mathrm{tar}}(U)=f_b(\xi)-\sum_{k\in U}\iota_{b\bar p_k}\bar v_k+\sum_{k\in U}\iota_{bp_k}q_k$, where $\bar p_k$ is the incumbent logical path of demand $k$. To satisfy minimum block loads, we iteratively update
\begin{equation}
U\gets\{k\in U:
f_b^{\mathrm{tar}}(U)=0
\text{ or }f_b^{\mathrm{tar}}(U)\ge\underline f_b,
\ \forall b\in\bar p_k\cup p_k\}.
\label{eq:ConstructionScreening}
\end{equation}
We repeat Eq.~\eqref{eq:ConstructionScreening} until $U$ remains unchanged. Each nonempty resulting subset forms a candidate assignment for \textsc{Serve}, which evaluates its feasibility and total cost.

\subsubsection{\textnormal{\textit{Assignment Evaluation.}}}
\label{subsubsec:ConstructionJointEvaluation}

For each candidate assignment $\{(k,p_k,q_k):k\in U\}$, \textsc{Serve} fixes the proposed logical paths of demands in $U$ and retains the incumbent assignments of all other demands. It also retains the incumbent physical paths of open blocks, selects feasible physical paths for newly opened blocks, and scales the target served volumes $q_k$ by a common factor $\theta\in[0,1]$. Given these path choices, the block load bounds and resource capacities determine a feasible interval $[\theta_-,\theta_+]$ for $\theta$. Since the total cost is affine in $\theta$, taking the form $z_0+\theta z_1$, the minimum-cost served volumes are
\begin{equation}
v_k=\theta^\star q_k\quad(k\in U),
\qquad
\theta^\star=
\begin{cases}
\theta_+,&z_1<0,\\
\theta_-,&z_1\ge0.
\end{cases}
\label{eq:ConstructionServiceVolumes}
\end{equation}
Using Eq.~\eqref{eq:ConstructionServiceVolumes}, \textsc{Serve} evaluates the feasible plans obtained from all candidate assignments and accepts the lowest-cost plan only if it improves the incumbent. The resulting blocking plan provides the starting point for \textit{destroy} and \textit{local repair}, described next.

\subsection{\textit{Destroy}}
\label{subsec:HeuristicDestroyRepair}

\textit{Destroy} selects a group of demands $\mathcal K_{\mathrm{free}}\subseteq\mathcal K$ for \textit{local repair}, which jointly reoptimizes their logical paths and served volumes while keeping all other demand assignments fixed. The selection proceeds in two steps. First, we choose a reference demand $k_0$. During the first outer cycle, successive calls select $k_0$ in decreasing order of its current unserved demand cost $(q_k-\bar v_k)c_k^{\mathrm U}$, while subsequent cycles select $k_0$ uniformly at random. Second, we rank demands according to their relationships with $k_0$ and their unserved demand costs using the priority score
\begin{equation}
\psi_k
=\omega_0\mathbf 1\{k=k_0\}
+\omega_1\mathbf 1\{o_k=o_{k_0}\ \mathrm{or}\ d_k=d_{k_0}\}
+\omega_2m_k
+\frac{(q_k-\bar v_k)c_k^{\mathrm U}}{\max\{1,\max_{j\in\mathcal K}(q_j-\bar v_j)c_j^{\mathrm U}\}},
\label{eq:RelatedDemandScore}
\end{equation}
where $m_k=\sum_{b\in\mathcal B_k}\bar x_{kb}^{\mathrm D}\mathbf1\{i_b=o_{k_0}\ \mathrm{or}\ j_b=d_{k_0}\}$ counts the blocks on the incumbent logical path of demand $k$ that depart from the reference origin or arrive at the reference destination. The weights $\omega_0,\omega_1,\omega_2$ set the selection priorities.

Given the current group size limit $N_{\mathrm D}\le|\mathcal K|$, we select the $N_{\mathrm D}$ demands with the highest scores in Eq.~\eqref{eq:RelatedDemandScore} to form $\mathcal K_{\mathrm{free}}$. Ties are broken by candidate order, which is randomized in subsequent outer cycles. The group size schedule and priority weights are specified in Section~\ref{ec:HeuristicSchedule} and Table~\ref{tab:HeuristicDefaults}.
\subsection{\textit{Local Repair}}
\label{subsec:HeuristicSmallMIP}

\textit{Local repair} constructs a restricted MIP of the DW formulation~\eqref{eq:DWMaster} to jointly reoptimize the logical paths and served volumes of demands in $\mathcal K_{\mathrm{free}}$, together with the physical paths of affected blocks. Assignments of other demands remain fixed.

For each demand $k\in\mathcal K_{\mathrm{free}}$, let $\widetilde{\mathcal P}_k\subseteq\mathcal P_k$ denote its candidate logical path set, including the incumbent path, if any, paths retained from \textit{path generation} and \textit{joint repair}, and low-EUC paths generated using the procedures in Section~\ref{subsubsec:ConstructionCandidateGeneration}. After consecutive calls without improvement, we increase the group size to allow more joint reassignments. We also consider bridge paths connecting chains of open blocks through additional blocks, allowing released demands to reuse existing services. The corresponding settings are given in Section~\ref{ec:HeuristicSchedule}.

The blocks whose physical paths may change are given by $\mathcal B_{\mathrm{free}}=\bigcup_{k\in\mathcal K_{\mathrm{free}}}\bigcup_{p\in\widetilde{\mathcal P}_k}\{b:\iota_{bp}=1\}$. Let $\mathcal K_{\mathrm{fixed}}=\mathcal K\setminus\mathcal K_{\mathrm{free}}$ and $\mathcal B_{\mathrm{fixed}}=\mathcal B\setminus\mathcal B_{\mathrm{free}}$ denote the remaining demands and blocks. For each block $b$, let $\widetilde{\mathcal R}_b\subseteq\mathcal R_b$ denote its retained physical path set. The path sets of fixed demands and blocks contain only their incumbent paths, when present. We substitute the fixed decisions into the DW block balance and resource capacity constraints to obtain
\begin{equation}
\begin{aligned}
\sum_{r\in\widetilde{\mathcal R}_b}v_{br}^{\mathrm B}
-\sum_{\substack{k\in\mathcal K_{\mathrm{free}}\\p\in\widetilde{\mathcal P}_k}}
\iota_{bp}v_{kp}^{\mathrm D}
&=\sum_{\substack{k\in\mathcal K_{\mathrm{fixed}}\\p\in\widetilde{\mathcal P}_k}}
\iota_{bp}\bar v_{kp}^{\mathrm D},
&&b\in\mathcal B_{\mathrm{free}},\\
\sum_{\substack{b\in\mathcal B_{\mathrm{free}}\\r\in\widetilde{\mathcal R}_b}}
\Bigl(\sum_{a\in\mathcal A_e}\delta_{ar}\Bigr)v_{br}^{\mathrm B}
&\le C_e-
\sum_{\substack{b\in\mathcal B_{\mathrm{fixed}}\\r\in\widetilde{\mathcal R}_b}}
\Bigl(\sum_{a\in\mathcal A_e}\delta_{ar}\Bigr)\bar v_{br}^{\mathrm B},
&&e\in\mathcal E,\\
\sum_{\substack{k\in\mathcal K_{\mathrm{free}}\\p\in\widetilde{\mathcal P}_k}}
\eta_{ip}v_{kp}^{\mathrm D}
&\le H_i-
\sum_{\substack{k\in\mathcal K_{\mathrm{fixed}}\\p\in\widetilde{\mathcal P}_k}}
\eta_{ip}\bar v_{kp}^{\mathrm D},
&&i\in\mathcal Y,\\
\sum_{\substack{b\in\mathcal B_{\mathrm{free}}\\i_b=i}}
\kappa_b\sum_{r\in\widetilde{\mathcal R}_b}x_{br}^{\mathrm B}
&\le T_i-
\sum_{\substack{b\in\mathcal B_{\mathrm{fixed}}\\i_b=i}}
\kappa_b\sum_{r\in\widetilde{\mathcal R}_b}\bar x_{br}^{\mathrm B},
&&i\in\mathcal Y.
\end{aligned}
\label{eq:LocalRepairCoupling}
\end{equation}
In Eq.~\eqref{eq:LocalRepairCoupling}, the physical path variables of each free block account for its full load, including contributions from fixed demands. All other applicable DW constraints, including service activation, remain enforced after substitution.

\subsection{\textit{Path Generation}}
\label{subsec:HeuristicLPSupplementation}

\textit{Path generation} identifies additional logical paths for \textit{local repair} and \textit{joint repair} using dual information from restricted DW LPs. It follows the column generation framework in Section~\ref{subsec:ExactPricing}, with heuristic logical pricing and a restricted search graph. All demands remain free in the RMP, while the retained physical path sets remain unchanged.

\subsubsection{\textnormal{\textit{Heuristic Logical Path Pricing.}}}
\label{subsubsec:HeuristicGuidancePrices}

\textit{Path generation} focuses on logical paths because additional physical paths offer limited benefit in our instances. Transportation costs are proportional to physical distance, making shortest physical paths preferable, while physical link capacities rarely become binding. Nevertheless, some physical paths may be useful when capacity constraints become binding. In this case, we generate them during \textit{joint repair}, restricting the search to open blocks to limit computational cost.

 We initialize the RMP with at most $N_{\mathrm{LP}}$ logical paths per demand, retaining the incumbent path, if any, and selecting the remaining paths in increasing EUC order, as defined in Section~\ref{subsubsec:ConstructionCandidateGeneration}. All applicable coupling constraints, including service activation, are enforced. To improve pricing efficiency, we modify the labeling algorithm in Section~\ref{subsubsec:ExactDemandPricing} in two ways. First, for labels $L_1=(i,S_1,c_1)$ and $L_2=(i,S_2,c_2)$ at the same yard, we use the heuristic dominance rule
\begin{equation}
L_1\preceq_{\mathrm H}L_2
\quad\Longleftrightarrow\quad c_1\le c_2.
\label{eq:HeuristicLogicalDominance}
\end{equation}
Unlike Proposition~\ref{prop:LogicalDominance}, Eq.~\eqref{eq:HeuristicLogicalDominance} omits the visited-set inclusion condition $S_1\subseteq S_2$. We still maintain visited-yard sets to ensure elementary paths, but the heuristic dominance rule no longer guarantees finding a minimum-JRC path.

Second, we restrict the logical pricing graph by selecting successor yards based on the retained physical paths. For demand $k$, define
\begin{equation}
d_k(i,j)=
\min_{\substack{b\in\mathcal B_k^+(i):\,j_b=j\\
r\in\widetilde{\mathcal R}_b}}
\sum_{a\in r}D_a,
\label{eq:HeuristicSuccessorDistance}
\end{equation}
with $d_k(i,j)=+\infty$ if no eligible block and physical path exist. At each yard $i$, we select up to $N_{\mathrm{nbr}}$ successor yards with the smallest distances in Eq.~\eqref{eq:HeuristicSuccessorDistance}, retaining all compatible blocks with physical path candidates leading to these yards. Compatible direct blocks and blocks on the incumbent logical path are also retained.

After each optimal RMP solve, heuristic pricing uses the current dual prices to identify logical paths with negative JRC, as defined in Eq.~\eqref{eq:LogicalFullLoadScore}. We add up to $N_{\mathrm{batch}}$ new paths across all demands, together with their selection and volume variables and the corresponding service-activation constraints, before reoptimizing the RMP. The process terminates when no new paths are found, LP improvement stagnates, or the computational limits are reached. Because pricing is heuristic, termination does not certify optimality of the complete DW relaxation.

\subsubsection{\textnormal{\textit{Logical Path Selection.}}}
\label{subsubsec:HeuristicPricingReturn}

After the final RMP solve, we consider the newly generated logical paths represented in the RMP that were absent from the shared path pool before the search. For each demand, we retain up to $N_{\mathrm{add}}$ paths with the lowest JRCs based on the final dual prices.

If computational limits terminate the search before convergence, we select paths using three criteria: JRC, resource diversity, and EUC. For each demand $k$, we first retain up to $N_{\mathrm{JRC}}$ paths with the lowest JRCs. Let $\mathcal C_k$ denote the newly generated paths in the final RMP and $\mathcal C_k^{\mathrm{sel}}\subseteq\mathcal C_k$ the paths selected so far. For each path $p$, let $\mathcal U_{\mathrm{res}}(p)$ contain its blocks, departure yards, and physical links along their reference physical paths. We measure resource diversity using the Jaccard distance~\citep{jaccard1912distribution} and select each additional path by
\begin{equation}
\begin{aligned}
d_{\mathrm{res}}(p,p')
&=1-\frac{|\mathcal U_{\mathrm{res}}(p)\cap\mathcal U_{\mathrm{res}}(p')|}
{|\mathcal U_{\mathrm{res}}(p)\cup\mathcal U_{\mathrm{res}}(p')|},\\
p^*
&\in\arg\max_{p\in\mathcal C_k\setminus\mathcal C_k^{\mathrm{sel}}}
\min_{p'\in\mathcal C_k^{\mathrm{sel}}}d_{\mathrm{res}}(p,p'),
\end{aligned}
\label{eq:GuidancePathDiversity}
\end{equation}
with lower JRC breaking ties. Using Eq.~\eqref{eq:GuidancePathDiversity}, we retain up to $N_{\mathrm{div}}$ such paths, followed by up to $N_{\mathrm{EUC}}$ additional paths in increasing EUC order, where $N_{\mathrm{JRC}}+N_{\mathrm{div}}+N_{\mathrm{EUC}}=N_{\mathrm{add}}$. The EUC search considers compatible logical paths supported by retained physical paths, including those absent from the RMP. All selected paths are added to the shared pool without removing existing paths, allowing subsequent \textit{local repair} and \textit{joint repair} calls to use them.

\subsection{\textit{Joint Repair}}
\label{subsec:HeuristicFullMIP}

\textit{Joint repair} combines logical and physical paths accumulated across different demand groups by jointly reoptimizing a larger set of demands. Before each call, we expand the retained path sets using the procedures of \textit{local repair}, including bridge paths, and generate additional physical paths for open blocks. We also consider \textit{partner chains}, which allow demands to jointly support the minimum load of an unopened block (Section~\ref{ec:PartnerChains}). To control the MIP size, at most \(N_{\mathrm J}\) demands are released. When \(|\mathcal K|\le N_{\mathrm J}\), we simply set \(\mathcal K_{\mathrm{free}}=\mathcal K\).

For larger instances, we select demands based on their unserved demand costs, potential EUC savings, and contributions to current block opening costs. Let \(\bar p_k\) denote the incumbent logical path of demand \(k\), with \(\bar p_k=\emptyset\) if unserved. For each demand with at least one retained logical path, we define
\begin{equation}
\rho_k=(q_k-\bar v_k)c_k^{\mathrm U}
+\bar v_k\max\left\{0,\widehat c_k(\bar p_k)
-\min_{p\in\widetilde{\mathcal P}_k}\widehat c_k(p)\right\}
+\sum_{b\in\bar p_k}\frac{c_b^{\mathrm O}\bar v_k}{f_b(\xi)},
\label{eq:JointRepairDemandScore}
\end{equation}
where the last two terms are zero for unserved demands.

We construct five types of candidate groups. The first three follow the grouping principles of \textit{destroy} and \textit{construction}, comprising current users of a common open block, demands whose minimum-EUC paths share a block with sufficient combined volume to meet its minimum load, and singleton demands with retained paths. The remaining two target shared resource constraints. One groups demands whose incumbent blocks use a common saturated physical link, while the other pairs a demand whose direct block is prevented by insufficient classification tracks with users of an open block whose closure could release the required tracks. For the latter, we consider the highest-scoring eligible demand at each yard and the blocks whose closure would provide sufficient track capacity.

To avoid leaving insufficient loads on blocks used by fixed demands, we complete each candidate group using a minimum-load closure. For a demand set \(U\), define \(f_b^{-}(U)=f_b(\xi)-\sum_{k\in U\cap\mathcal K_b}g_{kb}(\xi)\). Starting from \(U^{(0)}=\mathcal K_{\mathrm{cand}}\), we iteratively update
\begin{equation}
U^{(j+1)}
=U^{(j)}\cup
\bigcup_{\substack{b\in\mathcal B:\\0<f_b^{-}(U^{(j)})<\underline f_b}}
\{k\in\mathcal K_b:g_{kb}(\xi)>0\}
\label{eq:DestroyClosure}
\end{equation}
We repeat Eq.~\eqref{eq:DestroyClosure} until the set remains unchanged. Using the demand scores in Eq.~\eqref{eq:JointRepairDemandScore}, we rank the resulting groups by decreasing average score \(\rho(\mathcal K_{\mathrm{cand}})=|\mathcal K_{\mathrm{cand}}|^{-1}\sum_{k\in\mathcal K_{\mathrm{cand}}}\rho_k\), breaking ties by total score and demand index order. Starting with \(\mathcal K_{\mathrm{free}}=\emptyset\), we process groups in this order and add the closure of \(\mathcal K_{\mathrm{free}}\cup\mathcal K_{\mathrm{cand}}\) whenever its size does not exceed \(N_{\mathrm J}\); otherwise, the group is skipped.

Given \(\mathcal K_{\mathrm{free}}\), the joint MIP uses all retained logical paths of the released demands and physical paths of their associated blocks. Assignments outside the selected group remain fixed, and the model construction and solution procedures follow \textit{local repair} in Section~\ref{subsec:HeuristicSmallMIP}. When the model size permits, we additionally allow at most one released demand to select a logical path outside the retained pool through the compact constraints in Eqs.~\eqref{eq:LogicalBalance}--\eqref{eq:LogicalActivation}, using compatible blocks with retained physical paths. The free block set is expanded accordingly, and any newly selected path is retained if it improves the solution. Section~\ref{subsec:HeuristicFixedResources} describes this extension. After the joint MIP solve, the algorithm returns to \textit{destroy} and \textit{local repair}.

\section{Computational Experiments}
\label{sec:HeuristicExperiments}

Section~\ref{subsec:ExperimentProtocol} describes the common experimental setup. On the small generated instances (Section~\ref{subsec:genratedInstanceExperiments}), we compare solution quality and computational effort against exact benchmarks, evaluate the search components, and assess the robustness of DG-CGM.
Section~\ref{subsec:MediumExperiments} uses larger generated instances to study how performance varies with the number of demands and network size. On the large RAS instances (Section~\ref{subsec:RASExperiments}), we evaluate the search components at full scale and examine how demand volume, demand patterns, and resource limits affect service, consolidation, and resource use.

\subsection{Experimental Setup}
\label{subsec:ExperimentProtocol}

\subsubsection{Implementation and Computing Environment.}
All methods are implemented in C++17, compiled with Microsoft Visual C++ 14.44 (Visual Studio 2022), and run on a Windows workstation equipped with AMD EPYC 9654 processors at 2.4 GHz and 1.5 TiB of RAM. Gurobi 12.0.3 \citep{gurobi} solves the LPs and MIPs through the RouteOpt interface \citep{you2026routeopt}, and the local and joint MIPs of DG-CGM use \texttt{IntegralityFocus=1}. The ARC LP, ARC MIP, and DW LP benchmarks use eight solver threads per run, with eight pricing workers for the DW LP. Each benchmark run has a two-hour time limit. Unless stated otherwise, DG-CGM uses the default controls in Table~\ref{tab:HeuristicDefaults} and, on generated instances, one solver thread, random seed~1, and a 60-second time limit, including \textit{construction}. The six RAS scenarios use eight solver threads and a two-hour time limit, including \textit{construction}, and each scenario is run with random seeds 1, 2, and 3. All runs share the 384 logical processors of the workstation, and a new run starts only when at most 185 are busy. The ARC LP and DW LP times in Table~\ref{tab:SmallInstanceComparison} are measured separately, with at most four concurrent runs. All runtimes are reported as wall-clock times. The code, the generated instances, the configuration of every experiment, and the results of all runs are publicly available at \url{https://github.com/Zhengzhong-You/IRBP}.

ARC LP relaxes the binary variables of the ARC formulation, while ARC MIP uses the default Gurobi cut settings. We solve the LP relaxation of the DW formulation~\eqref{eq:DWMaster} using logical and physical pricing, as described in Section~\ref{subsec:ExactPricing}. Service-activation constraints~\eqref{eq:DWServiceActivation} are separated as needed, with zero dual variables assigned to constraints not yet added. All three exact benchmarks use the original instance data and do not receive solutions from DG-CGM.

\subsubsection{Instance Sets.}
For the operational study, we use data from the 2026 RAS Problem Solving Competition \citep{informs20262026}. The L1, L2, and L3 data sets share the same physical network with 47,193 nodes, 106,570 directed arcs, and identical resource capacities. L1 and L2 contain the same 2,044 demand records up to identifiers and ordering, so we select L1 and L3. L3 contains 2,427 demands, 18.7\% more than L1, and 0.3\% more total requested volume, so L1 and L3 have nearly the same volume but different demand patterns. For each data set, we apply demand multipliers \(\alpha\in\{0.5,1,2\}\) to every \(q_k\) and the corresponding block volume upper bounds derived from these quantities. Minimum block loads, resource capacities, and cost coefficients remain unchanged. Table~\ref{tab:ExperimentInstances} summarizes the instances.

The generated instances reproduce the regional connectivity, degree distribution, travel distances, yard operating characteristics, freight composition, and demand volume distributions of the RAS data. Section~\ref{ec:instance-generator} of the electronic companion describes the generation procedure. We generate 50 instances with 12 demands at each network size \(|\mathcal N|\in\{50,100,150\}\), denoted as N50, N100, and N150. To examine larger instances, we generate four additional sets following the same procedure and retain for each set the first ten valid instances.

\begin{table}[htbp]
\centering
\caption{Instance sets and their roles in the numerical study.}
\label{tab:ExperimentInstances}
\fontsize{10}{12}\selectfont\setlength{\tabcolsep}{4pt}\renewcommand{\arraystretch}{1.05}
\begin{tabular}{@{}clcccccc@{}}
\toprule
Source & Set & \#Inst. & \(|\mathcal N|\) & \(|\mathcal Y|\) & \(|\mathcal A|\) & \(|\mathcal K|\) & Role \\
\midrule
\multirow{7}{*}{Generated} & N50 & 50 & 50 & 10.4 & 153.1 & 12 & \multirow{3}{*}{\shortstack[c]{Mechanisms, solution quality,\\and computational effort}} \\
 & N100 & 50 & 100 & 18.3 & 318.4 & 12 & \\
 & N150 & 50 & 150 & 27.9 & 487.3 & 12 & \\
\cmidrule(l){2-8}
 & N150K24 & 10 & 150 & 27.4 & 485.0 & 24 & \multirow{4}{*}{Scalability} \\
 & N150K48 & 10 & 150 & 27.3 & 487.0 & 48 & \\
 & N300K48 & 10 & 300 & 55.0 & 985.0 & 48 & \\
 & N600K48 & 10 & 600 & 108.9 & 1,997.0 & 48 & \\
\midrule
\multirow{6}{*}{RAS} & L1\(\times0.5\) & 1 & 47,193 & 1,477 & 106,570 & 2,044 & Low demand \\
 & L1\(\times1\) & 1 & 47,193 & 1,477 & 106,570 & 2,044 & Base demand \\
 & L1\(\times2\) & 1 & 47,193 & 1,477 & 106,570 & 2,044 & High demand \\
\cmidrule(l){2-8}
 & L3\(\times0.5\) & 1 & 47,193 & 1,477 & 106,570 & 2,427 & Low demand \\
 & L3\(\times1\) & 1 & 47,193 & 1,477 & 106,570 & 2,427 & Base demand \\
 & L3\(\times2\) & 1 & 47,193 & 1,477 & 106,570 & 2,427 & High demand \\
\bottomrule
\end{tabular}\par\smallskip
\noindent\begin{minipage}{\linewidth}\fontsize{9}{11}\selectfont
\textit{Note.} For generated instances, \(|\mathcal Y|\) and \(|\mathcal A|\) are arithmetic means over the instances of each set.
\end{minipage}
\end{table}

\subsubsection{Performance Measures.}
For each of the 150 small generated instances, we compare the cost \(z\) of a feasible plan with a common reference cost \(z^{\mathrm{ref}}\), equal to the optimal cost when ARC MIP proves optimality and otherwise to the lowest feasible cost obtained by any ARC MIP or DG-CGM run in Section~\ref{subsec:genratedInstanceExperiments}. The relative difference \((z-z^{\mathrm{ref}})/z^{\mathrm{ref}}\) is an optimality gap when the reference is proven optimal and otherwise measures the deviation from the reference. For a valid lower bound \(\underline z\), the bound gap is \((z-\underline z)/z\). We denote the cost immediately after \textit{construction}, before any \textit{local repair}, by \(z^0\), and the final cost by \(z^{\mathrm H}\). The improvement over \textit{construction} is \((z^0-z^{\mathrm H})/z^0\). The volume-weighted service ratio is \(\sum_k v_k/\sum_k q_k\). Across instances, we report arithmetic means for costs, gaps, improvements, and service ratios, and geometric means for positive runtimes unless stated otherwise. For paired comparisons on the same instances, we report Win/Tie/Loss counts and, where indicated, the \(p\)-value from a two-sided Wilcoxon signed-rank test~\citep{wilcoxon1945individual} applied to paired cost differences. The test uses the signs and ranks of these differences to assess whether their distribution is symmetric about zero. Smaller \(p\)-values indicate stronger statistical evidence of a cost difference.
\subsection{Small Generated Instances}
\label{subsec:genratedInstanceExperiments}

\subsubsection{Bound Strength and DG-CGM Performance.}

Table~\ref{tab:SmallInstanceComparison} compares the lower bounds and computation times of the ARC and DW formulations. Across all 150 instances, the mean lower-bound gaps relative to \(z^{\mathrm{ref}}\) are 67.5\% for ARC and 3.1\% for DW, while the corresponding geometric mean computation times are 19.3 and 1.3 seconds. The substantially tighter DW bounds and shorter computation times support our use of the DW formulation to identify promising logical paths in Section~\ref{sec:PrimalHeuristic}.

\begin{table}[htbp]
\centering
\caption{LP bounds and computation times for the ARC and DW formulations.}
\label{tab:SmallInstanceComparison}
\fontsize{10}{12}\selectfont\setlength{\tabcolsep}{4pt}\renewcommand{\arraystretch}{1.05}
\begin{tabular*}{\linewidth}{@{\extracolsep{\fill}}lrrrrrrrrr@{}}
\toprule
 & & \multicolumn{3}{c}{ARC LP} & \multicolumn{3}{c}{DW LP} & \multicolumn{2}{c}{ARC MIP} \\
\cmidrule(lr){3-5}\cmidrule(lr){6-8}\cmidrule(l){9-10}
Set & \(z^{\mathrm{ref}}\)/\(10^6\) & \(\underline z_{\mathrm{arc}}\)/\(10^6\) & Gap/\% & Time/s & \(\underline z_{\mathrm{dw}}\)/\(10^6\) & Gap/\% & Time/s & Time/s & \#Opt. \\
\midrule
N50 & 17.2 & 4.9 & 69.8 & 2.1 & 16.8 & 2.6 & 0.1 & 52.1 & 50 \\
N100 & 17.5 & 6.3 & 64.6 & 23.9 & 16.9 & 3.4 & 1.5 & 2,790.1 & 32 \\
N150 & 17.3 & 5.5 & 68.0 & 146.7 & 16.7 & 3.3 & 10.0 & 6,365.9 & 9 \\
\hline
All & -- & -- & 67.5 & 19.3 & -- & 3.1 & 1.3 & 974.2 & 91 \\
\bottomrule
\end{tabular*}
\end{table}

\begin{table}[htbp]
\centering
\caption{Feasible solution quality and runtime of ARC MIP and DG-CGM.}
\label{tab:DGCGMComparison}
\fontsize{10}{12}\selectfont\setlength{\tabcolsep}{4pt}\renewcommand{\arraystretch}{1.05}
\begin{tabular*}{\linewidth}{@{\extracolsep{\fill}}llrrrrrcr@{}}
\toprule
Set & Method & \(z\)/\(10^6\) & Max. gap/\% & Mean gap/\% & \#Opt. & \shortstack[c]{Mean\\improvement/\%} & Win/Tie/Loss & Time/s \\
\midrule
N50 & ARC MIP & 17.2 & 0.0 & 0.0 & 50 & -- & -- & 52.1 \\
 & DG-CGM & 17.2 & 2.9 & 0.1 & 48 & \(-0.1\) & 0/48/2 & 36.7 \\
\addlinespace[2pt]
N100 & ARC MIP & 17.5 & 1.8 & 0.1 & 32 & -- & -- & 2,790.1 \\
 & DG-CGM & 17.5 & 1.5 & \(<0.1\) & 32 & 0.0 & 4/44/2 & 59.5 \\
\addlinespace[2pt]
N150 & ARC MIP & 19.6 & 67.9 & 13.5 & 9 & -- & -- & 6,365.9 \\
 & DG-CGM & 17.3 & 5.8 & 0.3 & 9 & 10.2 & 39/11/0 & 60.1 \\
\addlinespace[2pt] \hline
All & ARC MIP & -- & 67.9 & 4.5 & 91 & -- & -- & 974.2 \\
 & DG-CGM & -- & 5.8 & 0.1 & 89 & 3.4 & 43/103/4 & 50.8 \\
\bottomrule
\end{tabular*}
\par\smallskip
\noindent\begin{minipage}{\linewidth}\fontsize{9}{11}\selectfont
% \textit{Note.} For the N100 instance generated with seed 2026102348 (\texttt{n100-s2026102348}), Gurobi 12.0.3 incorrectly reports a solution with cost 12,625,591.91 as optimal at the root node, after exploring only one node. However, DG-CGM finds a feasible solution with cost 12,405,275.34, which is 1.7\% lower and therefore invalidates the reported optimality claim. The incorrect certification disappears when any one of the parameters is changed to \texttt{Aggregate=0}, \texttt{Presolve=0}, \texttt{Symmetry=0}, or \texttt{NumericFocus=3}. Under each of these four settings, Gurobi also finds a feasible solution with cost 12,405,275.34. These observations suggest a solver-specific issue with optimality certification under the default settings. We therefore do not count this instance as solved to proven optimality by ARC MIP.
\textit{Note.} One N100 instance is excluded from the proven-optimality count because a lower-cost feasible solution was found than Gurobi's reported optimum.
\end{minipage}
\end{table}

Table~\ref{tab:DGCGMComparison} compares the best feasible solutions obtained by ARC MIP and DG-CGM, with gaps relative to \(z^{\mathrm{ref}}\).
Among the 91 instances with proven optima, DG-CGM achieves the optimal cost on 89. Across all 150 instances, DG-CGM lowers the ARC MIP cost by 3.4\% on average while using 5.2\% of its geometric mean runtime. DG-CGM obtains lower costs on 43 instances, by up to 40.4\%, and higher costs on four, by at most 2.9\%. These four instances belong to N50 and N100, where exact optimization is most effective; on N50, ARC MIP proves optimality for every instance with a geometric mean runtime of 52.1 seconds. The two optima missed by DG-CGM are both on N50. It reaches one with a 300-second time limit but remains 0.7\% above the other across all runs and time limits, suggesting that this optimal plan uses a path outside the candidate pools of DG-CGM. The advantage of DG-CGM grows with instance size. On the largest set, N150, DG-CGM lowers the ARC MIP cost by 10.2\% on average while using 0.9\% of the ARC MIP runtime, with 39 wins, 11 ties, and no losses.

\begin{figure}[htbp]
\caption{Mean incumbent cost relative to the final DG-CGM cost over time.}
\label{fig:SmallInstanceQualityProgress}
\centering
{\normalfont\input{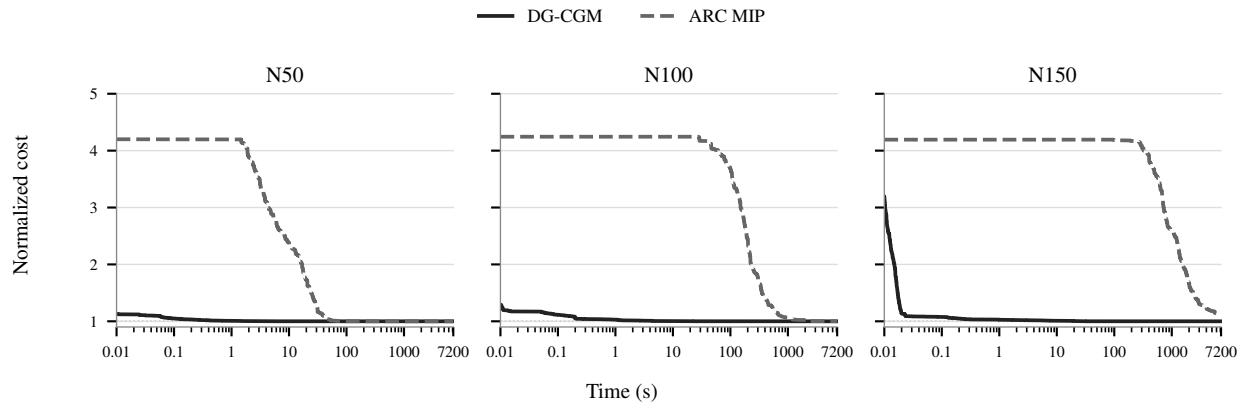}}
\end{figure}

Figure~\ref{fig:SmallInstanceQualityProgress} shows the solution trajectories of DG-CGM and ARC MIP for the three instance sets. Most of the cost reduction by DG-CGM occurs within the first 0.1 seconds, with further improvements as the search continues. After one second, the mean excess over the final DG-CGM cost is 0.7\%, 3.0\%, and 2.7\% for N50, N100, and N150, respectively. On N150, this one-second excess is smaller than the 13.1\% excess of ARC MIP at the two-hour time limit. Thus, DG-CGM obtains solutions close to its final cost within one second.

\subsubsection{Contributions of the Search Components.}
\label{subsubsec:SmallComponentStudy}

We compare the complete DG-CGM with four variants, \textit{Construction only}, \textit{Without path generation}, \textit{Without local repair}, and \textit{Without joint repair}. The \textit{Construction only} variant terminates after \textit{construction}. Each of the other three variants skips the corresponding component whenever it would otherwise be invoked. Table~\ref{tab:SmallComponents} reports the same measures as Table~\ref{tab:DGCGMComparison}, along with Win/Tie/Loss counts comparing DG-CGM with each variant.

\begin{table}[htbp]
\centering
\caption{Component comparisons on the 150 small generated instances with a 60-second time limit.}
\label{tab:SmallComponents}
\fontsize{10}{12}\selectfont\setlength{\tabcolsep}{4pt}\renewcommand{\arraystretch}{1.05}
\begin{tabular*}{\linewidth}{@{\extracolsep{\fill}}lrrrrrc@{}}
\toprule
Variant & \(z\)/\(10^6\) & Max. gap/\% & Mean gap/\% & \#Opt. & Time/s & \shortstack[c]{DG-CGM\\Win/Tie/Loss} \\
\midrule
\textit{Construction only} & 19.1 & 174.0 & 12.6 & 12 & \(<0.1\) & 135/15/0 \\
\textit{Without path generation} & 17.4 & 5.8 & 0.2 & 85 & 45.6 & 10/135/5 \\
\textit{Without local repair} & 17.3 & 5.8 & 0.1 & 89 & 46.6 & 6/139/5 \\
\textit{Without joint repair} & 17.4 & 9.0 & 0.4 & 71 & 11.7 & 42/107/1 \\
\hline
DG-CGM & 17.3 & 5.8 & 0.1 & 89 & 50.8 & -- \\
\bottomrule
\end{tabular*}
\end{table}

\textit{Joint repair} has the largest effect. With 12 demands, each \textit{joint repair} call releases all demands and solves the all-demand restricted MIP over the accumulated paths (Section~\ref{subsec:HeuristicFullMIP}). Compared with the variant without \textit{joint repair}, DG-CGM obtains lower costs on 42 instances and higher costs on one, while the variant matches 71 proven optima compared with 89 for DG-CGM (Wilcoxon \(p<0.001\)). The effects of \textit{path generation} and \textit{local repair} depend on the time limit. With a 60-second time limit, the variant without \textit{path generation} still matches 85 proven optima, with DG-CGM obtaining lower costs on ten instances and higher costs on five (\(p=0.45\)); without \textit{local repair}, the costs coincide on 139 instances (\(p=0.58\)).

\begin{figure}[htbp]
\caption{Mean gap of DG-CGM and its variants under different time limits.}
\label{fig:SmallBudgetProfile}
\centering
{\normalfont\input{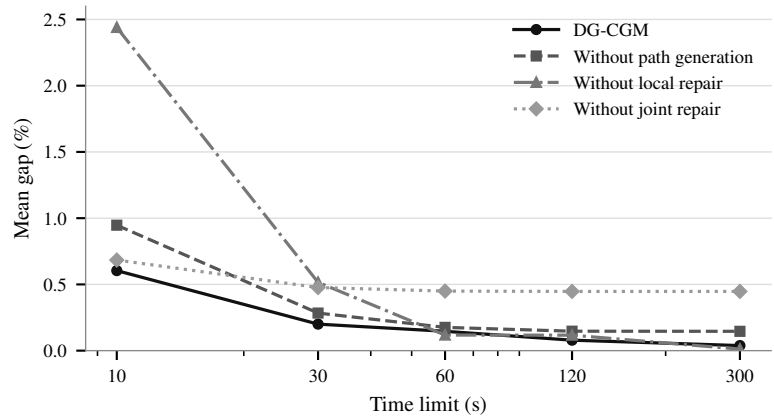}}
\end{figure}

We then repeat the comparison with time limits of 10, 30, 120, and 300 seconds, keeping all other settings fixed. Figure~\ref{fig:SmallBudgetProfile} reports the mean gap. Removing \textit{path generation} increases the mean gap at every time limit, from 0.60\% to 0.95\% at 10 seconds and from 0.04\% to 0.15\% at 300 seconds, with little further improvement after 120 seconds. Without \textit{joint repair}, the search terminates after two consecutive pricing rounds without improvement; with a 300-second time limit, 136 of the 150 runs terminate early. Consequently, the mean gap of this variant remains between 0.45\% and 0.48\% from 30 seconds onward, whereas DG-CGM continues to improve, reaching 0.04\% at 300 seconds. \textit{Local repair} drives early improvement. With a 10-second time limit, removing it increases the mean gap from 0.60\% to 2.44\%. From 60 seconds onward, this variant performs about as well as DG-CGM; with 300 seconds, it obtains lower costs on eight instances and higher costs on one, with a mean gap of 0.01\% instead of 0.04\% (Wilcoxon \(p=0.01\)). Skipping \textit{local repair} leaves more time for \textit{path generation} and \textit{joint repair}, which this variant invokes on average 54.6 and 27.2 times per run with 300 seconds, compared with 42.6 and 20.9 times for DG-CGM. Because each \textit{joint repair} call already releases all 12 demands, the neighborhoods of \textit{local repair}, with eight or twelve released demands, largely overlap with those of \textit{joint repair} when the time limit allows both to be applied. \textit{Local repair} is therefore most useful under short time limits.

\textit{Path generation} is not the only source of new paths in DG-CGM. \textit{Joint repair} also adds candidate paths before and during each call (Section~\ref{subsec:HeuristicFullMIP}), which can substitute for paths supplied by \textit{path generation}. To isolate the contribution of \textit{path generation}, we restrict \textit{joint repair} to incumbent paths and paths already in the shared pool, excluding its candidate supplements and arc options while keeping all other settings unchanged. Under this restriction, \textit{path generation} becomes decisive. Removing it increases cost by 2.5\% on average and by up to 53.7\%, with DG-CGM obtaining lower costs on 56 instances and higher costs on none (Wilcoxon \(p<0.001\)); the restricted variant matches 71 proven optima with \textit{path generation} and 52 without. Without \textit{path generation}, the restricted search terminates after 0.4 seconds (geometric mean) when local search stagnates, whereas with \textit{path generation} it recovers 97\% of the improvement achieved by the complete DG-CGM over \textit{construction}.

\subsubsection{Value of Dual Guidance.}
\label{subsubsec:CGPathGenerationValue}
We compare three \textit{path generation} procedures on all 150 instances. Each starts from the same plan immediately after \textit{construction}, before \textit{local repair}, with identical retained logical and physical paths. \textit{DG-CGM path generation} follows the procedure in Section~\ref{subsec:HeuristicLPSupplementation}.

The two shortest-path procedures independently search the full compatible logical graph supported by the retained physical paths. \textit{Transportation cost} assigns block \(b\) the weight \(c_{b,r_b^{\mathrm{ref}}}^{\mathrm B}-c_b^{\mathrm I}\), equal to the sum of transportation costs along its reference physical path. \textit{EUC} uses the block weights in Eq.~\eqref{eq:ConstructionBlockUnitCost}, which also account for interchange cost, the allocated opening cost of an unopened block, and handling cost at intermediate yards. Both procedures enumerate elementary logical paths in increasing order of total weight. All three procedures use the same reference physical paths and fixed physical path pool, retain the initial logical paths, and may add at most twelve new paths per demand.

Each procedure is allotted five seconds for generation, followed by five seconds for integer repair, including model construction, optimization, and validation. Each repair releases all twelve demands and disables the arc options of \textit{joint repair}. Generation time includes RMP preparation, restricted DW LP solves, pricing, and path selection where applicable. A diagnostic LP over each final pool, with a separate ten-second time limit, measures the LP gain; every diagnostic LP is solved to optimality.

Table~\ref{tab:CGPathGenerationValue} compares the resulting feasible costs and paired integer and LP gains relative to EUC on the same instance. Each paired gain is the cost reduction divided by the corresponding EUC cost, while Win/Tie/Loss counts compare each procedure with EUC. Reported time, an arithmetic mean, includes generation and integer repair but excludes diagnostic LP time.

Dual-guided paths yield the lowest repaired costs. DG-CGM \textit{path generation} reduces the cost of the \textit{construction} plan by 4.9\% on average, compared with 3.1\% for EUC and 1.9\% for transportation cost. Relative to EUC paths on the same instance, dual-guided paths reduce the repaired cost by 1.9\% on average, with 51 wins, 84 ties, and 15 losses, whereas paths ranked by transportation cost, which ignores opening, handling, and interchange costs, increase it by 1.6\%. The average cost reductions on winning instances also exceed the increases on losing instances, 6.4\% versus 2.8\%. The losses are consistent with differences in how the two criteria account for opening costs. JRC evaluates a path using LP dual prices, which reflect fractional block openings and can understate fixed opening costs, whereas EUC charges the full opening cost of an unopened block to the demand that would open it. The LP objective over the final dual-guided pool is also 3.2\% lower on average than that over the EUC pool.

\begin{table}[htbp]
\centering
\caption{Logical paths and integer repair on the small generated instances.}
\label{tab:CGPathGenerationValue}
\fontsize{10}{12}\selectfont\setlength{\tabcolsep}{4pt}\renewcommand{\arraystretch}{1.05}
\begin{tabular}{@{}lrrrrrr@{}}
\toprule
Procedure & \shortstack[c]{Cost/\\\(10^6\)} & \shortstack[c]{Gain from\\\(z^0\)/\%} & \shortstack[c]{Paired integer\\gain/\%} & \shortstack[c]{Paired LP\\gain/\%} & Win/Tie/Loss & Time/s \\
\midrule
DG-CGM & 18.2 & 4.9 & 1.9 & 3.2 & 51/84/15 & 0.8 \\
Transportation cost & 18.8 & 1.9 & \(-1.6\) & \(-1.9\) & 4/113/33 & 0.1 \\
EUC & 18.6 & 3.1 & 0.0 & 0.0 & -- & \(<0.1\) \\
\bottomrule
\end{tabular}
\end{table}

\subsubsection{Robustness.}
\label{subsubsec:SmallRobustness}
We evaluate the robustness of DG-CGM across five random seeds and under parameter variations based on Table~\ref{tab:HeuristicDefaults}, retaining the 60-second time limit.
Across the five runs, the mean gap ranges from 0.10\% to 0.16\%; each run matches 89 proven optima, and 143 of the 150 instances have identical costs across all runs.
The largest coefficient of variation in cost for an instance is 2.3\%. Table~\ref{tab:ParameterSensitivity} in the electronic companion reports the 14 control variants. Each variant changes the mean gap by at most 0.08 percentage points, comparable to the variation across runs, and leaves costs unchanged on 141 to 150 instances. Among the tested parameters, only the maximum number of pricing batches per \textit{path generation} invocation, \(I_{\mathrm{PG}}\), shows statistically significant differences from the default setting at the 5\% level.
Reducing \(I_{\mathrm{PG}}\) from its default value of 10 to 5 decreases the mean gap from 0.15\% to 0.09\%, while increasing it to 20 raises the gap to 0.17\%, indicating the robustness of DG-CGM.

\subsection{Larger Generated Instances}
\label{subsec:MediumExperiments}
The 12-demand instances are small enough for ARC MIP to prove optimality on many of them, whereas practical blocking problems involve many more demands. The four larger sets in Table~\ref{tab:ExperimentInstances} first increase the number of demands on 150-node networks and then increase the network size while keeping 48 demands. DG-CGM uses a 300-second time limit; we run the complete method with three different random seeds and each component variant with the default seed. Since proven optima are unavailable, we evaluate solution quality using the bound gap relative to the DW LP bound. We also run ARC LP and ARC MIP on all sets except N600K48, whose ARC formulation contains about \(2.3\times10^7\) block--arc pairs.

Table~\ref{tab:MediumScalability} summarizes the results. Where the DW bound is available, the mean gap of DG-CGM to this bound is 1.6\% on N150K24 and at most 0.6\% on N150K48 and N300K48, compared with 24.2\% to 78.8\% for ARC MIP. ARC MIP returns higher-cost solutions on 28 of the 30 instances. The larger gap on N150K24 likely reflects the bound more than the quality of the feasible plans. On the 12-demand instances, the DW bound alone lies 3.1\% below the reference cost (Table~\ref{tab:SmallInstanceComparison}), and with fewer demands per block, the LP relaxation has greater flexibility to satisfy minimum block loads through fractional block openings. DG-CGM reaches these gaps within 300 seconds, whereas exact computations take considerably longer at this scale. Within two hours, the DW LP returns a bound on 8 of the 10 N300K48 instances and on none of the N600K48 instances, while ARC LP requires 5,985.4 seconds on N300K48 (geometric mean).

The three DG-CGM runs yield nearly identical costs, with a mean coefficient of variation of at most 0.06\% in every set. The last three columns show that the contribution of \textit{joint repair} increases with instance size. Removing \textit{joint repair} raises the mean cost by 0.02\% on N150K24 and by 1.39\% on N600K48, while removing \textit{path generation} raises the mean cost by up to 0.27\% on the three sets with at most 300 nodes. Without \textit{path generation}, runs on these three sets terminate early, after 188 to 249 seconds on average, as on the RAS instances (Section~\ref{subsubsec:RASComponentStudy}). On N600K48, removing \textit{path generation} lowers the mean cost by 0.05\%, and the search then completes 7.9 instead of 5.6 joint MIP calls per run, consistent with \textit{joint repair} making the best use of the 300 seconds at this scale. Removing \textit{local repair} changes the mean cost by at most 0.08\%, consistent with each \textit{joint repair} call releasing all 24 or 48 demands, as on the 12-demand instances. The restriction of \textit{joint repair} described in Section~\ref{subsubsec:SmallComponentStudy} confirms the contribution of \textit{path generation} at this scale. Under this restriction, removing \textit{path generation} raises the mean cost by 0.61\%, 0.50\%, 0.65\%, and 0.86\% on N150K24, N150K48, N300K48, and N600K48, respectively, with 27 wins, 11 ties, and two losses across the 40 instances.

\begin{table}[htbp]
\centering
\caption{DG-CGM, DW LP, and ARC MIP on the larger generated instances.}
\label{tab:MediumScalability}
\fontsize{10}{12}\selectfont\setlength{\tabcolsep}{3pt}\renewcommand{\arraystretch}{1.05}
\begin{tabular*}{\linewidth}{@{\extracolsep{\fill}}lrrrrrrrrrr@{}}
\toprule
 & \multicolumn{2}{c}{DW LP} & \multicolumn{4}{c}{Bound gap/\%} & & \multicolumn{3}{c}{Increase without/\%} \\
\cmidrule(lr){2-3}\cmidrule(lr){4-7}\cmidrule(l){9-11}
Set & \#Bounds & Time/s & \(z^0\) & \shortstack[c]{DG-CGM\\mean} & \shortstack[c]{DG-CGM\\max.} & \shortstack[c]{ARC\\MIP} & \shortstack[c]{Improvement\\from \(z^0\)/\%} & \shortstack[c]{Path\\gen.} & Local & Joint \\
\midrule
N150K24 & 10 & 3.4 & 4.7 & 1.6 & 5.9 & 24.2 & 3.2 & 0.27 & 0.00 & 0.02 \\
N150K48 & 10 & 3.5 & 3.4 & 0.5 & 1.6 & 34.8 & 2.9 & 0.03 & 0.00 & 0.31 \\
N300K48 & 8 & 86.0 & 2.2 & 0.6 & 0.9 & 78.8 & 2.2 & 0.07 & 0.04 & 0.61 \\
N600K48 & 0 & -- & -- & -- & -- & -- & 3.9 & \(-0.05\) & 0.08 & 1.39 \\
\bottomrule
\end{tabular*}
\par\smallskip
\noindent\begin{minipage}{\linewidth}\fontsize{9}{11}\selectfont
\textit{Note.} Bound gaps and DW LP times average the instances with a DW bound. Improvement from \(z^0\) and the last three columns, the mean paired cost increase when a component is removed, use seed~1.
\end{minipage}
\end{table}

\subsection{RAS Instances}
\label{subsec:RASExperiments}
\subsubsection{\textit{Construction} and Complete Search.}
Table~\ref{tab:RASConstructionComparison} compares the cost \(z^0\) after \textit{construction} with the final cost \(z^{\mathrm H}\) obtained by DG-CGM, whose search uses the remainder of the two-hour time limit. Across the six scenarios, DG-CGM reduces the \textit{construction} cost by 7.2\% on average, with reductions ranging from 0.7\% on L3\(\times0.5\) to 26.4\% on L3\(\times2\). The \textit{construction} plan is identical across runs, and the three final costs differ by at most 0.6\% in every scenario except L3\(\times2\). On L3\(\times2\), the final costs range from 6,554.2 to 7,412.1, so even the most expensive run reduces the \textit{construction} cost by 22.0\%.

The larger spread on L3\(\times2\) is consistent with congestion. Its plans use all classification tracks at 26 yards (Table~\ref{tab:RASComparison}), so the final plan depends on which demands the search first admits to these tracks. The smallest reduction, on L3\(\times0.5\), is consistent with the design of this scenario. Halving demand volumes leaves minimum block loads unchanged, and L3 demands are already 15.5\% smaller on average than L1 demands. Consequently, many demands may be unable to reach a minimum block load even jointly; L3\(\times0.5\) keeps the lowest service of all scenarios, 95.9\%, after the search (Table~\ref{tab:RASComparison}).

\begin{table}[htbp]
\centering
\caption{\textit{Construction} and complete DG-CGM on the RAS instances.}
\label{tab:RASConstructionComparison}
\fontsize{10}{12}\selectfont\setlength{\tabcolsep}{4pt}\renewcommand{\arraystretch}{1.05}
\begin{tabular*}{\linewidth}{@{\extracolsep{\fill}}lrrrrrr@{}}
\toprule
 & \multicolumn{2}{c}{\textit{Construction}} & \multicolumn{4}{c}{DG-CGM} \\
\cmidrule(lr){2-3}\cmidrule(l){4-7}
Scenario & \(z^0\)/\(10^6\) & Time/s & \(z^{\mathrm H}\)/\(10^6\) & \shortstack[c]{\(z^{\mathrm H}\) range/\(10^6\)} & Improvement/\% & Total time/s \\
\midrule
L1\(\times0.5\) & 1,631.4 & 631.7 & 1,568.9 & 1,565.9--1,574.9 & 3.8 & 7,209.3 \\
L1\(\times1\) & 3,036.4 & 1,787.2 & 2,974.5 & 2,972.9--2,975.5 & 2.0 & 7,210.0 \\
L1\(\times2\) & 6,145.4 & 2,306.5 & 5,882.4 & 5,882.1--5,883.1 & 4.3 & 7,206.0 \\
\midrule
L3\(\times0.5\) & 1,787.6 & 1,022.1 & 1,775.7 & 1,775.1--1,776.6 & 0.7 & 7,205.3 \\
L3\(\times1\) & 3,327.3 & 1,768.0 & 3,123.9 & 3,121.5--3,126.6 & 6.1 & 7,204.9 \\
L3\(\times2\) & 9,498.2 & 3,715.3 & 6,989.7 & 6,554.2--7,412.1 & 26.4 & 7,203.7 \\
\hline
Mean & -- & -- & -- & -- & 7.2 & -- \\
\bottomrule
\end{tabular*}
\par\smallskip
\noindent\begin{minipage}{\linewidth}\fontsize{9}{11}\selectfont
\textit{Note.} \(z^{\mathrm H}\), \textit{construction} time, and total time are means over three runs; the improvement uses the mean \(z^{\mathrm H}\).
\end{minipage}
\end{table}

\subsubsection{Contributions of the Search Components.}
\label{subsubsec:RASComponentStudy}
To examine the components at full scale, we run DG-CGM and each variant from Section~\ref{subsubsec:SmallComponentStudy} starting from the \textit{construction} plan of each scenario, allowing 7,200 seconds of search after \textit{construction}.
Table~\ref{tab:RASComponents} reports the final cost of DG-CGM and the relative cost increase of each variant. \textit{Joint repair} again has the largest contribution. Without it, the cost increases by 0.4\% to 6.6\% on five scenarios and by 37.1\% on L3\(\times2\), where most of the improvement over \textit{construction} comes from serving more demand. Removing \textit{path generation} increases the cost on every scenario, by 0.3\% to 1.7\%; without it, the search terminates early, after 341 to 757 seconds, indicating that the paths generated by this component sustain the search over the two-hour time limit. The restriction of \textit{joint repair} described in Section~\ref{subsubsec:SmallComponentStudy} makes this contribution more pronounced. Under this restriction, the search without \textit{path generation} terminates after 1 to 41 seconds and changes the \textit{construction} cost by less than 0.001\% on four of the six scenarios, whereas with \textit{path generation} it runs for the full 7,200 seconds and reduces the \textit{construction} cost by 3.2\% on average. Removing \textit{path generation} in this setting increases the cost by 0.7\% to 4.3\%, exceeding the increase without the restriction on every scenario. \textit{Local repair} matters mainly on L3\(\times2\), where its removal increases the cost by 7.9\%; elsewhere, the increase is at most 0.3\%. This pattern is consistent with the initial demand selection in \textit{destroy}, which prioritizes demands in decreasing order of unserved demand cost. The \textit{construction} plan of L3\(\times2\) serves 85.7\% of the requested volume, compared with 98.8\% for L1\(\times2\) and 97.2\% for L3\(\times1\) (Table~\ref{tab:RASConstructionDiagnosis}), leaving many large unserved demands for \textit{local repair} to recover.

\begin{table}[htbp]
\centering
\caption{Search components on the RAS instances.}
\label{tab:RASComponents}
\fontsize{10}{12}\selectfont\setlength{\tabcolsep}{4pt}\renewcommand{\arraystretch}{1.05}
\begin{tabular*}{\linewidth}{@{\extracolsep{\fill}}lrrrrr@{}}
\toprule
 & & & \multicolumn{3}{c}{Cost increase without/\%} \\
\cmidrule(l){4-6}
Scenario & \(z^0\)/\(10^6\) & \(z^{\mathrm H}\)/\(10^6\) & \textit{Path generation} & \textit{Local repair} & \textit{Joint repair} \\
\midrule
L1\(\times0.5\) & 1,631.4 & 1,566.1 & 1.4 & 0.3 & 3.9 \\
L1\(\times1\) & 3,036.4 & 2,970.0 & 0.3 & 0.2 & 2.1 \\
L1\(\times2\) & 6,145.4 & 5,880.7 & 0.3 & 0.0 & 2.8 \\
\midrule
L3\(\times0.5\) & 1,787.6 & 1,776.1 & 0.4 & 0.2 & 0.4 \\
L3\(\times1\) & 3,327.3 & 3,118.1 & 0.3 & 0.2 & 6.6 \\
L3\(\times2\) & 9,498.2 & 6,429.1 & 1.7 & 7.9 & 37.1 \\
\bottomrule
\end{tabular*}
\end{table}

\subsubsection{Demand Volume and Operating Plans.}

Table~\ref{tab:RASComparison} reports total cost, cost per requested unit \(z/\sum_k q_k\), volume-weighted service, and the number of open blocks. Table~\ref{tab:RASCostComponents} reports the cost components. For each scenario, we analyze the lowest-cost plan among the three runs. We measure congestion by the number of resources operating at capacity, i.e., physical links whose loads equal their capacities, yards whose handling volumes equal their capacities, and yards with all classification tracks in use.

On L1, increasing the demand multiplier from 0.5 to 2 raises the number of open blocks from 759 to 1,249 and reduces the cost per requested unit from 1,002.8 to 941.7, consistent with consolidation over larger volumes, while service remains at or above 99.3\%. L3 opens more blocks than L1 at every multiplier and serves 95.9\% to 98.2\% of the requested volume, making unserved demand a larger cost component. At \(\alpha=2\), the cost per requested unit on L3 rises again to 1,046.1, with an unserved demand cost of 755.5 in \(10^6\) model cost units and 26 yards using all their classification tracks, compared with 28.6 and 19, respectively, on L1\(\times2\).

\begin{table}[htbp]
\centering
\caption{Cost, service, and resource use under the six RAS demand scenarios.}
\label{tab:RASComparison}
\fontsize{10}{12}\selectfont\setlength{\tabcolsep}{4pt}\renewcommand{\arraystretch}{1.05}
\begin{tabular*}{\linewidth}{@{\extracolsep{\fill}}lrrrrrrr@{}}
\toprule
 & & & & & \multicolumn{3}{c}{Resources at capacity} \\
\cmidrule(l){6-8}
Scenario & Total cost/\(10^6\) & \shortstack[c]{Cost/\\requested unit} & Service/\% & Open blocks & Links & \shortstack[c]{Yards\\(handling)} & \shortstack[c]{Yards\\(tracks)} \\
\midrule
L1\(\times0.5\) & 1,565.9 & 1,002.8 & 99.3 & 759 & 0 & 0 & 2 \\
L1\(\times1\) & 2,972.9 & 951.9 & 99.6 & 1,001 & 0 & 0 & 7 \\
L1\(\times2\) & 5,882.1 & 941.7 & 99.6 & 1,249 & 3 & 0 & 19 \\
\midrule
L3\(\times0.5\) & 1,775.1 & 1,133.3 & 95.9 & 789 & 0 & 0 & 7 \\
L3\(\times1\) & 3,121.5 & 996.4 & 98.2 & 1,109 & 0 & 0 & 12 \\
L3\(\times2\) & 6,554.2 & 1,046.1 & 96.6 & 1,293 & 0 & 2 & 26 \\
\bottomrule
\end{tabular*}
\end{table}

\begin{table}[htbp]
\centering
\caption{Cost components of the RAS plans.}
\label{tab:RASCostComponents}
\fontsize{10}{12}\selectfont\setlength{\tabcolsep}{4pt}\renewcommand{\arraystretch}{1.05}
\begin{tabular}{@{}lrrrrrr@{}}
\toprule
Scenario & Opening & Transportation & Handling & Interchange & Unserved & Total/\(10^6\) \\
\midrule
L1\(\times0.5\) & 1.9 & 1,422.0 & 40.6 & 52.7 & 48.7 & 1,565.9 \\
L1\(\times1\) & 2.5 & 2,795.0 & 46.5 & 101.9 & 27.0 & 2,972.9 \\
L1\(\times2\) & 3.1 & 5,572.4 & 75.7 & 202.3 & 28.6 & 5,882.1 \\
\midrule
L3\(\times0.5\) & 2.0 & 1,474.5 & 48.6 & 54.2 & 195.9 & 1,775.1 \\
L3\(\times1\) & 2.8 & 2,808.0 & 71.0 & 104.2 & 135.5 & 3,121.5 \\
L3\(\times2\) & 3.2 & 5,474.6 & 119.7 & 201.1 & 755.5 & 6,554.2 \\
\bottomrule
\end{tabular}
\end{table}

\subsubsection{Resource Limits.}
\label{subsubsec:RASResourceLimits}
To relate the operating plans to individual resource limits, we tighten one limit at a time on L1\(\times1\) and L3\(\times1\) and rerun DG-CGM. We multiply yard handling capacities by 0.75, reduce the number of classification tracks \(T_i\) at each yard to \(\lceil0.75T_i\rceil\), multiply physical link capacities by 0.75, or increase minimum block loads by a factor of 1.5 without exceeding the block upper bounds. We also vary the \textit{joint repair} cap from 250 to 100 and 500. Table~\ref{tab:RASSensitivity} compares each variant with the unchanged scenario under the same seed.

The minimum block load has the largest effect. Increasing it by 50\% raises cost by 3.1\% on L1\(\times1\) and by 10.4\% on L3\(\times1\), reduces the number of open blocks by 16\% to 20\%, and increases unserved demand cost by a factor of about two on L1\(\times1\) and three on L3\(\times1\). Reducing the number of classification tracks raises cost by 0.7\% on L1\(\times1\) but by 7.9\% on L3\(\times1\), where more yards already use all their classification tracks (Table~\ref{tab:RASComparison}). Lower handling and link capacities change the cost of L1\(\times1\) by at most 0.5\%; on L3\(\times1\), reducing link capacities raises the mean cost by 2.4\%, mainly due to one run. These two limits rarely bind, as no physical link reaches capacity and no yard reaches its handling capacity in the plans for L1\(\times1\) and L3\(\times1\) reported in Table~\ref{tab:RASComparison}. The 0.1\% cost decrease on L1\(\times1\) under tighter link capacities is consistent with the load limits in Section~\ref{subsec:HeuristicJRCProof}, which tighten as link capacities decrease and can strengthen the restricted MIPs, and with the different physical candidates subsequently generated by \textit{construction}.

The default \textit{joint repair} cap of 250 balances the size of each joint MIP against the time available to solve it. A cap of 100 increases cost by 0.7\% to 0.9\% in both scenarios, whereas a cap of 500 reduces the cost of L1\(\times1\) by 0.3\% but leaves L3\(\times1\) close to its \textit{construction} cost (3,324.4 versus 3,327.3 in \(10^6\) model cost units). The latter result indicates that MIPs with up to 500 released demands yield little improvement within the 30 seconds allotted to each \textit{joint repair} invocation (Table~\ref{tab:HeuristicDefaults}). L1\(\times1\), with fewer demands and fewer yards at track capacity, still benefits from the larger cap.

\begin{table}[htbp]
\centering
\caption{Resource limits and the \textit{joint repair} cap on the RAS instances.}
\label{tab:RASSensitivity}
\fontsize{10}{12}\selectfont\setlength{\tabcolsep}{3pt}\renewcommand{\arraystretch}{1.05}
\begin{tabular*}{\linewidth}{@{\extracolsep{\fill}}llrrrrrr@{}}
\toprule
Scenario & Variant & \(z\)/\(10^6\) & Change/\% & \shortstack[c]{Change by\\seed/\%} & Service/\% & \shortstack[c]{Open\\blocks} & \shortstack[c]{Unserved/\\\(10^6\)} \\
\midrule
L1\(\times1\) & Yard handling capacity \(\times0.75\) & 2,988.0 & 0.5 & 0.3 to 0.5 & 99.7 & 1,039 & 43.5 \\
 & Classification tracks \(\lceil0.75T_i\rceil\) & 2,995.2 & 0.7 & 0.6 to 0.8 & 99.5 & 971 & 36.0 \\
 & Link capacity \(\times0.75\) & 2,970.5 & \(-0.1\) & \(-0.2\) to \(-0.1\) & 99.6 & 1,002 & 24.1 \\
 & Minimum block load \(\times1.5\) & 3,067.0 & 3.1 & 2.8 to 3.5 & 99.4 & 841 & 64.1 \\
 & \textit{Joint repair} cap 100 & 2,994.1 & 0.7 & 0.5 to 0.8 & 99.6 & 1,019 & 34.3 \\
 & \textit{Joint repair} cap 500 & 2,966.5 & \(-0.3\) & \(-0.4\) to \(-0.2\) & 99.6 & 1,021 & 21.5 \\
\midrule
L3\(\times1\) & Yard handling capacity \(\times0.75\) & 3,130.8 & 0.2 & \(-0.3\) to 0.7 & 98.1 & 1,111 & 150.7 \\
 & Classification tracks \(\lceil0.75T_i\rceil\) & 3,371.6 & 7.9 & 2.8 to 17.3 & 96.0 & 992 & 444.6 \\
 & Link capacity \(\times0.75\) & 3,199.0 & 2.4 & 0.0 to 6.5 & 97.8 & 1,103 & 219.2 \\
 & Minimum block load \(\times1.5\) & 3,448.6 & 10.4 & 8.2 to 14.7 & 96.2 & 885 & 399.3 \\
 & \textit{Joint repair} cap 100 & 3,151.1 & 0.9 & 0.8 to 0.9 & 98.1 & 1,103 & 160.5 \\
 & \textit{Joint repair} cap 500 & 3,324.4 & 6.4 & 6.3 to 6.5 & 97.2 & 1,105 & 362.8 \\
\bottomrule
\end{tabular*}
\par\smallskip
\noindent\begin{minipage}{\linewidth}\fontsize{9}{11}\selectfont
\textit{Note.} Means over three runs. Change compares the mean cost with that of the unchanged scenario; change by seed reports the smallest and largest paired cost changes across the three runs.
\end{minipage}
\end{table}
\subsubsection{Service Recovery on \texorpdfstring{L3\(\times2\)}{L3 × 2}.}
\label{subsubsec:RASL3HighDemandConstruction}
The largest reduction, on L3\(\times2\), comes mainly from serving demand left unserved by the \textit{construction} plan. Table~\ref{tab:RASConstructionDiagnosis} compares the \textit{construction} plans of L1\(\times2\), L3\(\times1\), and L3\(\times2\) in terms of unserved demand cost, service, and the number of resources operating at capacity, counted as in Table~\ref{tab:RASComparison}. In the L3\(\times2\) \textit{construction} plan, ten yards operate at full handling capacity, and the unserved demand cost of 4,433.4 accounts for 46.7\% of the total cost, compared with 10.9\% for L3\(\times1\) and 6.0\% for L1\(\times2\). DG-CGM serves much of this unserved demand. The resulting L3\(\times2\) plan in Table~\ref{tab:RASCostComponents} reduces the unserved demand cost from 4,433.4 to 755.5 and has only two yards operating at handling capacity.

\begin{table}[htbp]
\centering
\caption{Service and resource use in the \textit{construction} plans.}
\label{tab:RASConstructionDiagnosis}
\fontsize{10}{12}\selectfont\setlength{\tabcolsep}{4pt}\renewcommand{\arraystretch}{1.05}
\begin{tabular*}{\linewidth}{@{\extracolsep{\fill}}lrrrrrr@{}}
\toprule
 & & & & \multicolumn{3}{c}{Resources at capacity} \\
\cmidrule(l){5-7}
Scenario & \shortstack[c]{Unserved demand\\cost/\(10^6\)} & \shortstack[c]{Share of total\\cost/\%} & Service/\% & Links & \shortstack[c]{Yards\\(handling)} & \shortstack[c]{Yards\\(tracks)} \\
\midrule
L1\(\times2\) & 366.6 & 6.0 & 98.8 & 3 & 2 & 17 \\
L3\(\times1\) & 363.8 & 10.9 & 97.2 & 0 & 0 & 15 \\
L3\(\times2\) & 4,433.4 & 46.7 & 85.7 & 3 & 10 & 26 \\
\bottomrule
\end{tabular*}
\end{table}

\section{Conclusions}
\label{sec:Conclusion}

We develop DG-CGM for the IRBP, which jointly determines block design, blocking sequences, and physical block routes under minimum-load and shared-resource constraints. We formulate the problem as ARC and derive its exact DW reformulation with a stronger LP relaxation. The proposed JRC enables exact pricing of logical and physical paths without path-specific dual variables. Building on these results, DG-CGM combines dual-guided \textit{path generation} with restricted MIP-based \textit{local and joint repair} to progressively improve feasible blocking plans.

Computational experiments demonstrate the effectiveness of the proposed approach. On 150 generated instances, the DW formulation reduces the mean lower-bound gap from 67.5\% for ARC to 3.1\%. DG-CGM matches 89 of the 91 proven optima within 60 seconds and reduces the feasible costs obtained by ARC MIP by 3.4\% on average, while requiring only 5.2\% of its geometric mean runtime. On larger generated instances, the mean gap of DG-CGM to the DW lower bound is at most 1.6\% where this bound is available. On six full-scale RAS scenarios with 47,193 physical nodes, DG-CGM reduces the \textit{construction} cost by 7.2\% on average and serves 95.9\%--99.6\% of the requested freight volume within two hours.

These findings highlight the interplay between \textit{path generation} and coordinated optimization in the IRBP. Newly generated logical paths may require multiple demands to jointly satisfy minimum block loads, while changes in demand assignments can affect block openings, physical paths, and shared resource capacities. Conversely, restricted MIP-based repair cannot explore logical paths absent from its candidate pool. Combining \textit{path generation} with \textit{local repair} and \textit{joint repair} allows the search to expand its routing alternatives and coordinate their use in feasible blocking plans.

We conclude with three directions for future research. First, learning-based methods could be incorporated into DG-CGM to improve search efficiency and adapt the algorithm to different network characteristics and demand patterns. Second, the DW formulation and exact pricing algorithms could be extended to develop a branch-and-price algorithm for solving larger IRBP instances to proven optimality. High-quality feasible solutions obtained by DG-CGM could provide initial upper bounds to accelerate the exact search. Finally, the DG-CGM framework could be extended to address other railway planning problems involving time-dependent freight demands and train schedules.

\bibliographystyle{informs2014-no-identifiers}
\bibliography{references}

\clearpage
\ECSwitch
\ECHead{Electronic Companion}
\section{Proofs}
\label{ec:Proofs}

\subsection{Proof of Theorem~\ref{thm:DantzigWolfeReformulation}}
\label{ec:DantzigWolfeReformulation}

\proof{Proof.}
We partition the constraints of ARC into three groups. The \textit{demand subsystem} of demand \(k\) consists of Eqs.~\eqref{eq:LogicalBalance}--\eqref{eq:ServedVolume} for \(k\), the bounds on \(\ell_{ki}^{\mathrm D}\) in Eq.~\eqref{eq:CompactOrderDomains}, and the binary domains of \(x_k^{\mathrm D}\) and \(x_{kb}^{\mathrm D}\). Its variables other than the order variables form the vector \(u_k^{\mathrm D}=(x_k^{\mathrm D},(x_{kb}^{\mathrm D})_{b\in\mathcal B_k},v_k,(g_{kb})_{b\in\mathcal B_k})\). The \textit{block subsystem} of block \(b\) consists of Eqs.~\eqref{eq:PhysicalBalance}--\eqref{eq:PhysicalRouteLength} and~\eqref{eq:PhysicalOrder}--\eqref{eq:BlockVolumeBound} for \(b\), the bounds on \(\ell_{bn}^{\mathrm B}\), and the binary domains of \(y_b\) and \(x_{ba}^{\mathrm B}\), with \(u_b^{\mathrm B}=(y_b,(x_{ba}^{\mathrm B})_{a\in\mathcal A},f_b,(w_{ba})_{a\in\mathcal A})\). The remaining Eqs.~\eqref{eq:LogicalActivation}, \eqref{eq:BlockVolume}, and~\eqref{eq:ArcHandlingVolume}--\eqref{eq:ArcHandlingCapacity} are the \textit{coupling constraints}. For \(p\in\mathcal P_k\) and \(r\in\mathcal R_b\), we write \(\iota_p=(\iota_{bp})_{b\in\mathcal B_k}\) and \(\delta_r=(\delta_{ar})_{a\in\mathcal A}\).

\textit{Step 1.} We show that the set \(\mathcal X_k^{\mathrm D}\) of vectors \(u_k^{\mathrm D}\) for which some order vector \((\ell_{ki}^{\mathrm D})_{i\in\mathcal Y}\) satisfies the demand subsystem of \(k\) is
\begin{equation}
\mathcal X_k^{\mathrm D}=\{0\}\cup\bigcup_{p\in\mathcal P_k}\bigl\{(1,\iota_p,t,t\iota_p):0\le t\le q_k\bigr\}.
\label{eq:DemandSubsystemSet}
\end{equation}
Let \(u_k^{\mathrm D}\in\mathcal X_k^{\mathrm D}\). If \(x_k^{\mathrm D}=0\), then, since every \(b\in\mathcal B_k\) belongs to \(\mathcal B_k^+(i_b)\), Eq.~\eqref{eq:LogicalDegrees} at \(i=i_b\) gives \(x_{kb}^{\mathrm D}=0\), and Eqs.~\eqref{eq:LogicalFlowBounds} and~\eqref{eq:ServedVolume} give \(g_{kb}=0\) and \(v_k=0\). Hence \(u_k^{\mathrm D}=0\).

Now let \(x_k^{\mathrm D}=1\), and call the blocks \(b\in\mathcal B_k\) with \(x_{kb}^{\mathrm D}=1\) \textit{selected}. By Eq.~\eqref{eq:LogicalDegrees}, at most one selected block leaves each yard and at most one enters it. Since \(o_k\ne d_k\), Eq.~\eqref{eq:LogicalBalance} states that the number of selected blocks leaving a yard minus the number entering it is one at \(o_k\), minus one at \(d_k\), and zero elsewhere. Together with these upper bounds, it implies that exactly one selected block leaves \(o_k\) and none enters it, that exactly one enters \(d_k\) and none leaves it, and that every other yard is left by a selected block if and only if it is entered by one. Starting at \(o_k\), we repeatedly follow the selected block that leaves the current yard. Every yard reached in this way other than \(d_k\) is entered by a selected block and is therefore also left by one, so the walk can stop only at \(d_k\). The walk never returns to a yard. A return to \(o_k\) would require a selected block entering \(o_k\). At the first return to another yard, that yard would be entered by two distinct selected blocks, because the blocks used for its two entries leave two different yards, as no yard is repeated before this return. Because \(\mathcal Y\) is finite, the walk reaches \(d_k\) and traces an elementary logical path \(p\in\mathcal P_k\) whose blocks are selected.

We next show that no other block is selected. Each yard of \(p\) is left by at most one selected block and entered by at most one, and these blocks lie on \(p\), because no selected block enters \(o_k\) or leaves \(d_k\). Hence every selected block off \(p\) leaves and enters yards off \(p\), and such a yard is left by a selected block off \(p\) if and only if it is entered by one. If some selected block lies off \(p\), following selected blocks off \(p\) from the tail of any of them therefore never stops, so some yard is eventually reached a second time. The blocks traversed between the two visits of the first such yard form a directed cycle \(C\) of selected blocks in which every yard is the head of exactly one block and the tail of exactly one block. Summing Eq.~\eqref{eq:LogicalOrder} over \(b\in C\), for which \(x_{kb}^{\mathrm D}=1\), gives \(\sum_{b\in C}(\ell_{kj_b}^{\mathrm D}-\ell_{ki_b}^{\mathrm D})\ge|C|\ge1\), whereas the left side equals zero. Hence the selected blocks are exactly the blocks of \(p\), that is, \((x_{kb}^{\mathrm D})_{b\in\mathcal B_k}=\iota_p\).

For the volumes, Eq.~\eqref{eq:LogicalFlowBounds} gives \(g_{kb}=0\) for \(b\notin p\). Let \(b_1,\ldots,b_m\) be the blocks of \(p\) in order, so that \(i_{b_1}=o_k\), \(j_{b_s}=i_{b_{s+1}}\) for \(s<m\), and \(j_{b_m}=d_k\). Eq.~\eqref{eq:LogicalFlowBalance} at \(o_k\) gives \(g_{kb_1}=v_k\), because \(b_1\) is the only block of \(p\) incident to \(o_k\), and at \(j_{b_s}\) with \(s<m\) it gives \(g_{kb_{s+1}}=g_{kb_s}\). Hence \(g_{kb}=v_k\iota_{bp}\) for every \(b\in\mathcal B_k\), and Eq.~\eqref{eq:ServedVolume} gives \(0\le v_k\le q_k\). Thus \(u_k^{\mathrm D}=(1,\iota_p,v_k,v_k\iota_p)\) belongs to the right side of Eq.~\eqref{eq:DemandSubsystemSet}.

Conversely, the zero vector satisfies the demand subsystem with \(\ell_{ki}^{\mathrm D}=0\) for all \(i\), since Eq.~\eqref{eq:LogicalOrder} then reads \(0\ge1-|\mathcal Y|\). For \(p\in\mathcal P_k\) with blocks \(b_1,\ldots,b_m\) and \(t\in[0,q_k]\), set \(\ell_{ko_k}^{\mathrm D}=0\), \(\ell_{kj_{b_s}}^{\mathrm D}=s\) for \(s=1,\ldots,m\), and \(\ell_{ki}^{\mathrm D}=0\) at the yards off \(p\). These positions lie in \([0,|\mathcal Y|-1]\) because \(p\) visits \(m+1\le|\mathcal Y|\) distinct yards. Eq.~\eqref{eq:LogicalOrder} holds with equality for the blocks of \(p\), and for every other block it reads \(\ell_{kj_b}^{\mathrm D}-\ell_{ki_b}^{\mathrm D}\ge1-|\mathcal Y|\), which holds because both positions lie in \([0,|\mathcal Y|-1]\). Eqs.~\eqref{eq:LogicalBalance} and~\eqref{eq:LogicalDegrees} hold because \(p\) is an elementary path from \(o_k\) to \(d_k\), and Eqs.~\eqref{eq:LogicalFlowBalance}, \eqref{eq:LogicalFlowBounds}, and~\eqref{eq:ServedVolume} hold for \((1,\iota_p,t,t\iota_p)\) because \(0\le t\le q_k\). The two inclusions prove Eq.~\eqref{eq:DemandSubsystemSet}.

\textit{Step 2.} Similarly, the set \(\mathcal X_b^{\mathrm B}\) of vectors \(u_b^{\mathrm B}\) for which some order vector \((\ell_{bn}^{\mathrm B})_{n\in\mathcal N}\) satisfies the block subsystem of \(b\) is
\begin{equation}
\mathcal X_b^{\mathrm B}=\{0\}\cup\bigcup_{r\in\mathcal R_b}\bigl\{(1,\delta_r,t,t\delta_r):\underline f_b\le t\le\overline f_b\bigr\}.
\label{eq:BlockSubsystemSet}
\end{equation}
If \(y_b=0\), Eq.~\eqref{eq:PhysicalDegrees} at \(n=i_a\) gives \(x_{ba}^{\mathrm B}=0\) for every arc \(a\), Eq.~\eqref{eq:BlockVolumeBound} gives \(f_b=0\), and Eq.~\eqref{eq:PhysicalFlowUpper} gives \(w_{ba}=0\). If \(y_b=1\), Eq.~\eqref{eq:BlockVolumeBound} gives \(\underline f_b\le f_b\le\overline f_b\), so \(\overline f_b\ge\underline f_b\). The argument of Step 1 applies to the arcs \(a\) with \(x_{ba}^{\mathrm B}=1\), with nodes in place of yards, Eqs.~\eqref{eq:PhysicalBalance}, \eqref{eq:PhysicalDegrees}, and~\eqref{eq:PhysicalOrder} in place of Eqs.~\eqref{eq:LogicalBalance}, \eqref{eq:LogicalDegrees}, and~\eqref{eq:LogicalOrder}, and \(i_b\ne j_b\) in place of \(o_k\ne d_k\). Hence the selected arcs form an elementary physical path \(r\) from \(i_b\) to \(j_b\), and \((x_{ba}^{\mathrm B})_{a\in\mathcal A}=\delta_r\). Eq.~\eqref{eq:PhysicalRouteLength} gives \(\sum_{a\in\mathcal A}D_a\delta_{ar}\le\overline D_b\), so \(r\in\mathcal R_b\) because \(\overline f_b\ge\underline f_b\). For \(a\notin r\), Eq.~\eqref{eq:PhysicalFlowUpper} gives \(0\le w_{ba}\le0\). For \(a\in r\), it gives \(w_{ba}\le f_b\), while Eq.~\eqref{eq:PhysicalFlowLower} gives \(w_{ba}\ge f_b-\overline f_b(1-1)=f_b\). Thus \(u_b^{\mathrm B}=(1,\delta_r,f_b,f_b\delta_r)\). Conversely, the zero vector satisfies the block subsystem with zero positions, since Eq.~\eqref{eq:PhysicalOrder} then reads \(0\ge1-|\mathcal N|\) and Eq.~\eqref{eq:PhysicalFlowLower} reads \(0\ge0\). For \(r\in\mathcal R_b\) and \(t\in[\underline f_b,\overline f_b]\), positions \(0,1,\ldots\) along \(r\) and zero elsewhere lie in \([0,|\mathcal N|-1]\) and satisfy Eq.~\eqref{eq:PhysicalOrder} as in Step 1, Eqs.~\eqref{eq:PhysicalBalance} and~\eqref{eq:PhysicalDegrees} hold because \(r\) is an elementary path from \(i_b\) to \(j_b\), Eq.~\eqref{eq:PhysicalRouteLength} holds because \(r\in\mathcal R_b\), and Eqs.~\eqref{eq:PhysicalFlowUpper}--\eqref{eq:BlockVolumeBound} hold for \((1,\delta_r,t,t\delta_r)\); in particular, Eq.~\eqref{eq:PhysicalFlowLower} reads \(t\ge t\) for \(a\in r\) and \(0\ge t-\overline f_b\) for \(a\notin r\). When \(\overline f_b<\underline f_b\), Eq.~\eqref{eq:BlockVolumeBound} excludes \(y_b=1\), in agreement with \(\mathcal R_b=\emptyset\).

\textit{Step 3.} We relate the path variables of DW to these sets. Let \(\xi=(x^{\mathrm D},x^{\mathrm B},y,v,f,g,w,v^{\mathrm H})\) collect the ARC variables other than the order variables. The linear map \(\Phi\) from the DW variables \((x,v)\) to \(\xi\) is
\begin{equation}
\begin{aligned}
x_k^{\mathrm D}&=\sum_{p\in\mathcal P_k}x_{kp}^{\mathrm D},
&
y_b&=\sum_{r\in\mathcal R_b}x_{br}^{\mathrm B},\\
x_{kb}^{\mathrm D}&=\sum_{p\in\mathcal P_k}\iota_{bp}x_{kp}^{\mathrm D},
&
x_{ba}^{\mathrm B}&=\sum_{r\in\mathcal R_b}\delta_{ar}x_{br}^{\mathrm B},\\
v_k&=\sum_{p\in\mathcal P_k}v_{kp}^{\mathrm D},
&
f_b&=\sum_{r\in\mathcal R_b}v_{br}^{\mathrm B},\\
g_{kb}&=\sum_{p\in\mathcal P_k}\iota_{bp}v_{kp}^{\mathrm D},
&
w_{ba}&=\sum_{r\in\mathcal R_b}\delta_{ar}v_{br}^{\mathrm B},\\
v_i^{\mathrm H}&=\sum_{k\in\mathcal K}\sum_{p\in\mathcal P_k}\eta_{ip}v_{kp}^{\mathrm D}.
\end{aligned}
\label{eq:CompactRecovery}
\end{equation}
Let \(\Phi_k^{\mathrm D}(x,v)=\sum_{p\in\mathcal P_k}(x_{kp}^{\mathrm D},x_{kp}^{\mathrm D}\iota_p,v_{kp}^{\mathrm D},v_{kp}^{\mathrm D}\iota_p)\) and \(\Phi_b^{\mathrm B}(x,v)=\sum_{r\in\mathcal R_b}(x_{br}^{\mathrm B},x_{br}^{\mathrm B}\delta_r,v_{br}^{\mathrm B},v_{br}^{\mathrm B}\delta_r)\) denote the components of \(\Phi(x,v)\) that form \(u_k^{\mathrm D}\) and \(u_b^{\mathrm B}\). The \textit{local constraints} of demand \(k\) in the DW LP relaxation are Eqs.~\eqref{eq:DWDemandVolume} and~\eqref{eq:DWDemandConvexity} for \(k\) together with \(x_{kp}^{\mathrm D}\ge0\), that is, \(x_{kp}^{\mathrm D}\ge0\), \(0\le v_{kp}^{\mathrm D}\le q_kx_{kp}^{\mathrm D}\), and \(\sum_px_{kp}^{\mathrm D}\le1\). Those of block \(b\), from Eqs.~\eqref{eq:DWBlockPathVolume} and~\eqref{eq:DWBlockConvexity}, are \(x_{br}^{\mathrm B}\ge0\), \(\underline f_bx_{br}^{\mathrm B}\le v_{br}^{\mathrm B}\le\overline f_bx_{br}^{\mathrm B}\), and \(\sum_rx_{br}^{\mathrm B}\le1\). We show two properties.

\textit{(a)} If the local constraints of \(k\) hold and every \(x_{kp}^{\mathrm D}\) is binary, then \(\Phi_k^{\mathrm D}(x,v)\in\mathcal X_k^{\mathrm D}\), and every point of \(\mathcal X_k^{\mathrm D}\) arises in this way. The same holds for blocks. The convexity constraint allows at most one path \(p\) with \(x_{kp}^{\mathrm D}=1\), and \(v_{kp'}^{\mathrm D}\le q_kx_{kp'}^{\mathrm D}=0\) for every other path \(p'\). Hence \(\Phi_k^{\mathrm D}(x,v)\) is either zero or \((1,\iota_p,v_{kp}^{\mathrm D},v_{kp}^{\mathrm D}\iota_p)\) with \(0\le v_{kp}^{\mathrm D}\le q_k\), which belongs to \(\mathcal X_k^{\mathrm D}\) by Eq.~\eqref{eq:DemandSubsystemSet}. Conversely, zero arises from zero path variables, and \((1,\iota_p,t,t\iota_p)\) arises from \(x_{kp}^{\mathrm D}=1\), \(v_{kp}^{\mathrm D}=t\), and zero for all other paths. For a block, the same argument uses \(\underline f_b\le v_{br}^{\mathrm B}\le\overline f_b\) when \(x_{br}^{\mathrm B}=1\) and Eq.~\eqref{eq:BlockSubsystemSet}; when \(\mathcal R_b=\emptyset\), both \(\Phi_b^{\mathrm B}(x,v)\) and the only point of \(\mathcal X_b^{\mathrm B}\) are zero.

\textit{(b)} If the local constraints of \(k\) hold, then \(\Phi_k^{\mathrm D}(x,v)\in\operatorname{conv}\mathcal X_k^{\mathrm D}\), and every point of \(\operatorname{conv}\mathcal X_k^{\mathrm D}\) arises in this way. The same holds for blocks. For the first part, under the local constraints,
\[
\Phi_k^{\mathrm D}(x,v)=\Bigl(1-\sum_{p\in\mathcal P_k}x_{kp}^{\mathrm D}\Bigr)\cdot0+\sum_{p\in\mathcal P_k}\Bigl[\Bigl(x_{kp}^{\mathrm D}-\frac{v_{kp}^{\mathrm D}}{q_k}\Bigr)(1,\iota_p,0,0)+\frac{v_{kp}^{\mathrm D}}{q_k}(1,\iota_p,q_k,q_k\iota_p)\Bigr],
\]
where the weights are nonnegative by the local constraints and sum to one, and the points belong to \(\mathcal X_k^{\mathrm D}\). For a block with \(\underline f_b<\overline f_b\), the corresponding representation places weight \((\overline f_bx_{br}^{\mathrm B}-v_{br}^{\mathrm B})/(\overline f_b-\underline f_b)\) on \((1,\delta_r,\underline f_b,\underline f_b\delta_r)\), weight \((v_{br}^{\mathrm B}-\underline f_bx_{br}^{\mathrm B})/(\overline f_b-\underline f_b)\) on \((1,\delta_r,\overline f_b,\overline f_b\delta_r)\), and weight \(1-\sum_rx_{br}^{\mathrm B}\) on zero. These weights are nonnegative, and those of path \(r\) sum to \(x_{br}^{\mathrm B}\) and reproduce the volume \(v_{br}^{\mathrm B}\). When \(\underline f_b=\overline f_b\), the local constraints give \(v_{br}^{\mathrm B}=\underline f_bx_{br}^{\mathrm B}\), and the weight \(x_{br}^{\mathrm B}\) on \((1,\delta_r,\underline f_b,\underline f_b\delta_r)\) suffices. When \(\overline f_b<\underline f_b\), \(\mathcal R_b=\emptyset\) and \(\Phi_b^{\mathrm B}(x,v)=0\in\mathcal X_b^{\mathrm B}\). For the second part, the path variables that satisfy the local constraints form a convex set and \(\Phi_k^{\mathrm D}\) is linear, so their images form a convex set, which contains \(\mathcal X_k^{\mathrm D}\) by (a) and hence contains \(\operatorname{conv}\mathcal X_k^{\mathrm D}\); the same argument applies to blocks.

\textit{Step 4.} We show that, for \(\xi=\Phi(x,v)\), the coupling constraints and the objective of ARC coincide with the remaining constraints and the objective of DW. Eq.~\eqref{eq:BlockVolume} becomes \(\sum_rv_{br}^{\mathrm B}=\sum_{k\in\mathcal K_b}\sum_p\iota_{bp}v_{kp}^{\mathrm D}\), which is Eq.~\eqref{eq:DWBlockBalance}. Eq.~\eqref{eq:LogicalActivation} becomes \(\sum_p\iota_{bp}x_{kp}^{\mathrm D}\le\sum_rx_{br}^{\mathrm B}\), which is Eq.~\eqref{eq:DWServiceActivation}. Eq.~\eqref{eq:ArcTrackCapacity} becomes Eq.~\eqref{eq:DWTrackCapacity}, and, since \(\sum_{a\in\mathcal A_e}w_{ba}=\sum_r(\sum_{a\in\mathcal A_e}\delta_{ar})v_{br}^{\mathrm B}\), Eq.~\eqref{eq:ArcLinkCapacity} becomes Eq.~\eqref{eq:DWLinkCapacity}. For handling, an elementary path \(p\in\mathcal P_k\) leaves each yard at most once, so \(\sum_{b\in\mathcal B_k^+(i)}\iota_{bp}\) equals one if \(p\) leaves yard \(i\) and zero otherwise. Since \(o_k\ne d_k\), every \(p\in\mathcal P_k\) leaves \(o_k\), and therefore \(\sum_{b\in\mathcal B_k^+(i)}\iota_{bp}-\mathbf1\{i=o_k\}=\eta_{ip}\) for every yard \(i\). Substituting \(g_{kb}=\sum_p\iota_{bp}v_{kp}^{\mathrm D}\) and \(v_k=\sum_pv_{kp}^{\mathrm D}\) into the right side of Eq.~\eqref{eq:ArcHandlingVolume} thus gives \(\sum_k\sum_p\eta_{ip}v_{kp}^{\mathrm D}\), the value of \(v_i^{\mathrm H}\) in Eq.~\eqref{eq:CompactRecovery}. Hence Eq.~\eqref{eq:ArcHandlingVolume} holds, and Eq.~\eqref{eq:ArcHandlingCapacity} becomes Eq.~\eqref{eq:DWHandlingCapacity}. Finally, substituting Eq.~\eqref{eq:CompactRecovery} into Objective~\eqref{eq:CompactObjective} gives
\[
\begin{aligned}
&\sum_bc_b^{\mathrm O}\sum_rx_{br}^{\mathrm B}+\sum_{b,r}\Bigl(\sum_ac_{ba}^{\mathrm T}\delta_{ar}\Bigr)v_{br}^{\mathrm B}+\sum_{k,p}\Bigl(\sum_ic_i^{\mathrm H}\eta_{ip}\Bigr)v_{kp}^{\mathrm D}+\sum_bc_b^{\mathrm I}\sum_rv_{br}^{\mathrm B}+\sum_kc_k^{\mathrm U}\Bigl(q_k-\sum_pv_{kp}^{\mathrm D}\Bigr)\\
&\qquad=\sum_kq_kc_k^{\mathrm U}+\sum_{k,p}c_{kp}^{\mathrm D}v_{kp}^{\mathrm D}+\sum_{b,r}\bigl(c_b^{\mathrm O}x_{br}^{\mathrm B}+c_{br}^{\mathrm B}v_{br}^{\mathrm B}\bigr)
\end{aligned}
\]
by the definitions of \(c_{kp}^{\mathrm D}\) and \(c_{br}^{\mathrm B}\) in Section~\ref{subsec:MasterReformulation}, and the right side is Objective~\eqref{eq:DWObjective}.

\textit{Step 5.} We prove \(z_{\mathrm{arc}}^*=z_{\mathrm{dw}}^*\). Let \((x,v)\) be feasible for DW. By Step 3(a), \(\Phi_k^{\mathrm D}(x,v)\in\mathcal X_k^{\mathrm D}\) and \(\Phi_b^{\mathrm B}(x,v)\in\mathcal X_b^{\mathrm B}\), so order vectors exist that satisfy every demand and block subsystem. By Step 4, \(\Phi(x,v)\) satisfies the coupling constraints and has the DW objective value. Hence \(\Phi(x,v)\), together with these order vectors, is feasible for ARC with the same cost, and \(z_{\mathrm{arc}}^*\le z_{\mathrm{dw}}^*\). Conversely, let \(\xi\) and its order vectors be feasible for ARC. By Eqs.~\eqref{eq:DemandSubsystemSet} and~\eqref{eq:BlockSubsystemSet}, each \(u_k^{\mathrm D}\) is zero or equals \((1,\iota_p,t,t\iota_p)\) for a path \(p\in\mathcal P_k\), which is unique because it consists of the blocks with \(x_{kb}^{\mathrm D}=1\), and each \(u_b^{\mathrm B}\) is zero or equals \((1,\delta_r,t,t\delta_r)\) for a unique \(r\in\mathcal R_b\). Setting the selection and volume variables of these paths to one and \(t\), and all other path variables to zero, gives binary path variables that satisfy the local constraints and whose image under \(\Phi\) is \(\xi\), where \(v^{\mathrm H}\) agrees by Step 4. By Step 4, these path variables also satisfy the remaining DW constraints and have the cost of \(\xi\), so \(z_{\mathrm{dw}}^*\le z_{\mathrm{arc}}^*\). Therefore \(z_{\mathrm{arc}}^*=z_{\mathrm{dw}}^*\).

Steps 3 and 4 also identify DW as a Dantzig--Wolfe reformulation of ARC. Under the linear map \(\Phi\), the local constraints of the DW LP relaxation describe exactly \(\operatorname{conv}\mathcal X_k^{\mathrm D}\) and \(\operatorname{conv}\mathcal X_b^{\mathrm B}\) by Step 3(b), and the coupling constraints and objective of ARC become the remaining constraints and objective of DW by Step 4. The DW LP relaxation is therefore the master problem obtained by convexifying every demand and block subsystem while retaining the coupling constraints, and requiring binary path selections recovers exactly \(\mathcal X_k^{\mathrm D}\) and \(\mathcal X_b^{\mathrm B}\) by Step 3(a), which makes the reformulation exact. Formally, \(\Phi\) maps the feasible set of DW onto the set of vectors \(\xi\) that, together with some order vectors, are feasible for ARC, and it preserves the objective value.

\textit{Step 6.} We prove \(z_{\mathrm{arc}}^{\mathrm{lp}}\le z_{\mathrm{dw}}^{\mathrm{lp}}\). Let \((x,v)\) be feasible for the DW LP relaxation. By Step 3(b), \(\Phi_k^{\mathrm D}(x,v)\) is a convex combination of points of \(\mathcal X_k^{\mathrm D}\), and by Step 1, each of these points satisfies the demand subsystem together with some order vector. Once the binary domains are relaxed to \([0,1]\), every constraint of the demand subsystem is linear, so the same convex combination of these points and their order vectors satisfies the LP relaxation of the demand subsystem, and its components other than the order vector equal \(\Phi_k^{\mathrm D}(x,v)\). The same argument applies to every block. By Step 4, \(\Phi(x,v)\) also satisfies the coupling constraints and has the DW objective value. Hence every feasible solution of the DW LP relaxation maps to a feasible solution of the ARC LP relaxation with the same objective value, and \(z_{\mathrm{arc}}^{\mathrm{lp}}\le z_{\mathrm{dw}}^{\mathrm{lp}}\).

\textit{Step 7.} We show that the inequality can be strict. Consider the instance with yards \(\mathcal Y=\{o,d\}\), physical nodes \(\mathcal N=\{o,n_1,n_2,d\}\), and four one-way links of capacity 1, each consisting of a single arc. Arcs \((o,n_1)\) and \((n_1,d)\) have length 1 and unit transportation cost 5, and arcs \((o,n_2)\) and \((n_2,d)\) have length 2 and unit transportation cost 1. There is one block \(b\) from \(o\) to \(d\) with \(\overline D_b=3\), \(\underline f_b=1\), \(\kappa_b=1\), \(c_b^{\mathrm O}=1\), and \(c_b^{\mathrm I}=0\), and one demand \(k\) from \(o\) to \(d\) with \(q_k=1\), \(\mathcal B_k=\{b\}\), and \(c_k^{\mathrm U}=100\), so that \(\overline f_b=q_k=1\). Both yards have \(T_i=H_i=1\) and \(c_i^{\mathrm H}=0\). The path \(o\to n_1\to d\) has length 2 and the path \(o\to n_2\to d\) has length 4, so \(\mathcal R_b\) contains only the first.

Consider any feasible solution of the ARC LP relaxation. Since \(\mathcal B_k=\{b\}\), Eq.~\eqref{eq:LogicalFlowBalance} at \(o\) gives \(g_{kb}=v_k\), Eq.~\eqref{eq:BlockVolume} gives \(f_b=g_{kb}\), and Eq.~\eqref{eq:BlockVolumeBound} with \(\underline f_b=\overline f_b=1\) gives \(f_b=y_b\), so \(v_k=y_b\). For every arc \(a\), Eq.~\eqref{eq:PhysicalFlowUpper} gives \(w_{ba}\le x_{ba}^{\mathrm B}\), and Eq.~\eqref{eq:PhysicalFlowLower} gives \(w_{ba}\ge f_b-(y_b-x_{ba}^{\mathrm B})=x_{ba}^{\mathrm B}\), so \(w_{ba}=x_{ba}^{\mathrm B}\). Eq.~\eqref{eq:PhysicalBalance} at \(n_1\), \(n_2\), and \(o\) gives \(x_{b,(o,n_1)}^{\mathrm B}=x_{b,(n_1,d)}^{\mathrm B}=:\chi_1\), \(x_{b,(o,n_2)}^{\mathrm B}=x_{b,(n_2,d)}^{\mathrm B}=:\chi_2\), and \(\chi_1+\chi_2=y_b\). Eq.~\eqref{eq:PhysicalRouteLength} gives \(2\chi_1+4\chi_2\le3y_b\), and substituting \(\chi_2=y_b-\chi_1\) yields \(\chi_1\ge y_b/2\). Since \(c_i^{\mathrm H}=c_b^{\mathrm I}=0\), the objective equals \(y_b+10\chi_1+2\chi_2+100(1-v_k)=100-97y_b+8\chi_1\ge100-93y_b\ge7\), where the last step uses \(y_b\le1\). The value 7 is attained by \(x_k^{\mathrm D}=x_{kb}^{\mathrm D}=y_b=v_k=g_{kb}=f_b=1\) and \(x_{ba}^{\mathrm B}=w_{ba}=1/2\) on all four arcs, with logical positions \(\ell_{ko}^{\mathrm D}=0\) and \(\ell_{kd}^{\mathrm D}=1\), physical positions \(\ell_{bn}^{\mathrm B}=0\) for every node \(n\), and \(v_o^{\mathrm H}=v_d^{\mathrm H}=0\). For this point, Eqs.~\eqref{eq:PhysicalBalance} and~\eqref{eq:PhysicalDegrees} hold because each of the two paths carries \(1/2\), Eq.~\eqref{eq:PhysicalRouteLength} holds because \((1+1)/2+(2+2)/2=3\), Eq.~\eqref{eq:PhysicalOrder} holds because its right side equals \(1-4(1-1/2)=-1\), Eq.~\eqref{eq:PhysicalFlowUpper} holds because \(w_{ba}=1/2=\overline f_bx_{ba}^{\mathrm B}\le f_b\), and Eq.~\eqref{eq:PhysicalFlowLower} holds because its right side equals \(1-(1-1/2)=1/2\). Eqs.~\eqref{eq:LogicalBalance}--\eqref{eq:LogicalActivation} hold because the single block carries the whole demand from \(o\) to \(d\) with positions 0 and 1, Eqs.~\eqref{eq:BlockVolume} and~\eqref{eq:BlockVolumeBound} hold because \(f_b=g_{kb}=1\), each link carries \(1/2\le1\), the track constraint at \(o\) reads \(1\le1\), and Eq.~\eqref{eq:ArcHandlingVolume} gives \(v_o^{\mathrm H}=g_{kb}-v_k=0\) and \(v_d^{\mathrm H}=0\), within \(H_o=H_d=1\). All positions lie within the bounds in Eq.~\eqref{eq:CompactOrderDomains}, since \(|\mathcal Y|-1=1\) and \(|\mathcal N|-1=3\), and all relaxed binary variables lie in \([0,1]\). Hence \(z_{\mathrm{arc}}^{\mathrm{lp}}=7\).

In the DW LP relaxation, \(\mathcal P_k\) contains the single path \(p=(b)\) and \(\mathcal R_b\) contains the single path \(r=o\to n_1\to d\), with \(c_{kp}^{\mathrm D}=-100\) and \(c_{br}^{\mathrm B}=10\). Eq.~\eqref{eq:DWBlockPathVolume} with \(\underline f_b=\overline f_b=1\) gives \(v_{br}^{\mathrm B}=x_{br}^{\mathrm B}\), and Eq.~\eqref{eq:DWBlockBalance} gives \(v_{kp}^{\mathrm D}=v_{br}^{\mathrm B}\). The objective therefore equals \(100-100x_{br}^{\mathrm B}+x_{br}^{\mathrm B}+10x_{br}^{\mathrm B}=100-89x_{br}^{\mathrm B}\ge11\), because \(x_{br}^{\mathrm B}\le1\) by Eq.~\eqref{eq:DWBlockConvexity}. Setting all four path variables to one attains 11 and satisfies every DW constraint. Eqs.~\eqref{eq:DWBlockBalance} and~\eqref{eq:DWDemandVolume}--\eqref{eq:DWServiceActivation} hold with all four values equal to one, the two links of \(r\) carry \(1\le1\) and the other two links carry \(0\le1\), the track constraint at \(o\) reads \(1\le1\), and the handling constraints read \(0\le1\) because \(\eta_{ip}=0\) for the direct path \(p\). Hence \(z_{\mathrm{dw}}^{\mathrm{lp}}=11>7=z_{\mathrm{arc}}^{\mathrm{lp}}\). \Halmos
\endproof

\subsection{Proofs of Propositions~\ref{prop:LogicalDominance} and~\ref{prop:PhysicalDominance}}
\label{ec:LogicalDominance}
\label{ec:PhysicalDominance}

\proof{Proof.}
We first verify that label costs equal path JRCs. For \(p\in\mathcal P_k\), substituting \(c_{kp}^{\mathrm D}=\sum_ic_i^{\mathrm H}\eta_{ip}-c_k^{\mathrm U}\) into Eqs.~\eqref{eq:LogicalAdjustedCost} and~\eqref{eq:LogicalFullLoadScore} gives
\[
\widetilde c_{kp}^{\mathrm D}=-q_kc_k^{\mathrm U}-\pi_k+\sum_{b\in\mathcal B_k}\iota_{bp}\bigl(q_k\sigma_b-\gamma_{kb}\bigr)+q_k\sum_{i\in\mathcal Y}\bigl(c_i^{\mathrm H}-\mu_i^{\mathrm H}\bigr)\eta_{ip}.
\]
Each block \(b\in\mathcal B_k\) belongs to exactly one set \(\mathcal B_k^+(i)\), the one with \(i=i_b\), so the definition of \(\eta_{ip}\) gives \(\sum_i(c_i^{\mathrm H}-\mu_i^{\mathrm H})\eta_{ip}=\sum_b\iota_{bp}\mathbf1\{i_b\ne o_k\}(c_{i_b}^{\mathrm H}-\mu_{i_b}^{\mathrm H})\). Hence \(\widetilde c_{kp}^{\mathrm D}=\widetilde c_k^{\mathrm D,\mathrm{init}}+\sum_b\iota_{bp}\widetilde c_{kb}^{\mathrm D}\), with the extension costs in Eq.~\eqref{eq:LogicalExtensionCost}. Similarly, substituting \(c_{br}^{\mathrm B}=c_b^{\mathrm I}+\sum_ac_{ba}^{\mathrm T}\delta_{ar}\) into Eqs.~\eqref{eq:PhysicalReducedCost} and~\eqref{eq:PhysicalFullLoadScore} gives \(\widetilde c_{br}^{\mathrm B}(t)=\widetilde c_{bt}^{\mathrm B,\mathrm{init}}+\sum_a\delta_{ar}\widetilde c_{bta}^{\mathrm B}\), with the extension costs in Eq.~\eqref{eq:BlockExtensionCost}. By the extension rules in Sections~\ref{subsubsec:ExactDemandPricing} and~\ref{subsubsec:ExactBlockPricing}, each logical label represents an elementary prefix from \(o_k\) that uses compatible blocks, each physical label represents an elementary prefix from \(i_b\) of the recorded length within the length budget, and the cost of a label equals the initial cost plus the extension costs along its prefix. Physical labels are generated only for blocks with \(\overline f_b\ge\underline f_b\) (Section~\ref{subsubsec:ExactBlockPricing}), so a complete physical path within the length budget belongs to \(\mathcal R_b\). A \textit{completion} of a label is a possibly empty sequence of blocks or arcs whose concatenation with the prefix is a path in \(\mathcal P_k\) or \(\mathcal R_b\), respectively. A label at the destination has only the empty completion, since an elementary path ends there. If \(\Delta(Q)\) denotes the sum of the extension costs along a completion \(Q\), the decompositions above show that the completed path has JRC \(c+\Delta(Q)\), where \(c\) is the cost of the label. Let \(c^\star(L)\) be the minimum of \(c+\Delta(Q)\) over the completions \(Q\) of label \(L\). Here and below, a minimum over the empty set equals \(+\infty\).

Now let \(L_1\preceq L_2\) be two logical labels at yard \(i\), as in Proposition~\ref{prop:LogicalDominance}, and let \(Q\) be a completion of \(L_2\). Because the path completed from \(L_2\) is elementary, each yard of \(Q\) other than \(i\) lies outside \(S_2\), and therefore outside \(S_1\subseteq S_2\). Appending \(Q\) to the prefix of \(L_1\) thus yields an elementary path from \(o_k\) to \(d_k\) with compatible blocks, so \(Q\) is also a completion of \(L_1\). The extension costs of \(Q\) depend only on \(k\) and the blocks of \(Q\), so the two completed paths have JRCs \(c_1+\Delta(Q)\le c_2+\Delta(Q)\), whatever the sign of \(\Delta(Q)\). Hence \(c^\star(L_1)\le c^\star(L_2)\), which is the sense in which \(L_1\) dominates \(L_2\). For two physical labels \(L_1\preceq L_2\) at node \(n\), as in Proposition~\ref{prop:PhysicalDominance}, the same argument applies with nodes in place of yards and with extension costs that depend only on \(b\), \(t\), and the arcs of \(Q\). In addition, the path completed from \(L_1\) has length \(d_1+D(Q)\le d_2+D(Q)\le\overline D_b\), where \(D(Q)\) is the length of \(Q\). Hence \(Q\) is a completion of \(L_1\) within the length budget, and again \(c^\star(L_1)\le c^\star(L_2)\).

It remains to show that discarding dominated labels preserves the minimum JRC computed by the labeling algorithm. A label at the destination is stored and not extended. Every other label is eventually extended unless it is discarded, and a label \(L\) is discarded only when a stored label \(L_1\preceq L\) exists that has not been discarded, whether or not \(L_1\) has been extended. The algorithm stops because every extension enlarges the visited set, which is contained in \(\mathcal Y\) or \(\mathcal N\). Let \(\mathcal D\) be the set of labels at the destination that are never discarded. Each label in \(\mathcal D\) represents a path in \(\mathcal P_k\) or \(\mathcal R_b\), so its cost is at least the minimum JRC. We show, by induction on the number of blocks or arcs of \(Q\), that for every generated label \(L\) and every completion \(Q\) of \(L\), some label in \(\mathcal D\) has cost at most \(c+\Delta(Q)\), where \(c\) is the cost of \(L\). Applied to \(L_0\), whose completions are all paths in \(\mathcal P_k\) or \(\mathcal R_b\), this claim shows that the minimum cost over \(\mathcal D\) equals the minimum JRC.

Assume the claim for all completions shorter than \(Q\). Suppose first that \(L\) is never discarded. If \(Q\) is empty, then \(L\) lies at the destination and belongs to \(\mathcal D\). Otherwise, \(L\) is not at the destination and is therefore extended. Its extension \(L'\) along the first block or arc of \(Q\) is generated, because the head of that block or arc lies outside the visited set of \(L\) and, for a physical label, the extended prefix stays within the length budget, since the remaining arcs of \(Q\) have positive lengths. The rest \(Q'\) of \(Q\) is a completion of \(L'\), the cost \(c'\) of \(L'\) satisfies \(c'+\Delta(Q')=c+\Delta(Q)\), and \(Q'\) is shorter than \(Q\), so the induction hypothesis gives the claim. Suppose next that \(L\) is discarded. At that time, some label \(L_1\preceq L\) exists that has not been discarded. If \(L_1\) is discarded later, some label \(L_2'\preceq L_1\) exists that has not been discarded, and so on. The discard times increase strictly along this chain, so the chain ends at a label \(L^\ast\) that is never discarded. By transitivity of \(\preceq\), \(L^\ast\preceq L\), so the preceding paragraph shows that \(Q\) is a completion of \(L^\ast\) with cost \(c^\ast+\Delta(Q)\le c+\Delta(Q)\), where \(c^\ast\) is the cost of \(L^\ast\). The first case, applied to \(L^\ast\) and \(Q\), gives the claim. \Halmos
\endproof

\subsection{Proof of Proposition~\ref{prop:CompletionBounds}}
\label{ec:CompletionBounds}

\proof{Proof.}
Fix a block \(b\) with \(\overline f_b\ge\underline f_b\). Every completion \(Q\) of a physical label \(L=(n,S,d,c)\) is a directed path in \(G^{\mathrm P}\) from \(n\) to \(j_b\), possibly empty, so its length satisfies \(D(Q)\ge\underline D_b(n)\), the shortest length from \(n\) to \(j_b\) without the visited-node restriction. If no directed path from \(n\) to \(j_b\) exists, then \(L\) has no completion and \(\underline D_b(n)=+\infty\). The distances \(\underline D_b(n)\) are therefore lower bounds on the lengths of all completions. Dijkstra's algorithm from \(j_b\) on the reversed graph computes \(\underline D_b(n)\) for every node \(n\), because reversing the arcs maps the directed paths from \(n\) to \(j_b\) onto the directed paths from \(j_b\) to \(n\) with the same lengths, and all lengths \(D_a\) are positive. These distances depend only on \(G^{\mathrm P}\), the arc lengths, and \(j_b\), so they can be reused for every load and every set of dual prices.

If \(d+\underline D_b(n)>\overline D_b\), every completion \(Q\) of \(L\) would give a complete path of length \(d+D(Q)\ge d+\underline D_b(n)>\overline D_b\), which lies outside \(\mathcal R_b\), so \(L\) has no completion. Similarly, the extension of \(L\) along an arc \(a\in\mathcal A^+(n)\) is a label at \(j_a\) of length \(d+D_a\), and each of its completions has length at least \(\underline D_b(j_a)\). If \(d+D_a+\underline D_b(j_a)>\overline D_b\), this extension therefore has no completion. Both tests thus remove only labels without completions. In the induction at the end of Section~\ref{ec:LogicalDominance}, the claim is vacuous for such labels, and the extension \(L'\) used there is never rejected, because it has the completion \(Q'\). The induction therefore holds with the tests, and the labeling algorithm still returns the minimum physical path JRC. \Halmos
\endproof

\subsection{Proof of Proposition~\ref{prop:ExactPricing}}
\label{ec:ExactPricing}

\proof{Proof.}
We use LP duality for the DW LP relaxation. We call its constraints other than the path volume bounds in Eqs.~\eqref{eq:DWDemandVolume} and~\eqref{eq:DWBlockPathVolume} the \textit{master rows}, i.e., Eqs.~\eqref{eq:DWBlockBalance}--\eqref{eq:DWTrackCapacity} and~\eqref{eq:DWDemandConvexity}--\eqref{eq:DWServiceActivation}. The RMP contains the master rows, restricted to its paths, and the volume bounds of its paths. When rows of Eq.~\eqref{eq:DWServiceActivation} are separated as in Section~\ref{subsec:ExperimentProtocol}, the argument below applies once no omitted row is violated, because the omitted rows can then be added with zero dual variables without changing the optimal primal and dual solutions or their objective values. In a minimization LP, the dual variable of a \(\le\) constraint is nonpositive, that of an equality is free, and the reduced cost of a variable equals its objective coefficient minus the sum of its constraint coefficients weighted by the dual variables. The dual objective is the sum of the right-hand sides weighted by the dual variables, plus the objective constant \(\sum_kq_kc_k^{\mathrm U}\).

\textit{Step 1.} We derive the reduced costs. The selection variable \(x_{kp}^{\mathrm D}\) has objective coefficient zero, coefficient one in Eq.~\eqref{eq:DWDemandConvexity} for \(k\), coefficient \(\iota_{bp}\) in Eq.~\eqref{eq:DWServiceActivation} for each pair \((k,b)\), and coefficient \(-q_k\) in its row \(v_{kp}^{\mathrm D}-q_kx_{kp}^{\mathrm D}\le0\) of Eq.~\eqref{eq:DWDemandVolume}, whose dual variable is \(\lambda_{kp}\le0\). The volume variable \(v_{kp}^{\mathrm D}\) has objective coefficient \(c_{kp}^{\mathrm D}\), coefficient \(-\iota_{bp}\) in Eq.~\eqref{eq:DWBlockBalance} for each block \(b\), coefficient \(\eta_{ip}\) in Eq.~\eqref{eq:DWHandlingCapacity} for each yard \(i\), and coefficient one in the same bound row. Subtracting the dual-weighted coefficients from the objective coefficients gives the reduced costs in Eq.~\eqref{eq:LogicalReducedCost}, with the coefficients in Eq.~\eqref{eq:LogicalAdjustedCost}. For a physical path \(r\in\mathcal R_b\), we write the bounds in Eq.~\eqref{eq:DWBlockPathVolume} as \(\underline f_bx_{br}^{\mathrm B}-v_{br}^{\mathrm B}\le0\) and \(v_{br}^{\mathrm B}-\overline f_bx_{br}^{\mathrm B}\le0\), with dual variables \(\underline\lambda_{br},\overline\lambda_{br}\le0\). The variable \(x_{br}^{\mathrm B}\) has objective coefficient \(c_b^{\mathrm O}\), coefficient one in Eq.~\eqref{eq:DWBlockConvexity} for \(b\), coefficient \(\kappa_b\) in Eq.~\eqref{eq:DWTrackCapacity} for \(i_b\), coefficient \(-1\) in Eq.~\eqref{eq:DWServiceActivation} for each \(k\in\mathcal K_b\), and coefficients \(\underline f_b\) and \(-\overline f_b\) in its two bound rows. The variable \(v_{br}^{\mathrm B}\) has objective coefficient \(c_{br}^{\mathrm B}\), coefficient one in Eq.~\eqref{eq:DWBlockBalance} for \(b\), coefficient \(\sum_{a\in\mathcal A_e}\delta_{ar}\) in Eq.~\eqref{eq:DWLinkCapacity} for each link \(e\), and coefficients \(-1\) and \(1\) in the two bound rows. Since each arc belongs to exactly one link, \(\sum_e\mu_e\sum_{a\in\mathcal A_e}\delta_{ar}=\sum_a\mu_{e_a}\delta_{ar}\), and the reduced costs are
\[
\operatorname{rc}(x_{br}^{\mathrm B})=\bar c_{br}^{\mathrm B,x}-\underline f_b\underline\lambda_{br}+\overline f_b\overline\lambda_{br},
\qquad
\operatorname{rc}(v_{br}^{\mathrm B})=\bar c_{br}^{\mathrm B,v}+\underline\lambda_{br}-\overline\lambda_{br},
\]
with the coefficients in Eq.~\eqref{eq:PhysicalReducedCost}.

\textit{Step 2.} We fix the dual variables of the master rows and characterize when the bound rows of a single path admit dual variables that make the reduced costs of both of its variables nonnegative. For a logical path, such a \(\lambda_{kp}\le0\) exists if and only if \(\bar c_{kp}^{\mathrm D,x}\ge0\) and \(\widetilde c_{kp}^{\mathrm D}\ge0\). If such a \(\lambda_{kp}\) exists, then \(\bar c_{kp}^{\mathrm D,x}=\operatorname{rc}(x_{kp}^{\mathrm D})-q_k\lambda_{kp}\ge0\) and \(\widetilde c_{kp}^{\mathrm D}=\operatorname{rc}(x_{kp}^{\mathrm D})+q_k\operatorname{rc}(v_{kp}^{\mathrm D})\ge0\). Conversely, \(\lambda_{kp}=-\bar c_{kp}^{\mathrm D,x}/q_k\le0\) gives \(\operatorname{rc}(x_{kp}^{\mathrm D})=0\) and \(\operatorname{rc}(v_{kp}^{\mathrm D})=\widetilde c_{kp}^{\mathrm D}/q_k\ge0\). Moreover, \(\bar c_{kp}^{\mathrm D,x}=-\pi_k-\sum_b\gamma_{kb}\iota_{bp}\ge0\) holds for every logical path because \(\pi_k,\gamma_{kb}\le0\), so the condition reduces to \(\widetilde c_{kp}^{\mathrm D}\ge0\). For a physical path, such \(\underline\lambda_{br},\overline\lambda_{br}\le0\) exist if and only if \(\widetilde c_{br}^{\mathrm B}(\underline f_b)\ge0\) and \(\widetilde c_{br}^{\mathrm B}(\overline f_b)\ge0\). By Step 1, for every \(t\),
\[
\operatorname{rc}(x_{br}^{\mathrm B})+t\operatorname{rc}(v_{br}^{\mathrm B})=\widetilde c_{br}^{\mathrm B}(t)+(t-\underline f_b)\underline\lambda_{br}+(\overline f_b-t)\overline\lambda_{br}.
\]
If both reduced costs are nonnegative, then for \(t\in\{\underline f_b,\overline f_b\}\) the left side is nonnegative and the last two terms are nonpositive, so \(\widetilde c_{br}^{\mathrm B}(t)\ge0\). Conversely, if \(\bar c_{br}^{\mathrm B,v}\ge0\), the choice \(\underline\lambda_{br}=-\bar c_{br}^{\mathrm B,v}\) and \(\overline\lambda_{br}=0\) gives \(\operatorname{rc}(v_{br}^{\mathrm B})=0\) and \(\operatorname{rc}(x_{br}^{\mathrm B})=\widetilde c_{br}^{\mathrm B}(\underline f_b)\ge0\). If \(\bar c_{br}^{\mathrm B,v}<0\), the choice \(\underline\lambda_{br}=0\) and \(\overline\lambda_{br}=\bar c_{br}^{\mathrm B,v}\) gives \(\operatorname{rc}(v_{br}^{\mathrm B})=0\) and \(\operatorname{rc}(x_{br}^{\mathrm B})=\widetilde c_{br}^{\mathrm B}(\overline f_b)\ge0\).

\textit{Step 3.} We prove the proposition. Take an optimal primal solution of the RMP and the optimal RMP dual solution whose master-row dual variables define the JRCs. By hypothesis and Step 2, every path in \(\mathcal P_k\) or \(\mathcal R_b\), whether or not it belongs to the RMP, admits bound dual variables that make the reduced costs of its two variables nonnegative under these master-row dual variables. Together, the master-row dual variables and these bound dual variables form a feasible dual solution of the DW LP relaxation. Its objective value equals that of the optimal RMP dual solution, because the master rows coincide and the bound rows have zero right-hand sides, and therefore equals the optimal RMP value by LP duality. Extending the RMP primal solution with zero values for all paths outside the RMP gives a feasible solution of the DW LP relaxation with the same objective value, because the added variables contribute nothing to the objective or to the master rows and satisfy their bound rows as \(0\le0\). By weak duality, the extended solution is optimal for the DW LP relaxation, so the RMP value equals \(z_{\mathrm{dw}}^{\mathrm{lp}}\). Since the DW LP relaxation relaxes the binary domains of DW, \(z_{\mathrm{dw}}^{\mathrm{lp}}\le z_{\mathrm{dw}}^*\). \Halmos
\endproof

Algorithm~\ref{alg:PricingFramework} terminates, and by Proposition~\ref{prop:ExactPricing} its final RMP solution is optimal for the DW LP relaxation. Every RMP is feasible, since setting all path variables to zero satisfies its rows, the capacities being nonnegative, and Eqs.~\eqref{eq:DWDemandVolume}--\eqref{eq:DWBlockConvexity} bound all variables, so every RMP has an optimal solution. For a path in the RMP, the optimal RMP dual solution provides bound dual variables that make the reduced costs of its two variables nonnegative, so Step 2 shows that its JRC is nonnegative, at both load bounds for a physical path. Exact pricing returns only paths with negative JRC, so every iteration that does not stop adds paths outside the RMP. Since every \(\mathcal P_k\) and \(\mathcal R_b\) is finite, exact pricing eventually returns no path. By Propositions~\ref{prop:LogicalDominance}--\ref{prop:CompletionBounds} and the single search for both load bounds in Section~\ref{subsubsec:ExactBlockPricing}, exact pricing minimizes the logical JRC over every \(\mathcal P_k\) and the physical JRC over every \(\mathcal R_b\) at both load bounds, so an empty return means that every JRC is nonnegative, and Proposition~\ref{prop:ExactPricing} applies. If rows of Eq.~\eqref{eq:DWServiceActivation} are separated, each RMP solve adds violated rows until none remains, which takes finitely many solves because these rows are finite in number.
\par

\section{Heuristic Details}
\label{sec:HeuristicCompanion}

\subsection{\textit{Construction} Details}
\label{ec:IncrementalEvaluation}

The initial restricted block domain in Section~\ref{subsec:HeuristicCandidates} uses compatible direct connections and at most eight hub proposals per demand. Hubs are ranked by compatible volume support divided by summed path length budgets. The full physical graph, original capacities, demands, and individual compatibility are retained.

The physical candidates \(\mathcal R_b^{\mathrm C}\subseteq\mathcal R_b\) of a block \(b\) are generated by shortest-path searches when \(b\) is first considered. Let \(q_b^{\max}=\max\{\underline f_b,\max_{k\in\mathcal K_b}q_k\}\). For each load level \(t\in\{\underline f_b,\min\{\overline f_b,q_b^{\max}\},\min\{\overline f_b,2q_b^{\max}\}\}\), in increasing order, we remove the arcs of the links with \(C_e<t\) and run Dijkstra's algorithm from \(i_b\) to \(j_b\) with each of the arc weights \(D_a\), \(c_{ba}^{\mathrm T}\), \(c_{ba}^{\mathrm T}+D_a\), and \(c_{ba}^{\mathrm T}+4D_a\). If no path exists or the minimum-length path exceeds \(\overline D_b\), no path within the length budget exists at this or any higher level, and the searches for \(b\) stop. If the minimum-cost path exceeds \(\overline D_b\), we also examine the two cheapest elementary paths with respect to \(c_{ba}^{\mathrm T}\) (Table~\ref{tab:HeuristicDefaults}). When the unit transportation costs of \(b\) are proportional to the arc lengths, all four weights yield the same path, and only the length search is run. A path joins \(\mathcal R_b^{\mathrm C}\) if it is new, satisfies the length budget, and can carry the minimum block load, i.e., \(\overline f_{br}\ge\underline f_b\) in Eq.~\eqref{eq:HeuristicRouteVolume}.

When \textsc{Serve} evaluates a candidate assignment (Section~\ref{subsubsec:ConstructionJointEvaluation}), the open blocks keep their incumbent physical paths, and the newly opened blocks receive paths from their candidates one at a time, in increasing block index. \textsc{Serve} compares two selection rules, both of which account for the paths selected earlier. The first selects the cheapest candidate that keeps every link within capacity at the smallest value of \(\theta\) allowed by the block load bounds and handling capacities. The second selects the candidate that keeps the largest value of \(\theta\) within the link capacities, breaking ties by unit cost. For each rule, \textsc{Serve} computes \([\theta_-,\theta_+]\) with the selected paths and applies Eq.~\eqref{eq:ConstructionServiceVolumes}, and it keeps the cheaper of the two resulting plans.
\par

\subsection{Search Schedule}
\label{ec:HeuristicSchedule}

After \textit{construction}, DG-CGM repeats search cycles until the time limit, using the controls in Table~\ref{tab:HeuristicDefaults}. Each cycle starts with a local search stage. It then performs pricing rounds, each consisting of one \textit{path generation} invocation followed by a local search stage, until two consecutive rounds each fail to lower the incumbent cost by more than \(\varepsilon_{\mathrm H}\max\{1,|z|\}\), where \(z\) is the incumbent cost at the start of the round. \textit{Joint repair} then runs once, and the next cycle starts.

A local search stage alternates \textit{destroy} and \textit{local repair} for at most \(|\mathcal K|\) calls within 30 seconds, and each local MIP has a time limit of one second. A call is improving if it lowers the incumbent cost by more than \(\varepsilon_{\mathrm H}\max\{1,|z|\}\), where \(z\) is the incumbent cost before the call. Every strict feasible cost decrease is retained. In each stage, the group size limit \(N_{\mathrm D}\) of Section~\ref{subsec:HeuristicDestroyRepair} starts at eight, rises to twelve after two consecutive nonimproving calls, and returns to eight after an improving call. The stage ends after three consecutive nonimproving calls.

A \textit{path generation} invocation is skipped when the incumbent and the candidate path pools coincide with those at the start of the most recent invocation that ended with an optimal LP and a complete pricing sweep, since it would return the same paths. Its round then consists of the local search stage alone. The 30-second allowance of a \textit{joint repair} invocation starts before the candidate searches of Section~\ref{subsec:HeuristicFullMIP}, so the joint MIP receives the remaining time. The search for alternative physical paths visits the open blocks in turn and stops after 0.1 seconds in total. For each retained physical path of a block and each node of this path, it keeps the prefix up to the node, removes the prefix nodes and the arcs by which retained paths with the same prefix continue, and computes a minimum-cost suffix with arc costs \(c_{ba}^{\mathrm T}\) over the links with \(C_e\ge\underline f_b\). If the resulting path violates the length budget, revisits a node, or cannot carry \(\underline f_b\), the suffix search is repeated with arc costs \(c_{ba}^{\mathrm T}+D_a\), \(c_{ba}^{\mathrm T}+4D_a\), and \(D_a\). The two cheapest new paths of each block join its pool. Section~\ref{ec:PartnerChains} describes the search for partner chains. The joint MIP is skipped when it coincides with the most recent joint MIP whose search was completed, that is, which was solved to optimality or shown to contain no plan cheaper than the incumbent.

The search stops at the time limit or after a cycle in which every \textit{path generation} invocation and the joint MIP are skipped and the candidate searches before \textit{joint repair} add no path.
\par

\subsection{Minimum-Load Partner Chains}
\label{ec:PartnerChains}

Partner chains are generated before each joint MIP from the incumbent and the shared logical pools \(\widetilde{\mathcal P}_j\). Let \(\mathcal B^{\mathrm O}\) be the set of open blocks, \(c_b^{\min}=\min_{r\in\widetilde{\mathcal R}_b}c_{br}^{\mathrm B}\), and \(\bar c_j\) the incumbent unit cost of demand \(j\) in Eq.~\eqref{eq:ConstructionUnitSaving}. For a path \(p\in\widetilde{\mathcal P}_j\), let \(\mathcal D_{jp}\) contain its blocks \(b\notin\mathcal B^{\mathrm O}\) with \(q_j<\underline f_b\), which demand \(j\) cannot open alone, and let
\[
\begin{aligned}
\Delta_{jp}={}&\bar v_j\bar c_j+(q_j-\bar v_j)c_j^{\mathrm U}+\sum_{b\in\mathcal B_j^{\mathrm{excl}}\setminus p}c_b^{\mathrm O}\\
&-q_j\sum_{b\in p}\bigl(c_b^{\min}+\mathbf1\{i_b\ne o_j\}c_{i_b}^{\mathrm H}\bigr)-\sum_{b\in p\setminus\mathcal B^{\mathrm O}}c_b^{\mathrm O}
\end{aligned}
\]
be the estimated saving from moving all of \(q_j\) to \(p\), where \(\mathcal B_j^{\mathrm{excl}}\) contains the blocks of \(\bar p_j\) that carry no other served demand. For a demand \(k\) and a block \(b\in\mathcal B_k\) with \(\widetilde{\mathcal R}_b\ne\emptyset\), the reward
\begin{equation}
R_{kb}=\max\Bigl(\{0\}\cup\Bigl\{\frac{\Delta_{jp}}{|\mathcal D_{jp}|}:j\ne k,\ p\in\widetilde{\mathcal P}_j,\ \Delta_{jp}>0,\ b\in\mathcal D_{jp},\ q_j+q_k\ge\underline f_b\Bigr\}\Bigr)
\label{eq:SupportReward}
\end{equation}
assigns to \(b\) the largest share of such a saving among the paths of other demands that need \(b\) and can open it together with \(k\). For each demand \(k\) with a positive reward in Eq.~\eqref{eq:SupportReward} for some \(b\), a labeling search finds the cheapest elementary paths from \(o_k\) to \(d_k\) over the compatible blocks with retained physical paths, where block \(b\) costs
\[
q_k\bigl(c_b^{\min}+\mathbf1\{i_b\ne o_k\}c_{i_b}^{\mathrm H}\bigr)-R_{kb}+\mathbf1\{b\notin\mathcal B^{\mathrm O},\,R_{kb}=0\}c_b^{\mathrm O}.
\]
An unopened block is charged its opening cost only when it carries no reward, because a positive reward is computed after the opening costs in \(\Delta_{jp}\) have been deducted. The three cheapest paths that are not yet in \(\widetilde{\mathcal P}_k\) join the pool, and the searches for all demands stop after 0.1 seconds in total. The rewards only guide this search, and the joint MIP evaluates the resulting paths with the original objective and constraints.
\par

\subsection{Arc Options in \textit{Joint Repair}}
\label{subsec:HeuristicFixedResources}

The arc option of Section~\ref{subsec:HeuristicFullMIP} gives each released demand \(k\) a copy \((x_k^{\mathrm{D,arc}},x_{kb}^{\mathrm{D,arc}},g_{kb}^{\mathrm{arc}},v_k^{\mathrm{arc}},\ell_{ki}^{\mathrm{D,arc}})\) of its logical variables in ARC. The copy satisfies Eqs.~\eqref{eq:LogicalBalance}--\eqref{eq:ServedVolume} over the compatible blocks with retained physical paths, and \(x_{kb}^{\mathrm{D,arc}}\le\sum_{r\in\widetilde{\mathcal R}_b}x_{br}^{\mathrm B}\) replaces Eq.~\eqref{eq:LogicalActivation}. The constraints \(\sum_{p\in\widetilde{\mathcal P}_k}x_{kp}^{\mathrm D}+x_k^{\mathrm{D,arc}}\le1\) for each released demand and \(\sum_{k\in\mathcal K_{\mathrm{free}}}x_k^{\mathrm{D,arc}}\le1\) allow each demand at most one of its options and at most one demand the arc option. The copy enters the other constraints and the objective as a logical path does. The volume \(g_{kb}^{\mathrm{arc}}\) is added to the volumes of the released demands on block \(b\) in the block balance of Eq.~\eqref{eq:LocalRepairCoupling}, the reclassified volume \(\sum_{b\in\mathcal B_k^+(i)}g_{kb}^{\mathrm{arc}}-\mathbf1\{i=o_k\}v_k^{\mathrm{arc}}\) is added at each yard \(i\) to the handling capacity and the handling cost, and the served volume of \(k\) becomes \(\sum_{p\in\widetilde{\mathcal P}_k}v_{kp}^{\mathrm D}+v_k^{\mathrm{arc}}\) in the unserved demand cost. The size estimate \(s_{\mathrm{arc}}\) is the number of variables of the augmented MIP plus an upper bound on its number of constraints, computed before the model is built.
\par

\subsection{Physical Path Load Limits}
\label{subsec:HeuristicJRCProof}

An elementary physical path visits no node twice, so it traverses each link at most once, that is, \(\sum_{a\in\mathcal A_e}\delta_{ar}\in\{0,1\}\). A path \(r\) selected for block \(b\) therefore carries the full block load on every link that it traverses, and the restricted MIPs of \textit{local repair} and \textit{joint repair} use the load limit
\begin{equation}
\overline f_{br}=\min\Bigl\{\overline f_b,\min_{e\in\mathcal E:\,\sum_{a\in\mathcal A_e}\delta_{ar}=1}C_e\Bigr\}
\label{eq:HeuristicRouteVolume}
\end{equation}
in place of \(\overline f_b\) in Eq.~\eqref{eq:DWBlockPathVolume}, while the physical candidate sets exclude paths with \(\overline f_{br}<\underline f_b\). These changes remove no integer feasible solution of a restricted MIP. If such a solution selects path \(r\) for block \(b\), then for every link \(e\) traversed by \(r\), the link row of Eq.~\eqref{eq:LocalRepairCoupling}, with the fixed loads moved back to its left side, contains the term \(v_{br}^{\mathrm B}\) and further nonnegative terms, so \(v_{br}^{\mathrm B}\le C_e\). Together with \(v_{br}^{\mathrm B}\le\overline f_b\), this gives \(v_{br}^{\mathrm B}\le\overline f_{br}\). If the solution does not select \(r\), then \(v_{br}^{\mathrm B}=0\). A path with \(\overline f_{br}<\underline f_b\) can therefore never be selected, because it would require \(\underline f_b\le v_{br}^{\mathrm B}\le\overline f_{br}\). Since \(\overline f_{br}\le\overline f_b\), the load limits can only shrink the LP relaxations of the restricted MIPs, and \(\overline f_{br}\) can only decrease when link capacities are reduced.

In each \textit{local repair} call, the candidate set \(\widetilde{\mathcal R}_b\) of a block \(b\in\mathcal B_{\mathrm{free}}\) contains up to three of its cheapest retained physical paths, a retained path with the largest \(\overline f_{br}\), and its incumbent path.
\par

\subsection{Default Search Controls}
\label{subsec:HeuristicDefaults}

Table~\ref{tab:HeuristicDefaults} lists the default controls of DG-CGM, which all experiments use unless stated otherwise.

\begin{table}[t]
\centering
\fontsize{10}{12}\selectfont\setlength{\tabcolsep}{4pt}\renewcommand{\arraystretch}{1.05}
\caption{Default search controls.}
\label{tab:HeuristicDefaults}
\begin{tabular}{@{}lr@{}}
\toprule
Control & Default \\
\midrule
Related-demand score weights \((\omega_0,\omega_1,\omega_2)\) & \((10^4,10^3,500)\) \\
Initial demand-group size limit & \(8\) \\
Consecutive nonimprovements triggering expansion & \(2\) \\
Expanded demand-group size limit & \(12\) \\
Stage-control relative-improvement threshold \(\varepsilon_{\mathrm H}\) & \(10^{-4}\) \\
Consecutive nonimprovements ending a local search stage & \(3\) \\
Maximum number of local MIP calls per local search stage & \(|\mathcal K|\) \\
Time allowance per local search stage & \(30\,\mathrm{s}\) \\
Nonimproving completed pricing rounds triggering \textit{joint repair} & \(2\) \\
Time limit per local MIP & \(1\,\mathrm{s}\) \\
Cheapest retained physical paths per released block & \(3\) \\
Alternative physical paths per open block before \textit{joint repair} & \(2\) \\
Time limit for this physical path search & \(0.1\,\mathrm{s}\) \\
Partner chains per demand before \textit{joint repair} & \(3\) \\
Time limit for the partner chain search & \(0.1\,\mathrm{s}\) \\
Regenerated lowest-EUC logical paths per demand in each repair & \(16\) \\
\textit{Construction} stopping tolerance \(\varepsilon_{\mathrm C}\) & \(10^{-6}\) \\
Hub yards for two-hub groups \(N_{\mathrm{hub}}\) & \(16\) \\
Candidate blocks above which \textit{construction} restricts the block set & \(10{,}000\) \\
New blocks permitted in an expanded logical chain & \(2\) \\
Initial logical-path limit per demand \(N_{\mathrm{LP}}\) & \(4\) \\
Return cap per demand \(N_{\mathrm{add}}\) & \(12\) \\
Mixed-selection JRC quota \(N_{\mathrm{JRC}}\) & \(8\) \\
Mixed-selection diversity quota \(N_{\mathrm{div}}\) & \(2\) \\
Mixed-selection EUC quota \(N_{\mathrm{EUC}}\) & \(2\) \\
Successor-yard neighborhood size \(N_{\mathrm{nbr}}\) & \(32\) \\
Soft time allowance per RMP iteration \(\tau_{\mathrm{PG}}\) & \(10\,\mathrm{s}\) \\
New logical paths per LP batch \(N_{\mathrm{batch}}\) & \(1{,}000\) \\
Maximum pricing batches per invocation \(I_{\mathrm{PG}}\) & \(10\) \\
LP plateau window \(W_{\mathrm{PG}}\) & \(3\) \\
Relative LP plateau tolerance \(\varepsilon_{\mathrm{PG}}\) & \(10^{-5}\) \\
Pricing-batch share of remaining overall time & \(0.2\) \\
Maximum demands released in \textit{joint repair} \(N_{\mathrm J}\) & \(250\) \\
Joint MIP size limit for arc options \(N_{\mathrm{arc}}\) & \(10^5\) \\
Cost-ranked shortest paths examined after an initial length violation & \(2\) \\
Time allowance per \textit{joint repair} invocation & \(30\,\mathrm{s}\) \\
\bottomrule
\end{tabular}
\end{table}

\subsection{Sensitivity to Search Controls}
\label{subsec:ParameterSensitivity}
Table~\ref{tab:ParameterSensitivity} reports the one-at-a-time variants of Section~\ref{subsubsec:SmallRobustness}. Each variant changes one control of Table~\ref{tab:HeuristicDefaults} and keeps all others, the 60-second time limit, one solver thread, and seed~1 on the 150 generated instances. Diff.\ is the mean paired relative cost difference from the default configuration, W/T/L counts the instances on which the variant is cheaper, equal, or more expensive, and \(p\) is the two-sided Wilcoxon signed-rank \(p\)-value, reported when at least six instances differ.

\begin{table}[htbp]
\centering
\caption{One-at-a-time sensitivity of DG-CGM to its search controls.}
\label{tab:ParameterSensitivity}
\fontsize{10}{12}\selectfont\setlength{\tabcolsep}{3pt}\renewcommand{\arraystretch}{1.05}
\begin{tabular*}{\linewidth}{@{\extracolsep{\fill}}lccrrrrrcr@{}}
\toprule
Control & Default & Value & \shortstack[c]{Mean\\gap/\%} & \shortstack[c]{Max.\\gap/\%} & \#Opt. & Time/s & Diff./\% & W/T/L & \(p\) \\
\midrule
\(N_{\mathrm{LP}}\) & 4 & 2 & 0.10 & 5.8 & 89 & 49.3 & \(-0.05\) & 5/142/3 & 0.148 \\
\(N_{\mathrm{LP}}\) & 4 & 8 & 0.10 & 5.8 & 89 & 51.1 & \(-0.04\) & 4/144/2 & 0.219 \\
\((N_{\mathrm{add}},N_{\mathrm{JRC}},N_{\mathrm{div}},N_{\mathrm{EUC}})\) & (12,8,2,2) & (6,4,1,1) & 0.06 & 2.8 & 89 & 50.7 & \(-0.08\) & 6/142/2 & 0.078 \\
\((N_{\mathrm{add}},N_{\mathrm{JRC}},N_{\mathrm{div}},N_{\mathrm{EUC}})\) & (12,8,2,2) & (24,16,4,4) & 0.09 & 4.5 & 89 & 50.6 & \(-0.06\) & 5/144/1 & 0.094 \\
\(N_{\mathrm{nbr}}\) & 32 & 16 & 0.12 & 5.8 & 89 & 50.7 & \(-0.03\) & 4/143/3 & 0.297 \\
\(N_{\mathrm{nbr}}\) & 32 & 64 & 0.15 & 5.8 & 89 & 50.8 & 0.00 & 0/150/0 & -- \\
\(I_{\mathrm{PG}}\) & 10 & 5 & 0.09 & 5.8 & 91 & 51.4 & \(-0.05\) & 6/142/2 & 0.039 \\
\(I_{\mathrm{PG}}\) & 10 & 20 & 0.17 & 5.8 & 88 & 50.7 & 0.03 & 1/142/7 & 0.016 \\
\(\tau_{\mathrm{PG}}\)/s & 10 & 5 & 0.12 & 4.2 & 88 & 50.7 & \(-0.03\) & 1/148/1 & -- \\
\(\tau_{\mathrm{PG}}\)/s & 10 & 20 & 0.15 & 5.8 & 89 & 51.0 & 0.00 & 0/150/0 & -- \\
Group size limits & (8,12) & (4,8) & 0.14 & 5.8 & 89 & 50.1 & \(-0.01\) & 5/141/4 & 0.820 \\
Group size limits & (8,12) & (12,16) & 0.14 & 5.8 & 88 & 48.1 & \(-0.01\) & 1/147/2 & -- \\
\(\varepsilon_{\mathrm H}\) & \(10^{-4}\) & \(10^{-3}\) & 0.14 & 5.8 & 89 & 50.5 & \(-0.01\) & 1/148/1 & -- \\
\(\varepsilon_{\mathrm H}\) & \(10^{-4}\) & \(10^{-5}\) & 0.14 & 5.8 & 88 & 50.7 & \(-0.01\) & 1/147/2 & -- \\
\hline
Default & -- & -- & 0.15 & 5.8 & 89 & 50.8 & -- & -- & -- \\
\bottomrule
\end{tabular*}
\end{table}

\section{Details of Instance Generation}
\label{ec:instance-generator}

The generator takes a random seed and the target numbers of nodes and demands as input. It estimates all distributions from data set L1 (Section~\ref{subsec:ExperimentProtocol}) and keeps the network-wide parameters of L1, such as the cost coefficients, unchanged. Each L1 node is a yard, a gate, or a micro node, each yard has one of three levels, and every L1 demand connects yards of levels 1 and 2. We treat the L1 network as an undirected simple graph that keeps the shortest of any parallel links, so that the degree of a node is its number of neighbors. In its largest connected component, we replace every maximal path whose interior nodes are micro nodes of degree two by one link of the same total length, again keeping the shortest of any parallel links, and retain all other nodes. In the resulting simplified real network, we partition the nodes once, for all instances, into six regions by \(k\)-means \citep{macqueen1967methods} on their coordinates in miles under an equirectangular projection. Region \(g\) has \(\widehat n_g\) nodes, and \(\widehat m_{gh}\) and \(\widehat s_{gh}\) denote the mean and standard deviation, over these nodes, of the fraction of a node's neighbors that lie in region \(h\). The demands are apportioned by the largest-remainder method among commodities and then among pairs of origin and destination yard levels, in proportion to the numbers of L1 demand records, and each commodity and level pair with demands forms a demand category.

Each generated node joins region \(g\) with probability \(\widehat n_g/\sum_h\widehat n_h\), and a node \(i\) in region \(g\) draws a preference \(\varphi_{ih}\sim\operatorname{Normal}(\widehat m_{gh},\widehat s_{gh}^2)\) for links to each region \(h\), clipped to \([0,1]\) and normalized to sum to one. If \(n_g\) generated nodes lie in region \(g\), the nodes of this region in the simplified real network, sorted by degree, supply reference records at positions \(1+\lfloor(t-\tfrac12)\widehat n_g/n_g\rfloor\), \(t=1,\ldots,n_g\), which are assigned to the generated nodes of the region in random order. GenCAT \citep{maekawa2023gencat} then calibrates the preferences of each region to \(\widehat m_{gh}\) and adds links according to them, up to a target degree for each node equal to the degree of its reference record. Each link can be traversed in both directions.

The reference record of a generated node also supplies its node type and coordinates, with each coordinate perturbed by independent Gaussian noise whose standard deviation is 0.035 times that of the coordinate over the region in the simplified real network; \(\mathbf r_i\) denotes the result for node \(i\). The nodes whose reference records are yards form \(\mathcal Y\). To preserve the observed combinations of yard characteristics, we sample \(|\mathcal Y|\) complete records of yards of the simplified real network with replacement within each level. The numbers of records per level are apportioned in the same way by the level shares of that network, after which levels 1 and 2 are enlarged at the expense of level 3 until the \(n\) records of each of these levels give at least as many distinct ordered pairs, \(n(n-1)\), as any category with both endpoints at that level has demands. A minimum-cost one-to-one matching assigns these records to the generated yards, where assigning the record of yard \(j\) to generated yard \(i\) costs \(3|\log[(1+\deg(i))/(1+\deg(j))]|+\|\mathbf r_i-\mathbf r_j^{\mathrm{ref}}\|_2/1{,}500\), with degrees in the respective networks and the coordinates \(\mathbf r_j^{\mathrm{ref}}\) of yard \(j\) in miles. The assigned record provides the operator, level, classification track slots \(T_i\), handling capacity \(H_i\), and handling cost \(c_i^{\mathrm H}\) of yard \(i\).

The length of each arc \(a\) is first set to \(D_a=\max\{0.01,\|\mathbf r_{i_a}-\mathbf r_{j_a}\|_2\}\) miles. All arc lengths are then multiplied by \(\widehat D/D\), where \(\widehat D\) and \(D\) are the median shortest-path distances over unordered pairs of distinct yards in the simplified real network and in the generated network, so that the two medians coincide. As in L1, every link has capacity \(C_e=700{,}000\).

For each demand \(k\), we sample an L1 demand record uniformly from its category. This record supplies an initial volume \(q_k^{(0)}\), the shortest-path distance \(\widehat D_k\) between its origin and destination in the simplified real network, and whether these two yards have the same operator. Among the ordered pairs \((i,j)\) of distinct generated yards, with \(i\) at the origin level and \(j\) at the destination level of the category, that no earlier demand of the same commodity uses, we draw \((o_k,d_k)\) with probability proportional to \(\exp\{-3|\log(\max\{\operatorname{dist}(i,j),1\}/\max\{\widehat D_k,1\})|\}\), multiplied by \(\exp(-1.5)\) if the pair and the record differ in whether their two yards have the same operator. Here, \(\operatorname{dist}(i,j)\) is the shortest-path distance in miles in the generated network.

Within each commodity, the initial volumes divided by their mean define an outgoing and an incoming total at each yard. The conditional reconstruction method (CReMa) of \citet{parisi2020faster}, implemented in NEMtropy \citep{vallarano2021fast}, fits yard parameters \(\nu_i^+\) and \(\nu_i^-\) such that independent samples \(X_k\sim\operatorname{Exp}(\nu_{o_k}^++\nu_{d_k}^-)\) of the demands, where \(\operatorname{Exp}\) is parameterized by its rate, match these totals in expectation; every absolute deviation is at most \(10^{-5}\) times the largest total. Let \(\mathcal K'\) list the demands of the commodity in increasing order of \(X_k\), and let \(\widehat Q\) be the sample quantile function of the L1 volumes of the commodity, interpolated linearly between order statistics. The \(t\)-th demand \(k\) in \(\mathcal K'\) receives the requested volume \(q_k=\max\{1,\lfloor\widehat Q((t-\tfrac12)/|\mathcal K'|)\rceil\}\), where \(\lfloor\cdot\rceil\) rounds to the nearest integer, with ties to even. The final volumes thus follow the L1 volume distribution of the commodity, and the samples \(X_k\) determine only their order. An instance is valid if the network generated by GenCAT is connected, so that its largest connected component contains all generated nodes, and if it has at least eight yards, level 3 can supply the enlargement of levels 1 and 2, every category has at least as many ordered pairs of distinct yards at its levels as demands, and every fit meets the tolerance above.
\par

\end{document}